\documentclass[preprint,3p,authoryear]{elsarticle}
\usepackage{amsfonts,amssymb,amsmath,amsopn,float}
\usepackage{xifthen}
\usepackage{graphicx}
\usepackage{tikz}
\usepackage{pgfplots}
\usepgfplotslibrary{units}
\usepackage{epstopdf}
\usepackage{url}
\usepackage{subcaption}
\usepackage[pdftex,hypertexnames=true,plainpages=false,
naturalnames,colorlinks=true,linkcolor=blue,bookmarks]{hyperref}

\newcommand{\bbR}{\mathbb{R}}

\newcommand{\bbZ}{\mathbb{Z}}

\newcommand{\bzero}{\mathbf{0}}
\newcommand{\abs}[1]{\left\lvert #1 \right\rvert}

\newcommand{\normX}[2]{{\left\lVert #1 \right\rVert}_{#2}}

\newcommand{\bxi}{{\boldsymbol{\xi}}}

\newcommand{\dparder}[2]{\dfrac{\partial #1}{\partial #2}}

\newcommand{\dsecder}[2]{\dfrac{\partial^2 #1}{\partial #2^2}}

\newcommand{\Ltwo}[1]{%
\ifthenelse{\equal{#1}{}}{L^2}{L^2(#1)}%
}

\newcommand{\Ltwoz}[1]{%
\ifthenelse{\equal{#1}{}}{L^2_0}{L^2_0(#1)}%
}

\newcommand{\Ltwonorm}[2]{\left\lVert #1 \right\rVert_{\Ltwo{#2}}}

\newcommand{\Cone}[1]{%
\ifthenelse{\equal{#1}{}}{C^{1}}{C^{1}(#1)}%
}

\newcommand{\Conez}[1]{%
\ifthenelse{\equal{#1}{}}{C^{1}_{0}}{C^{1}_{0}(#1)}%
}

\newcommand{\Ctwo}[1]{%
\ifthenelse{\equal{#1}{}}{C^{2}}{C^2(#1)}%
}

\newcommand{\Ctwoz}[1]{%
\ifthenelse{\equal{#1}{}}{C^{2}_{0}}{C^{2}_{0}(#1)}%
}

\newcommand{\Cholder}[1]{%
\ifthenelse{\equal{#1}{}}{C^{0,\gamma}}{C^{0,\gamma}(#1)}%
}

\newcommand{\Cholderz}[1]{%
\ifthenelse{\equal{#1}{}}{C^{0,\gamma}_{0}}{C^{0,\gamma}_{0}(#1)}%
}

\newcommand{\Choldernorm}[2]{%
{\left\lVert #1 \right\rVert}_{\Cholder{#2}}%
}

\newcommand{\Ctwointime}[2]{C^{2}(#1;#2)}

\newcommand{\perdrd}{D;\bbR^{d}}

\newcommand{\bolds}[1]{\boldsymbol{#1}}
\newcommand{\ba}{\bolds{a}}
\newcommand{\bb}{\bolds{b}}

\newcommand{\be}{\bolds{e}}

\newcommand{\br}{\bolds{r}}

\newcommand{\bu}{\bolds{u}}
\newcommand{\bv}{\bolds{v}}

\newcommand{\bx}{\bolds{x}}
\newcommand{\by}{\bolds{y}}

\newcommand{\buhat}{\hat{\bu}}
\newcommand{\bubar}{\bar{\bu}}
\newcommand{\butilde}{\tilde{\bu}}
\newcommand{\bvhat}{\hat{\bv}}
\newcommand{\bvbar}{\bar{\bv}}
\newcommand{\bvtilde}{\tilde{\bv}}

\newcommand{\sautoref}[2]{\hyperref[#2]{#1 \ref*{#2}}}

\newtheorem{theorem}{Theorem}
\newtheorem{proposition}{Proposition}
\newtheorem{lemma}{Lemma}

\newdefinition{rmk}{Remark}
\newproof{proof}{Proof}

\ifpdf
  \DeclareGraphicsExtensions{.eps,.pdf,.png,.jpg}
\else
  \DeclareGraphicsExtensions{.eps}
\fi

\begin{document}

\begin{frontmatter}

\title{Numerical convergence of finite difference approximations for state based peridynamic fracture models\tnoteref{t1}}
\tnotetext[t1]{\textbf{Funding: }This material is based upon work supported by the U. S. Army Research Laboratory and the U. S. Army Research Office under contract/grant number W911NF1610456.}

\author[aa]{Prashant K. Jha\corref{cor1}\fnref{fn1}}
\ead{prashant.j16o@gmail.com}

\author[aa,bb]{Robert Lipton\fnref{fn2}}
\ead{lipton@lsu.edu}

\address[aa]{Department of Mathematics, Louisiana State University, Baton Rouge, LA}

\address[bb]{Center for Computation and Technology, Louisiana State University, Baton Rouge, LA}

\cortext[cor1]{Corresponding author (pjha4@lsu.edu)}

\fntext[fn1]{Orcid: https://orcid.org/0000-0003-2158-364X}
\fntext[fn2]{Orcid: https://orcid.org/0000-0002-1382-3204}

\begin{abstract}
In this work, we study the finite difference approximation for a class of nonlocal fracture models. 
The nonlocal model is initially elastic but beyond a critical strain the material softens with increasing strain. This model is formulated as a state-based peridynamic model using two potentials: one associated with hydrostatic strain and the other associated with tensile strain. We show that the dynamic evolution is well-posed in the space of H\"older continuous functions $C^{0,\gamma}$ with H\"older exponent $\gamma \in (0,1]$. Here the length scale of nonlocality is $\epsilon$, the size of time step is $\Delta t$ and the mesh size is $h$. The finite difference approximations are seen to  converge to the H\"older solution at the rate $C_t \Delta t + C_s h^\gamma/\epsilon^2$ where the  constants $C_t$ and $C_s$ are independent of the discretization.  The semi-discrete approximations are found to be stable with time. 
We present numerical simulations for crack propagation that computationally verify the theoretically predicted convergence rate. We also present numerical simulations for crack propagation in pre-cracked samples subject to a bending load.
\end{abstract}


\end{frontmatter}


\section{Introduction}

In \citet{CMPer-Silling} and \cite{States} a self consistent non-local continuum mechanics is proposed. This  formulation known as peridynamics has been employed in the computational reproduction of dynamic fracture as well as offering dynamically based explanations for features observed in fracture, see e.g.,\cite{CMPer-Silling5,CMPer-Silling7,CMPer-Lipton2,BobaruHu,HaBobaru,SillBob,CMPer-Agwai,CMPer-Ghajari}. These references are by no means complete and a recent review of this approach together with further references to the literature can be found in \cite{Handbook}.

The peridynamic formulation expresses internal forces as functions of displacement differences as opposed to  displacement gradients. This generalization allows for an extended  kinematics and provides a unified treatment of differentiable and non-differentiable displacements. 
The motion of a point $\bx$ is influenced by its neighbors through non-local forces. In its simplest formulation forces act within a horizon and only neighbors confined  to a ball of radius $\epsilon$ surrounding $\bx$ can influence the motion of $\bx$. The radius $\epsilon$ is referred to as the peridynamic horizon.
When the forces are linear in the strain and when length scale of nonlocality $\epsilon$ tends to zero the peridynamic models converge to the linear elastic model \cite{CMPer-Emmrich,CMPer-Silling4,AksoyluUnlu,CMPer-Mengesha2}. If one considers non-linear forces associated with two point interactions that are initially elastic and then soften after a critical strain, then the dynamic evolutions are found to converge to a different ``limiting'' dynamics  associated with a crack set and a displacement that satisfies the balance of linear momentum away from the crack set and has bounded elastic energy and Griffith surface energy, see  \cite{CMPer-Lipton,CMPer-Lipton3} and \cite{CMPer-JhaLipton}. A numerical analysis of this two-point interaction or  bond based peridynamic model is carried out in \cite{CMPer-JhaLipton,CMPer-JhaLipton3}. In these works the a-priori convergence rates for finite difference and finite element methods together with different time stepping schemes are reported.

This article focuses on the numerical analysis of a state based peridynamic fracture model governed by forces that are initially elastic and then soften for sufficiently large tensile and hydrostatic strains.  Attention is given to the prototypical  state-based peridynamic model  proposed in \cite{CMPer-Lipton4}. The analysis performed here provides a-priori upper bounds on the convergence rate for a numerical scheme that applies the finite difference approximation in space and the forward Euler discretization scheme in time.  The state based peridynamic model treated here has two components of non-local force acting on a point. The first force is due to tensile strains acting on $\bx$ by its neighbors $\by$, while the second force is due to the net hydrostatic strain on $\bx$ associated with the change in volume about $\bx$. In this article we analyze the convergence of the numerical scheme for two different cases of constitutive law relating non-local force to strain.  For the first case we take both  tensile and hydrostatic forces to be initially linear and increasing with the strain and then after reaching critical values of tensile and hydrostatic strain respectively the forces decrease to zero with strain, see figures \ref{ConvexConcave}(b) and \ref{ConvexConcaveFunctionG}(b).
For the second case we choose the hydrostatic force to be a linear function of the hydrostatic strain  (see dashed line \ref{ConvexConcaveFunctionG}(b)) while the tensile force is initially linear and then decreases to zero after a critical tensile strain is reached, see \autoref{ConvexConcave}(b). The choice of the two constitutive models studied here is motivated by the prospect of simulating materials that exhibit failure due to extreme local tensile stress or strain or materials that fail due to extreme local hydrostatic stress or strain.
Here the quadratic potential function for the dilatational strain can be associated with materials that fail under extreme local tensile loads while the convex-concave dilatational potential function can be associated with materials in which fail under extreme local hydrostatic loads.

The primary new contribution of this paper is that a-priori convergence rates are established for numerical schemes used for simulation using these prototypical state based peridynamic models. As mentioned earlier the constitutive behavior is non-linear, non-convex and material properties can degrade during the course of the evolution.
We consider the class of H\"older continuous displacement fields and show the existence of a unique H\"older continuous evolution for a prescribed H\"older continuous initial condition and body force, see \autoref{thm:existence over finite time domain}. 
To obtain a-priori bounds on the error, we develop an $L^2$ approximation theory for the finite difference approximation in the spatial variables and the forward Euler approximation in time, see \autoref{s:finite difference}. 
We show that discrete approximations converge to the exact H\"older continuous solution uniformly over finite time intervals with respect to the $L^2$ norm. The a-priori rate of convergence in the $L^2$ norm is given by $(C_t\Delta t + C_s h^\gamma/\epsilon^2)$, where $\Delta t$ is the size of the time step, $h$ is the size of spatial mesh discretization, $\gamma\in (0,1]$ is the H\"older exponent, and $\epsilon$ is the length scale of nonlocal interaction relative to the size of the domain, see \autoref{thm:convergence} The constant $C_t$ depends on the $L^2$ norm of the time derivatives of the solution,  $C_s$ depends on the H\"older norm of the solution and the Lipschitz constant of peridynamic force. 
We point out that the convergence results derived here can be extended to general single step time discretization using arguments provided in \cite{CMPer-JhaLipton}.  Although the constitutive law relating force to strain is nonlinear we are still able to establish stability for the semi-discrete approximation and it is shown that the energy at any given time $t$ is bounded above by the energy of the initial conditions and the total work done by the body force up to time $t$, see \autoref{thm:stab semi}.  
Our numerical simulations support the theoretical studies, see \autoref{s:numerical}. In the simulations we introduce a straight crack and it propagates in response  to applied boundary conditions. For these simulations we use piecewise constant interpolants and record the rate of convergence with respect to mesh size while keeping the horizon fixed. Our results show that convergence rate remains above the a-priori estimated rate of $1$ during the simulation. For illustration we also present numerical simulations for a  pre-cracked samples  subject to a bending load.

It is pointed out that there is now a significant number of investigations examining the numerical approximation of singular kernels for non-local problems with applications to nonlocal diffusion, advection, and continuum mechanics.  Numerical formulations and convergence theory for nonlocal $p$-Laplacian formulations are developed in \cite{DeEllaGunzberger}, \cite{Nochetto1}. Numerical analysis of nonlocal steady state diffusion is presented in \cite{CMPer-Du2} and \cite{CMPer-Du3}, and \cite{CMPer-Chen}.  The use of fractional Sobolev spaces for nonlocal problems is investigated and developed in \cite{CMPer-Du1}. Quadrature approximations and stability conditions for linear peridynamics are analyzed in \cite{CMPer-Weckner} and  \cite{CMPer-Silling8}.  The interplay between nonlocal interaction length and grid refinement for linear peridynamic models is presented in \cite{CMPer-Bobaru}. Analysis of adaptive refinement and domain decomposition for the linearized peridynamics are provided in \cite{AksoyluParks}, \cite{LindParks}, and \cite{AksMen}.  This list is by no means complete and the literature continues to  grow rapidly. 

The paper is organized as follows. In \autoref{s:nonlocal dynamics}, we describe the nonlocal model and state the peridynamic equation of motion. The Lipschitz continuity of the peridynamic force and global existence of unique solutions are presented in \autoref{ss:existence holder}. 
The finite difference discretization is introduced in \autoref{s:finite difference}. We demonstrate the energy stability of the semi-discrete approximation in \autoref{semidiscrete}. In \autoref{time discrete} we give the a-priori bound on the error for the fully discrete approximation, see \autoref{thm:convergence}. The numerical simulations are described and presented in \autoref{s:numerical}.  
The Lipschitz continuity of the peridynamic force and  stability of the semi-discrete approximation are proved in \autoref{s:proofs} and \autoref{ss:stab proof}. 
In \autoref{s:conclusions} we summarize our results.

\section{Nonlocal Dynamics}\label{s:nonlocal dynamics}
We now formulate the nonlocal dynamics. Let $D \subset \bbR^d$ denote the material domain of dimension $d=2,3$ and let the horizon be given by $\epsilon >0$. We make the assumption of small (infinitesimal) deformations so that the displacement field $\bu: [0,T]\times D \to \bbR^d$ is small compared to the size of $D$ and the deformed configuration is the same as the reference configuration.  We have $\bu=\bu(t,\bx)$ as a function of space and time but will  suppress the $\bx$ dependence when convenient and write $\bu(t)$. The tensile strain $S$ between two points $\bx,\by \in D$ along the direction $\be_{\by-\bx}$is defined as
\begin{align}\label{strain}
S(\by,\bx,\bu(t))=\frac{\bu(t,\by)-\bu(t,\bx)}{|\by-\bx|}\cdot \be_{\by-\bx},
\end{align}
where $ \be_{\by-\bx}=\frac{\by-\bx}{|\by-\bx|}$ is a unit vector and ``$\cdot$'' is the dot product. The influence function $J^\epsilon(|\by-\bx|)$ is a measure of the influence that the point $\by$ has on $\bx$. Only points inside the horizon can influence $x$ so $J^\epsilon(|\by-\bx|)$ nonzero for $|\by - \bx| < \epsilon$ and zero otherwise. We take $J^\epsilon$ to be of the form: $J^\epsilon(|\by - \bx|) = J(\frac{|\by - \bx|}{\epsilon})$ with $J(r) = 0$ for $r\geq 1$ and $0\leq J(r)\leq M < \infty$ for $r<1$. We also introduce the boundary function $\omega(\bx)$ providing the influence of the boundary on the non-local force. Here $\omega(\bx)$ takes the value $1$, for all $\bx \in D$, an $\epsilon$ distance away from $\partial D$. As $\bx $ approaches $\partial D$ from the interior, $\omega(\bx)$ smoothly decays from $1$ to $0$ on $\partial D$ and is extended by zero outside $D$.

The spherical or hydrostatic strain at $\bx$ is a measure of the volume change about $\bx$ and is given by
\begin{align}\label{sphericalstrain}
\theta(\bx,\bu(t))=\frac{1}{\epsilon^d \omega_d}\int_{H_\epsilon(\bx)} \omega(\by) J^\epsilon(|\by-\bx|)S(\by,\bx,\bu(t)){|\by-\bx|}\,d\by,
\end{align}
where $\omega_d$ is the volume of the unit ball in dimension $d=2,3$, and $H_\epsilon(\bx)$ denotes the ball of radius $\epsilon$ centered at $\bx$.

\subsection{The class of nonlocal potentials}
\begin{figure}
    \centering
    \begin{subfigure}{.45\linewidth}
        \begin{tikzpicture}[xscale=0.75,yscale=0.75]
			\draw [<->,thick] (0,6.5) -- (0,0) -- (2.5,0);
			\draw [-,thick] (0,0) -- (-2.5,0);
			\draw [-] (0,2.15) -- (2.5,2.15);
			\draw [-] (0,6.15) -- (-3.5,6.15);
			\draw [-,thick] (0,0) to [out=0,in=-175] (2.5,2);
			\draw [-,thick] (-3.5,6) to [out=-5,in=180] (0,0);
			

			\draw [-,thick, dashed] (0,0) parabola (2.5,6);
			\draw [-,thick, dashed] (0,0) parabola (-2.5,6);
			
			\draw (1,-0.2) -- (1, 0.2);
			\node [below] at (1,-0.2) {${r}_\ast^+$};
			
			\draw (-1.5,-0.2) -- (-1.5, 0.2);
			\node [below] at (-1.5,-0.2) {${r}_\ast^-$};

			\node [right] at (2.5,0) {$r$};
			\node [left] at (0,2.250) {$C^+_\ast$};
			\node [left] at (0.9,6.1) {$C^-_\ast$};
		  \end{tikzpicture}
		  \caption{}
		  \label{ConvexConcaveFunctionGa}
    \end{subfigure}
    \hskip2em
    \begin{subfigure}{.45\linewidth}
        \begin{tikzpicture}[xscale=0.4,yscale=0.4]
			\draw [<-,thick] (0,3) -- (0,-10);
			\draw [->,thick] (-6.5,0) -- (4,0);
			\draw [-,thick] (-5.5,-10.25) to [out=-30,in=180] (-4.5,-10.5) to [out=0, in=245 ] (0,0) to [out=65,in=180] (1.5,1.5)
			           to [out=0,in=165] (3,0.25);
			           
			\draw [-,thick,dashed] (-3,-9) -- (1,3);
			
			\draw (1.5,-0.2) -- (1.5, 0.2);
			\draw (-4.5,-0.2) -- (-4.5, 0.2);
			\node [below] at (1.5,-0.2) {${r}_\ast^+$};
			\node [below] at (-4.5,-0.2) {${r}_\ast^-$};
			\node [right] at (4,0) {${r}$};
		  \end{tikzpicture}
		  \caption{}
		  \label{ConvexConcaveFunctionGb}
    \end{subfigure}
    
    \caption{(a) Two types of potential function $g(r)$ for hydrostatic  force. The dashed line corresponds to the quadratic potential $g(r) = \beta r^2/2$. The solid line corresponds to the convex-concave type potential $g(r)$. For the convex-concave type potential, there are two special points $r^-_\ast$ and $r^+_\ast$ at which material points start to soften. $C^+_\ast$ and $C^-_\ast$ are two extreme values. (b) Hydrostatic forces. The dashed line corresponds to the quadratic potential and solid line corresponds to the convex-concave potential.}
   \label{ConvexConcaveFunctionG}
\end{figure}
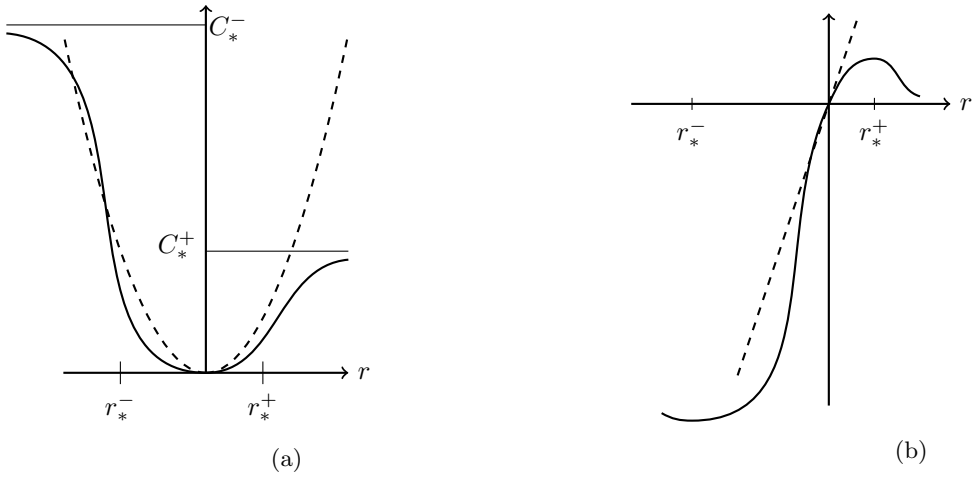

Motivated by potentials of Lennard-Jones type, the force potential for tensile strain is defined by
\begin{equation}\label{tensilepot}
\mathcal{W}^\epsilon (S(\by,\bx,\bu(t)))=\omega(\bx) \omega(\by) J^\epsilon(|\by-\bx|)\frac{1}{\epsilon|\by-\bx|}f(\sqrt{|\by-\bx|}S(\by,\bx,\bu(t)))
\end{equation}
and the potential for hydrostatic strain is defined as
\begin{equation}\label{hydropot}
\mathcal{V}^\epsilon(\theta(\bx,\bu(t)))= \omega(\bx)\frac{g(\theta (\bx,\bu(t)))}{\epsilon^2}\\
\end{equation}
where $\mathcal{W}^\epsilon (S(\by,\bx,\bu(t)))$ is the pairwise force potential per unit length between two points $\bx$ and $\by$ and $\mathcal{V}^\epsilon(\theta(\bx,\bu(t)))$ is the hydrostatic force potential density at $\bx$. They are described in terms of their potential functions $f$ and $g$, see \autoref{ConvexConcave} and \autoref{ConvexConcaveFunctionG}. 

The potential function $f$ represents a convex-concave potential such that the associated force acting between material points $\bx$ and $\by$ are initially elastic and then soften and decay to zero as the strain between points increases, see \autoref{ConvexConcave}. The first well for $\mathcal{W}^\epsilon (S(\by,\bx,\bu(t)))$ is at zero tensile strain and the potential function satisfies
\begin{align}
\label{choice at oregon}
f(0)=f'(0) = 0.
\end{align}
The behavior for infinite tensile strain is characterized by the horizontal asymptotes $\lim_{S\rightarrow \infty} f(S)=C^+$ and $\lim_{S\rightarrow -\infty} f(S)=C^-$ respectively, see \autoref{ConvexConcave}. The critical tensile strain $S^+_c>0$ for which the force begins to soften is given by the inflection point $r^+>0$ of $f$ and is 
\begin{equation}
S^+_c=\frac{r^+}{\sqrt{|y-x|}}.
\label{crittensileplus}
\end{equation}
The critical negative tensile strain is chosen much larger in magnitude than $S_c^+$ and is 
\begin{equation}
S^-_c=\frac{r^-}{\sqrt{|y-x|}},
\label{crittensileminus}
\end{equation}
with $r^-<0$ and $r^+<<|r^-|$.

We assume here that the all the potential functions are bounded and have bounded derivatives up to order $3$, We denote the $i^{th}$ derivative of the function $f$ by $f^{(i)}$, $i=1,2,3$.  Let $C^f_i$ for $i=0,1,2,3$ denote the bounds on the functions and derivatives given by
\begin{align}\label{eq:def Cfi}
C^f_0 := \sup_r |f(r)|, \qquad C^f_i := \sup_r |f^{(i)}(r)| \quad \text{for } i = 1,2,3,
\end{align}
and $C^f_i < \infty$ for $i=0,1,2,3$.

We will consider two types of potentials associated with hydrostatic strain. The first potential we consider is a quadratic potential characterized by a quadratic potential function $g$ with a minimum at zero strain. The second potential we consider is characterized by a convex-concave potential function $g$, see \autoref{ConvexConcaveFunctionG} . If $g$ is assumed to be quadratic then the force due to spherical strain is linear and there is no softening of the material. However, if $g$ is convex-concave the force internal to the material is initially linear and increasing but then becomes decreasing with strain as the hydrostatic strain exceeds a critical value. For the convex-concave $g$, the critical values $0<\theta_c^+$ and $\theta^-_c< 0$ beyond which the force begins to soften is related to the inflection point $r^+_\ast$ and $r_\ast^-$ of $g$ as follows
\begin{equation}
\theta^+_c={r^+_\ast}, \qquad \theta^-_c={r_\ast^-}.
\label{crittensileplus theta}
\end{equation}
The critical compressive hydrostatic strain where the force begins to soften for negative hydrostatic strain is chosen much larger in magnitude than $\theta_c^+$, i.e. $\theta^+_c << |\theta_c^-|$. When $g$ is convex-concave  we assume it is bounded and has bounded derivatives up to order three. These bounds are denoted by  $C^g_i< \infty$ for $i=0,1,2,3$ and,
\begin{align}\label{eq:def Cgi}
C^g_0 := \sup_r |g(r)|, \qquad C^g_i := \sup_r |g^{(i)}(r)| \quad \text{for } i = 1,2,3.
\end{align}



\begin{figure}
    \centering
    \begin{subfigure}{.45\linewidth}
        \begin{tikzpicture}[xscale=0.75,yscale=0.75]
		    \draw [<->,thick] (0,5) -- (0,0) -- (3.0,0);
			\draw [-,thick] (0,0) -- (-3.5,0);
			\draw [-,thin] (0,2.15) -- (2.5,2.15);
			\draw [-, thin] (-3.5,4.15) -- (0,4.15);
			\draw [-,thick] (0,0) to [out=0,in=-175] (2.5,2);
			\draw [-,thick] (-3.5,4) to [out=-5,in=180] (0,0);
			
			\draw (1.5,-0.2) -- (1.5, 0.2);
			\node [below] at (1.5,-0.2) {${r}^+$};
			
			\draw (-2.25,-0.2) -- (-2.25, 0.2);
			\node [below] at (-2.0,-0.2) {${r}^-$};

			\node [right] at (3,0) {$r$};
			\node [left] at (0,2.250) {$C^+$};
			\node [left] at (0.9,4.15) {$C^-$};
			\node [above] at (-2.5,2.20) {$f(r)$};
		  \end{tikzpicture}
		  \caption{}
		  \label{ConvexConcavea}
    \end{subfigure}
    \hskip2em
    \begin{subfigure}{.45\linewidth}
        \begin{tikzpicture}[xscale=0.6,yscale=0.6]
		    \draw [<-,thick] (0,3) -- (0,-3);
			\draw [->,thick] (-5,0) -- (3.5,0);
			\draw [-,thick] (0,0) to [out=65,in=180] (1.5,1.5) to [out=0,in=165] (3,0.25);
			
			\draw [-,thick] (-4.6,-0.5) to [out=-20,in=130] (-3.0,-2.0) to [out=-50, in=245] (0,0);
			
			\draw (1.5,-0.2) -- (1.5, 0.2);
			\draw (-2.25,-0.2) -- (-2.25, 0.2);
			\node [below] at (1.5,-0.2) {${r}^+$};
			\node [below] at (-2.0,-0.2) {${r}^-$};
			\node [right] at (3.5,0) {${r}$};
			\node [right] at (0,2.2) {$f'(r)$};
		  \end{tikzpicture}
		   \caption{}
		   \label{ConvexConcaveb}
    \end{subfigure}
    \caption{(a) The potential function $f(r)$  for tensile force. Here $C^+$ and $C^-$ are the two asymptotic values of $f$. (b) Cohesive tensile force.}\label{ConvexConcave}
\end{figure}
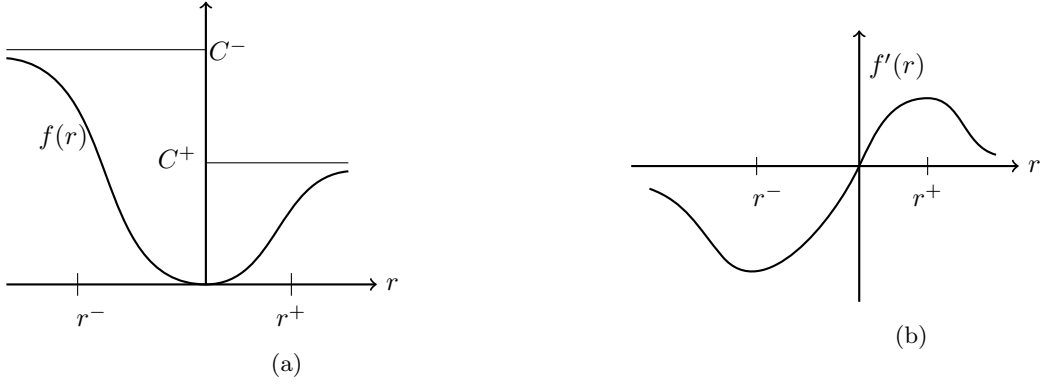

\subsection{Peridynamic equation of motion}
The potential energy of the motion is given by
\begin{equation}\label{new peri}
\begin{aligned}
PD^\epsilon(\bu)=\frac{1}{\epsilon^d \omega_d}\int_D \int_{H_\epsilon(\bx)} |\by-\bx|\mathcal{W}^\epsilon(S(\by,\bx,\bu(t)))\,d\by d\bx\\
+\int_D \mathcal{V}^\epsilon(\theta(\bx,\bu(t)))\,d\bx.
\end{aligned}
\end{equation}
In this treatment the material is assumed homogeneous and the density $\rho$ is constant. We denote the body force by $\bb(\bx,t)$ and define the {Lagrangian}  $${\rm{L}}(\bu,\partial_t \bu,t)=\frac{\rho}{2}||\dot \bu||^2 _{L^2 (D;\mathbb{R}^d)}-PD^\epsilon(\bu)+\int_D \bb\cdot \bu d\bx,$$
where $\dot \bu=\frac{\partial \bu}{\partial t}$ is the velocity and $\Vert \dot \bu\Vert_{L^2(D;\mathbb{R}^d)}$ denotes the $L^2$ norm of the vector field $\dot \bu: D\rightarrow \mathbb{R}^d$.
Applying the {principal of least action} gives the nonlocal dynamics 
\begin{equation}\label{energy based model2}
\begin{aligned}
\rho \ddot{\bu}(\bx,t)=\mathcal{L}^\epsilon(\bu)(\bx,t)+\bb(\bx,t),\hbox{  for  $\bx\in D$},
\end{aligned}
\end{equation}
where 
\begin{align}\label{eq:total peri force}
\mathcal{L}^\epsilon(\bu)(\bx,t) = \mathcal{L}^\epsilon_T(\bu)(\bx,t) + \mathcal{L}^\epsilon_D(\bu)(\bx,t).
\end{align}
Here $\mathcal{L}^\epsilon_T(\bu)$ is the peridynamic force due to the tensile strain and is given by
\begin{align}\label{nonlocforcetensite}
&\mathcal{L}^\epsilon_T(\bu)(\bx,t) \notag \\
&=\frac{2}{\epsilon^d \omega_d}\int_{H_\epsilon(\bx)} \omega(\bx) \omega(\by) \frac{J^\epsilon(|\by-\bx|)}{\epsilon|\by-\bx|}\partial_S f(\sqrt{|\by-\bx|}S(\by,\bx,\bu(t)))\be_{\by-\bx}\,d\by,
\end{align}
and $\mathcal{L}^\epsilon_D(\bu)$ is the peridynamic force due to the hydrostatic strain and is given by
\begin{align}\label{nonlocforcedevia}
&\mathcal{L}^\epsilon_D(\bu)(\bx,t)\notag \\
&=\frac{1}{\epsilon^d \omega_d}\int_{H_\epsilon(\bx)} \omega(\bx) \omega(\by) \frac{J^\epsilon(|\by-\bx|)}{\epsilon^2}\left[\partial_\theta g(\theta(\by,\bu(t)))+\partial_\theta g(\theta(\bx,\bu(t)))\right]\be_{\by-\bx}\,d\by.
\end{align}

The dynamics is complemented with the initial data 
\begin{align}\label{idata}
\bu(\bx,0)=\bu_0(\bx), \qquad  \partial_t \bu(\bx,0)=\bv_0(\bx),
\end{align}
and we prescribe zero Dirichlet boundary condition on the boundary $\partial D$
\begin{align}\label{eq:bc}
\bu(\bx) = \bzero \qquad \forall \bx \in \partial D.
\end{align}
The zero boundary value is extended outside $D$ by zero to  $\bbR^d$. Last we note that since the material is homogeneous we will divide both sides of the equation of motion by $\rho$ and assume, without loss of generality, that $\rho = 1$.



\section{Existence of solutions}\label{ss:existence holder}
Let $\Cholder{\perdrd}$ be the H\"{o}lder space with exponent $\gamma \in (0,1]$. We introduce $C_0^{0,\gamma}(D)=C^{0,\gamma}(D)\cap C_0(D)$ where $C_0(D)$ is the closure of continuous functions with compact support on $D$ in the supremum norm. Functions in $C_0(D)$ are uniquely extended to $\overline{D}$ and take zero values  on $\partial D$, see \cite{MA-Driver}. In this paper we extend all functions in $C_0^{0,\gamma}(D)$ by zero outside $D$. The norm of $\bu \in \Cholderz{\perdrd}$ is given by 
\begin{align*}
\Choldernorm{\bu}{\perdrd} &:= \sup_{\bx \in D} \abs{\bu(\bx)} + \left[\bu \right]_{\Cholder{\perdrd}},
\end{align*}
where $\left[\bu \right]_{\Cholder{\perdrd}}$ is the H\"{o}lder semi norm and given by
\begin{align*}
\left[\bu \right]_{\Cholder{\perdrd}} &:= \sup_{\substack{\bx\neq \by,\\
\bx,\by \in D}} \dfrac{\abs{\bu(\bx)-\bu(\by)}}{\abs{\bx - \by}^\gamma},
\end{align*}
and $\Cholderz{\perdrd}$ is a Banach space with this norm. Here we make the hypothesis that the domain function $\omega$ belongs to $C^{0,\gamma}_0(D;[0,1])$.

We consider the first order system of equations equivalent to \autoref{energy based model2}. Let $y_1(t) = \bu(t)$, $y_2(t) = {\bv}(t)$ with $\bv(t) = \dot{\bu}(t)$. We form the vector $y = (y_1,y_2)^T$ where $y_1,y_2 \in \Cholderz{\perdrd}$ and let $F^{\epsilon}(y,t) = (F^{\epsilon}_1(y,t), F^{\epsilon}_2(y,t))^T$ with
\begin{align}
F^\epsilon_1(y,t) &:= y_2 \label{eq:per first order eqn 1} \\
F^\epsilon_2(y, t) &:= \mathcal{L}^\epsilon(y_1(t)) + \bb(t). \label{eq:per first order eqn 2}
\end{align}
The initial boundary value associated with the evolution \autoref{energy based model2} is equivalent to the initial boundary  value problem for the first order system given by
\begin{align}\label{eq:per first order}
\dfrac{d}{dt}y = F^{\epsilon}(y,t),
\end{align}
with initial condition given by $y(0) = (\bu_0, \bv_0)^T \in \Cholderz{\perdrd}\times\Cholderz{\perdrd}$.

We next show that $F^\epsilon(y,t)$ is Lipschitz continuous. 

{\vskip 2mm}
\begin{proposition}\label{prop:lipschitz}
\textbf{Lipschitz continuity and bound}\\
Let $X = \Cholderz{\perdrd} \times \Cholderz{\perdrd}$. We suppose that the boundary function $\omega$ belongs to $C^{0,\gamma}_0(D;[0,1])$. Let $f$ be a convex-concave potential  function satisfying \autoref{eq:def Cfi} and let the potential function $g$ either be a quadratic function or be a convex-concave function satisfying  \autoref{eq:def Cgi}, then the function $F^\epsilon(y,t) = (F^\epsilon_1, F^\epsilon_2)^T$, as defined in \autoref{eq:per first order eqn 1} and \autoref{eq:per first order eqn 2}, is Lipschitz continuous in any bounded subset of $X$. We have, for any $y,z \in X$ and $t> 0$,
\begin{align}\label{eq:lipschitz property of F}
&\normX{F^{\epsilon}(y,t) - F^{\epsilon}(z,t)}{X} \notag \\
&\leq \frac{L_1 (1+||\omega||_{\Cholder{}}) (1+ ||y||_X + ||z||_X)}{\epsilon^{2+\alpha(\gamma)}} ||y - z||_X.
\end{align}
where $L_1$ is independent of $\bu,\bv$ and $\epsilon$, and depends on $f$, $J$, and $g$. The exponent $\alpha(\gamma)$ is $0$ if $\gamma \geq 1/2$ and is $1/2 - \gamma$ if $\gamma \leq 1/2$. 
Furthermore, for any $y \in X$ and any $t\in [0,T]$, we have the bound
\begin{align}\label{eq:bound on F}
\normX{F^\epsilon(y,t)}{X} &\leq \frac{L_2(1+||\omega||_{\Cholder{}}) (1+||y||_{X})}{\epsilon^2} + b,
\end{align}
where $b = \sup_{t} \Choldernorm{\bb(t)}{\perdrd}$ and $L_2$ is independent of $y$. 
\end{proposition}
{\vskip 2mm}

We easily see that on choosing $z=0$ in \autoref{eq:lipschitz property of F} that $\mathcal{L}^\epsilon(\bu) $ is in $C^{0,\gamma}(D;\mathbb{R}^d)$ provided that $\bu$ belongs to $C^{0,\gamma}(D;\mathbb{R}^3)$. Moreover since  $\mathcal{L}^\epsilon(\bu)$ takes the value $\bzero$ on $\partial D$ we can conclude that $\mathcal{L}^\epsilon(\bu)$ also belongs to  $C^{0,\gamma}_0(D;\mathbb{R}^d)$.


The following theorem gives the existence and uniqueness of solution in any given time domain $I_0 = (-T, T)$.

{\vskip 2mm}
\begin{theorem}\label{thm:existence over finite time domain} 
\textbf{Existence and uniqueness of H\"older solutions over finite time intervals}\\
Let $f$ be a convex-concave function satisfying \autoref{eq:def Cfi} and let $g$ either be a quadratic function or a convex-concave function satisfying \autoref{eq:def Cgi}. For any initial condition $x_0\in X ={ \Cholderz{\perdrd} \times \Cholderz{\perdrd}}$, time interval $I_0=(-T,T)$, and right hand side $\bb(t)$  continuous in time for $t\in I_0$ such that $\bb(t)$ satisfies $\sup_{t\in I_0} {||\bb(t)||_{\Cholder{}}}<\infty$, there is a unique solution $y(t)\in C^1(I_0;X)$ of
\begin{equation}
y(t)=x_0+\int_0^tF^\epsilon(y(\tau),\tau)\,d\tau,
\label{10}
\end{equation}
or equivalently
\begin{equation}
y'(t)=F^\epsilon(y(t),t),\hbox{with    $y(0)=x_0$},
\label{11}
\end{equation}
where $y(t)$ and $y'(t)$ are Lipschitz continuous in time for $t\in I_0$.
\end{theorem}
{\vskip 2mm}
The proof of this theorem follows directly from \sautoref{Proposition}{prop:lipschitz} and is established along the same lines as the existence proof for H\"older continuous solutions of bond based peridynamics given in  [Theorem 2,\cite{CMPer-JhaLipton}]. 

We conclude this section by stating following result which shows the Lipschitz bound of peridynamic force in $L^2$ norm for functions in $L^2_0(D; \bbR^d)$. Here $L^2_0(D;\bbR^d)$ denotes the space of functions $\bu \in L^2(D;\bbR^d)$ such that $\bu = \bzero$ on $\partial D$. We assume that functions in $L^2_0(D;\bbR^d)$ are extended to $\bbR^d$ by zero. 

{\vskip 2mm}
\begin{proposition}\label{prop:lipschitz L2}
\textbf{Lipschitz continuity of peridynamic force in $L^2$}\\
Let $f$ and $g$ satisfy the hypothesis of Proposition \sautoref{Proposition}{prop:lipschitz}, then for any $\bu,\bv \in L^2_0(D;\bbR^d)$ we have
\begin{align}\label{eq:lipschitz L2}
||\mathcal{L}^\epsilon(\bu) - \mathcal{L}^\epsilon(\bv)||_{L^2(D;\bbR^d)} &\leq \dfrac{L_3}{\epsilon^2} ||\bu - \bv||_{L^2(D;\bbR^d)},
\end{align}
where the constants $L_3$ and $L_4$ are independent of $\epsilon$, $\bu$ and $\bv$. Here  $L_3 = 4 (C^f_1 \bar{J}_1 + C^g_2 \bar{J}_0^2)$, for convex-concave $g$, and $L_3 = 4 (C^f_1 \bar{J}_1 + g''(0) \bar{J}_0^2)$, for quadratic $g$.  Here $\bar{J}_\alpha  = \frac{1}{\omega_d} \int_{H_1(\bzero)} J(|\bxi|) |\bxi|^{-\alpha} d\bxi$.
\end{proposition}
{\vskip 2mm}

The proofs of \sautoref{Proposition}{prop:lipschitz} and \sautoref{Proposition}{prop:lipschitz L2} are provided in \autoref{s:proofs}. We now describe the finite difference scheme and analyze the rate of convergence to H\"older continuous solutions of the peridynamic equation of motion.

\section{Finite difference approximation}
\label{s:finite difference}
In this section we consider the discrete approximation to the dynamics given by finite differences in space and the forward Euler discretization in time.
\begin{figure}[h]
\centering
\includegraphics[scale=0.5]{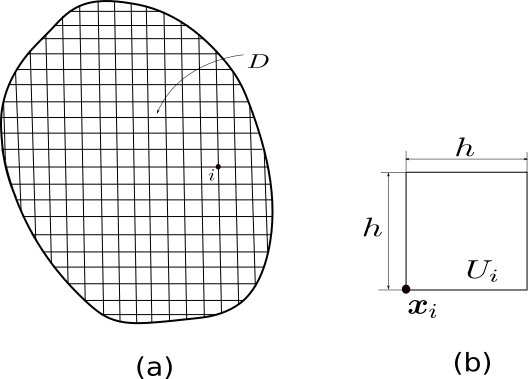}
\caption{(a) Typical mesh of size $h$. (b) Unit cell $U_i$ corresponding to material point $\bx_i$.}\label{fig:peridynamic mesh}
\end{figure} Let $h$ denote the mesh size and $D_h = D\cap(h \bbZ)^d$ be the associated discretization of the material domain $D$. In this paper we will keep the horizon length scale $\epsilon$ fixed and assume that the spatial discretization length satisfies $h< \epsilon<1$. Let $i\in \bbZ^d$ be the index such that $\bx_i = hi \in D$, see \autoref{fig:peridynamic mesh}. Let $U_i$ be a the cell of volume $h^d$ corresponding to the grid point $\bx_i$. The exact solution evaluated at grid points is denoted by $(\bu_i(t),\bv_i(t))$. Given any discrete set $\{ \hat{\bu}_i \}_{i, \bx_i\in D}$, where $i$ is index representing grid point of mesh, we define its piecewise constant extension as
\begin{align}\label{eq:const extension u}
\hat{\bu}(\bx) := \sum_{i, \bx_i \in D} \hat{\bu}_i \chi_{U_i}(\bx). 
\end{align}
In this way we have representation of the discrete set as a piecewise constant function.

We now describe the $L^2$-projection of the function $\bu : D\to \bbR^d$ onto the space of piecewise constant functions defined over the cells $U_i$. We denote the average of $\bu$ over the unit cell $U_i$ as  $\tilde{\bu}_i$ and
\begin{align}
\tilde{\bu}_i &:= \frac{1}{h^d} \int_{U_i} \bu(\bx) d\bx
\end{align}
and the $L^2$ projection of $\bu$ onto piecewise constant functions is $\tilde{\bu}$ given by
\begin{align}\label{eq:projn u}
\tilde{\bu}(\bx) &:= \sum_{i,\bx_i \in D} \tilde{\bu}_i \chi_{U_i}(\bx).
\end{align}

\begin{lemma}\label{lem:proj error}
Let $\bu \in C^{0,\gamma}_0(D;\bbR^d)$ and let $\tilde{\bu}$ be its $L^2$ projection defined in \autoref{eq:projn u}, then we have
\begin{align}\label{eq:proj error}
|\tilde{\bu}(\bx) - \bu(\bx)| &\leq \left[c^\gamma ||\bu||_{\Cholder{}}\right] h^\gamma, \forall \bx \in D, \notag \\
||\tilde{\bu}(\bx) - \bu(\bx)||_{L^2} &\leq \left[c^\gamma \sqrt{|D|} ||\bu||_{\Cholder{}} \right] h^\gamma,
\end{align}
where $c = \sqrt{2}$ for $d=2$ and $c = \sqrt{3}$ for $d=3$. 
\end{lemma}
This lemma can be demonstrated easily by substituting \autoref{eq:projn u} for $\tilde{\bu}$ and using the fact that $\bu \in \Cholderz{\perdrd}$. We also note that first line of  \autoref{eq:proj error} remains valid of $\bx$ in a layer of thickness $2\epsilon$ surrounding $D$.

\subsection{Stability of the semi-discrete approximation}
\label{semidiscrete}
We first introduce the semi-discrete boundary condition by setting $\hat{\bu}_i(t) = \bzero$ for all $t$ and for all $\bx_i \notin D$.
Let $\{\hat{\bu}_i(t)\}_{i,\bx_i\in D}$ denote the semi-discrete approximate solution which satisfies the following, for all $t\in [0,T]$ and $i$ such that $\bx_i\in D$,
\begin{align}\label{eq:fd semi discrete}
\ddot{\hat{\bu}}_i(t) = \mathcal{L}^\epsilon(\hat{\bu}(t))(\bx_i) + \bb(\bx_i,t),
\end{align}
where $\hat{\bu}(t)$ is the piecewise constant extension of discrete set $\{\hat{\bu}_i(t) \}_i$ and is defined as
\begin{equation}
\label{eq:def piecewise ext}
\hat{\bu}(\bx,t) =
\begin{cases}
 \sum_{i, \bx_i \in D} \hat{\bu}_i(t) \chi_{U_i}(\bx),\\
 \bzero,\,\,\,i,\hbox{ such that }\bx_i\not\in D.
\end{cases}
\end{equation}
The scheme is complemented with the discretized initial conditions $\hat{\bu}_i(0) = \bu_0(\bx_i)$ and $\hat{\bv}_i(0) =\bv_0(\bx_i)$. 

The total kinetic and potential energy is given by
\begin{align*}
\mathcal{E}^\epsilon(\bu)(t) = \frac{1}{2} ||\dot{\bu}(t)||^2_{L^2} + PD^\epsilon(\bu(t)),
\end{align*} 
and we introduce the augmented energy given by
\begin{align}\label{eq:per energy total new}
\bar{\mathcal{E}}^\epsilon(\bu)(t) := \mathcal{E}^\epsilon(\bu)(t) + \frac{1}{2}||\bu(t)||^2_{L^2}. 
\end{align}

We have the stability of the semi-discrete evolution. 

{\vskip 2mm}
\begin{theorem}\label{thm:stab semi}
\textbf{Energy stability of the semi-discrete approximation}\\
Let $\{\hat{\bu}_i(t) \}_{i,\bx_i\in D}$ be the solution to the semidiscrete initial boundary value problem \autoref{eq:fd semi discrete} and $\hat{\bu}(t)$ denote its piecewise constant extension. Similarly let $\hat{\bb}(t,\bx)$ denote the piecewise constant extension of $\{ \bb(t,\bx_i)\}_{i,\bx_i\in D}$. 
If  $f$ and $g$ are convex-concave type functions satisfying \autoref{eq:def Cfi} and \autoref{eq:def Cgi}, then the total energy $\mathcal{E}^\epsilon(\buhat)(t)$ satisfies,
\begin{align}\label{eq:inequal energy doublewell g}
\mathcal{E}^\epsilon(\hat{\bu})(t) &\leq \left( \sqrt{\mathcal{E}^\epsilon(\hat{\bu})(0)} + \dfrac{t C}{\epsilon^2} + \int_0^t ||\hat{\bb}(s)||_{L^2} ds \right)^2, \hbox{  $\forall t \in [0,T]$,}
\end{align}
and the constant $C$ is independent of $\epsilon$ and $h$. 

If  $f$ is  a convex-concave type function satisfying \autoref{eq:def Cfi} and $g$ is quadratic  then the augmented energy $\bar{\mathcal{E}}^\epsilon(\buhat)(t)$ satisfies, 
\begin{align}\label{eq:inequal energy quadratic g}
\bar{\mathcal{E}}^\epsilon(\hat{\bu})(t) &\leq \exp[3(C_2/\epsilon^2 +1) t] \bigg( \bar{\mathcal{E}}^\epsilon(\hat{\bu})(0) \notag \\
&\quad \quad + \int_0^T (\frac{C_1^2}{\epsilon^4} + ||\hat{\bb}(s)||^2_{L^2}) \exp[-3(C_2/\epsilon^2 +1) s] ds \bigg),\hbox{  $\forall t \in [0,T]$},
\end{align}
where the constants $C_1$ and $C_2$ are independent of $\epsilon$ and $h$. 
\end{theorem}
{\vskip 2mm}

We provide proof of \autoref{thm:stab semi} in \autoref{ss:stab proof}. We now discuss the fully discrete scheme.

\subsection{Time discretization}
\label{time discrete}
Let $\Delta t$ be the size of the time step and $[0,T] \cap (\Delta t \bbZ)$ be the discretization of the time domain. We denote the fully discrete solution at $(t^k = k\Delta t, \bx_i = ih)$ as $(\hat{\bu}^k_{i}, \hat{\bv}^k_i)$ and the exact solution as $(\bu^k_i,\bv^k_i)$. We enforce the boundary condition $\hat{\bu}^k_i = \bzero$ for all $\bx_i \notin D$ and for all $k$. The piecewise constant extension of $\{\hat{\bu}^k_i\}_{i\in \bbZ^d}$ and $\{\hat{\bv}^k_i\}_{i\in \bbZ^d}$ are denoted by $\hat{\bu}^k$ and $\hat{\bv}^k$ respectively. The $L^2$-projection of $\bu^k$ and $\bv^k$ onto piecewise constant functions are denoted by $\tilde{\bu}^k$ and $\tilde{\bv}^k$ respectively. 

The forward Euler time discretization, with respect to velocity, and the finite difference scheme for $(\hat{\bu}^k_{i}, \hat{\bv}^k_i)$ is written
\begin{align}
\dfrac{\buhat^{k+1}_i - \buhat^k_i}{\Delta t} &= \bvhat^{k+1}_i \label{eq:finite diff eqn u} \\
\dfrac{\bvhat^{k+1}_i - \bvhat^k_i}{\Delta t} &= \mathcal{L}^\epsilon(\hat{\bu}^k)(\bx_i) + \bb^k_i. \label{eq:finite diff eqn v}
\end{align}
The initial condition is enforced by setting $\buhat^{0}_i =(\buhat_0)_i$ and $\bvhat^{0}_i =(\bvhat_0)_i$. We note that the forward difference scheme for the system reduces to the central difference scheme for  the second order differential equation \autoref{energy based model2} on substitution of \autoref{eq:finite diff eqn u} into \autoref{eq:finite diff eqn v}. 

%

\subsubsection{Convergence of approximation}

In this section we  provide an upper bound on the convergence rate of the fully discrete approximation to the H\"older continuous solution as measured by the $L^2$ norm.
The $L^2$ approximation error $E^k$ at time $t^k$, for $0< t^k \leq T$, given by
\begin{align*}
E^k &:= \Ltwonorm{\buhat^k - \bu^k}{\perdrd} + \Ltwonorm{\bvhat^k- \bv^k}{\perdrd} .
\end{align*}
The following theorem gives an explicit a-priori upper bound on the convergence rate.

{\vskip 2mm}
\begin{theorem}\label{thm:convergence}
\textbf{Convergence of finite difference approximation (forward Euler time discretization)}\\
Let $\epsilon>0$ be fixed. Let $(\bu, \bv)$ be the solution of peridynamic equation \autoref{eq:per first order}. We assume $\bu, \bv \in \Ctwointime{[0,T]}{\Cholderz{\perdrd}}$.   Then the finite difference scheme given by \autoref{eq:finite diff eqn u} and \autoref{eq:finite diff eqn v} is consistent in both  time and spatial discretization and converges to the exact solution uniformly in time with respect to the $\Ltwo{\perdrd}$ norm. If we assume the error at the initial step is zero then the error $E^k$ at time $t^k$ is bounded and satisfies 
\begin{align}\label{eq: first est}
\sup_{0\leq k \leq T/\Delta t} E^k\leq O\left( C_t\Delta t + C_s\dfrac{h^\gamma}{\epsilon^2} \right),
\end{align}
where constant $C_s$ and $C_t$ are independent of $h$ and $\Delta t$ and $C_s$ depends on the H\"older norm of the solution and $C_t$ depends on the $L^2$ norms of time derivatives of the solution. 
\end{theorem}
{\vskip 2mm}
Here we have assumed the initial error is zero for ease of exposition only. 

We remark that the explicit constants leading to \autoref{eq: first est} can be large. The inequality that delivers \autoref{eq: first est}  is given by
\begin{align}\label{eq: fund est initial}
\sup_{0\leq k \leq T/\Delta t} E^k\leq \exp \left[T (1 + L_3/\epsilon^2) \right] T \left[ C_t \Delta t + (C_s/\epsilon^2) h^\gamma \right],
\end{align}
where the constants $L_3$, $C_t$ and $C_s$ are given by \autoref{eq:def L3},  \autoref{eq:const Ct}, and \autoref{eq:const Cs}.  The explicit constant $C_t$ depends on the spatial $L^2$ norm of the time derivatives of the solution and $C_s$ depends on the spatial H\"older continuity of the solution  and  the constant $L_3$. The constant $L_3$ is bounded independently of horizon $\epsilon$.  Although the constants are necessarily pessimistic they deliver a-priori error estimates. We carry out numerical simulations for different values of the horizon $\epsilon$ in \autoref{s:numerical}. We find that the convergence rate for piecewise constant finite difference interpolation functions is greater than or equal to $\gamma=1$ for simulations lasting in the tens of microseconds. These results are consistent with the a-priori estimates given in \autoref{thm:convergence} above.

\subsubsection{Error analysis}\label{ss:error analysis1}
We split the error between $(\buhat^k, \bvhat^k)^T$ and $(\bu^k,\bv^k)^T$ in two parts as follows
\begin{align}
E^k &= ||\buhat^k - \bu^k||_{L^2} + ||\bvhat^k- \bv^k||_{L^2} \notag \\
&\leq \left[ ||\butilde^k - \bu^k||_{L^2} + ||\bvtilde^k- \bv^k||_{L^2} \right] + \left[ ||\buhat^k - \butilde^k||_{L^2} + ||\bvhat^k- \bvtilde^k||_{L^2} \right].
\end{align}
In section \autoref{ss:error analysis} we will show that the error between the $L^2$ projections of the actual solution and the discrete approximation for both forward Euler and implicit one step methods decay according to
\begin{align}\label{eq:write estimate error ek}
\sup_{0\leq k \leq T/\Delta t} \left( ||\buhat^k - \butilde^k||_{L^2} + ||\bvhat^k- \bvtilde^k||_{L^2} \right)  &= O\left( \Delta t + \dfrac{h^\gamma}{\epsilon^2} \right).
\end{align}
And using \sautoref{Lemma}{lem:proj error} we have
\begin{align*}
& \sup_{k} \left( ||\butilde^k - \bu^k||_{L^2} + ||\bvtilde^k- \bv^k||_{L^2} \right) \notag \\
&= c^\gamma \sqrt{|D|} \left[\sup_{t\in [0,T]} ||\bu(t)||_{\Cholder{}} + \sup_{t\in [0,T]} ||\bv(t)||_{\Cholder{}} \right] h^\gamma.
\end{align*}

We now study the difference $\buhat^k - \butilde^k$ and $\bvhat^k - \bvtilde^k$. 


\subsubsection{Error analysis for approximation of $L^2$ projection of the exact solution}\label{ss:error analysis}

Let the differences be denoted by $\be^k(u) := \buhat^k - \butilde^k$ and  $\be^k(v ):= \bvhat^k - \bvtilde^k$ and their evaluation at grid points are $\be^k_i(u) := \buhat^k_i - \butilde^k_i$ and $\be^k_i(v) := \bvhat^k_i - \bvtilde^k_i$. We have the following lemma for the evolution of the differences in the discrete dynamics.

\begin{lemma}
\label{lem: discretedifferences}
The differences $\be^k_i(u)$ and $\be^k_i(v)$ discretely evolve according to the equations:
\begin{align}\label{eq:error eqn in u}
\be^{k+1}_i(u) &= \be^k_i(u) + \Delta t \be^{k+1}_i(v) + \Delta t\tau^{k}_i(u)
\end{align}
and
\begin{align}\label{eq:error eqn in v}
\be^{k+1}_i(v) &= \be^k_i(v) + \Delta t \left(\tau^k_i(v) + \sigma^k_i(u) + \sigma^k_i(v) \right) \notag \\
&\quad + \Delta t \left( \mathcal{L}^{\epsilon}(\buhat^k)(\bx_i) - \mathcal{L}^{\epsilon}(\butilde^k)(\bx_i) \right).
\end{align}
Here $\tau^k_i(u),\tau^k_i(v)$ and $\sigma^k_i(u),\sigma^k_i(v)$ are consistency error terms and are defined as
\begin{align}\label{eq:def consistency error}
\tau^k_i(u) &:= \dparder{\butilde^k_i}{t} - \dfrac{\butilde^{k+1}_i - \butilde^k_i}{\Delta t}, \notag  \\
\tau^k_i(v) &:= \dparder{\bvtilde^k_i}{t} - \dfrac{\bvtilde^{k+1}_i - \bvtilde^k_i}{\Delta t} , \notag \\
\sigma^k_i(u) &:= \left( \mathcal{L}^{\epsilon}(\butilde^k)(\bx_i) - \mathcal{L}^{\epsilon}(\bu^k)(\bx_i)  \right) \notag  \\
\sigma^k_i(v) &:=  \dparder{\bv^k_i}{t} - \dparder{\bvtilde^k_i}{t}.
\end{align}
\end{lemma}

To prove this we start by subtracting $(\butilde^{k+1}_i - \butilde^k_i)/\Delta t$ from \autoref{eq:finite diff eqn u} to get
\begin{align*}
& \dfrac{\buhat^{k+1}_i - \buhat^k_i}{\Delta t} - \dfrac{\butilde^{k+1}_i - \butilde^k_i}{\Delta t}  \\
&= \bvhat^{k+1}_i - \dfrac{\butilde^{k+1}_i - \butilde^k_i}{\Delta t} \notag \\
&= \bvhat^{k+1}_i - \bvtilde^{k+1}_i + \left( \bvtilde^{k+1}_i - \dparder{\butilde^{k+1}_i}{t} \right) + \left( \dparder{\butilde^{k+1}_i}{t} - \dfrac{\butilde^{k+1}_i - \butilde^k_i}{\Delta t} \right).
\end{align*}
Taking the average over unit cell $U_i$ of the exact peridynamic equation \autoref{eq:per first order} at time $t^k$, we will get $\bvtilde^{k+1}_i - \dparder{\butilde^{k+1}_i}{t}  = 0$ and we recover \autoref{eq:error eqn in u}. 

Next, we subtract $(\bvtilde^{k+1}_i - \bvtilde^k_i)/\Delta t$ from \autoref{eq:finite diff eqn v} and add and subtract terms to get
\begin{align}\label{eq:error in v 1}
\dfrac{\bvhat^{k+1}_i - \bvhat^k_i}{\Delta t} - \dfrac{\bvtilde^{k+1}_i - \bvtilde^k_i}{\Delta t} &=   \mathcal{L}^\epsilon(\hat{\bu}^k)(\bx_i) + \bb^k_i - \dparder{\bv^k_i}{t} + \left( \dparder{\bv^k_i}{t} - \dfrac{\bvtilde^{k+1}_i - \bvtilde^k_i}{\Delta t}\right) \notag \\
&=  \mathcal{L}^\epsilon(\hat{\bu}^k)(\bx_i)  + \bb^k_i - \dparder{\bv^k_i}{t} \notag \\
&\quad + \left( \dparder{\bvtilde^k_i}{t} - \dfrac{\bvtilde^{k+1}_i - \bvtilde^k_i}{\Delta t}\right) + \left( \dparder{\bv^k_i}{t} - \dparder{\bvtilde^k_i}{t}\right),
\end{align}
where from \autoref{eq:def consistency error}
\begin{align}\label{eq:consistency error in v}
\tau^k_i(v) &:= \dparder{\bvtilde^k_i}{t} - \dfrac{\bvtilde^{k+1}_i - \bvtilde^k_i}{\Delta t}.
\end{align}
Note that from the exact peridynamic equation, we have
\begin{align}\label{eq:exact per eqn v 1}
\bb^k_i - \dparder{\bv^k_i}{t} =-\mathcal{L}^\epsilon({\bu}^k)(\bx_i) .
\end{align}
Combining \autoref{eq:error in v 1}, \autoref{eq:consistency error in v}, and \autoref{eq:exact per eqn v 1}, gives
\begin{align*}
\be^{k+1}_i(v) &= \be^k_i(v) + \Delta t \tau^k_i(v) + \Delta t \left( \dparder{\bv^k_i}{t} - \dparder{\bvtilde^k_i}{t}\right)  \notag \\
&\quad + \Delta t \left( \mathcal{L}^\epsilon(\hat{\bu}^k)(\bx_i)  -\mathcal{L}^\epsilon({\bu}^k)(\bx_i)  \right) \notag \\
&= \be^k_i(v) + \Delta t \tau^k_i(v) + \Delta t \left( \dparder{\bv^k_i}{t} - \dparder{\bvtilde^k_i}{t}\right) \notag \\
&\quad + \Delta t \left( \mathcal{L}^\epsilon(\hat{\bu}^k)(\bx_i)  -\mathcal{L}^\epsilon(\tilde{\bu}^k)(\bx_i)  \right) \notag \\
&\quad + \Delta t \left( \mathcal{L}^\epsilon(\tilde{\bu}^k)(\bx_i)  - \mathcal{L}^\epsilon({\bu}^k)(\bx_i)  \right).
\end{align*}
and the lemma follows on applying the definitions given in \autoref{eq:def consistency error}.

\subsubsection{Consistency}\label{sss:consistency}
In this section we provide upper bounds on the consistency errors. This error is measured in the $L^2$ norm.  Here the upper bound on the consistency error with respect to time follows using Taylor's series expansion. The upper bound on the spatial consistency error is established using the H\"older continuity of nonlocal forces.

\textbf{Time discretization: } We apply a Taylor series expansion in time to estimate $\tau^k_i(u)$ as follows
\begin{align*}
\tau^k_i(u) &= \dfrac{1}{h^d} \int_{U_i} \left( \dparder{\bu^k(\bx)}{t} - \dfrac{\bu^{k+1}(\bx) - \bu^k(\bx)}{\Delta t} \right) d\bx \\
&= \dfrac{1}{h^d} \int_{U_i} \left( -\dfrac{1}{2} \dsecder{\bu^k(\bx)}{t} \Delta t + O((\Delta t)^2) \right) d\bx .
\end{align*}
We form the $\Ltwo{}$ norm of $\tau^k_i(u)$ and apply Jensen's inequality to get
\begin{align*}
\Ltwonorm{\tau^k(u)}{} &\leq \frac{\Delta t}{2} \Ltwonorm{\dsecder{\bu^k}{t}}{} + O((\Delta t)^2) \notag \\
&\leq \frac{\Delta t}{2} \sup_{t} \Ltwonorm{\dsecder{\bu(t)}{t}}{} + O((\Delta t)^2).
\end{align*}
A similar argument gives
\begin{align*}
\Ltwonorm{\tau^k(v)}{} = \frac{\Delta t}{2} \sup_{t} \Ltwonorm{\dsecder{\bv(t)}{t}}{} + O((\Delta t)^2).
\end{align*}

\

\textbf{Spatial discretization: }From \autoref{eq:def consistency error} one can write
\begin{align*}
\sigma^k_i(v) &= \frac{\partial \bv^k_i}{\partial t} - \frac{\partial \bvtilde^k_i}{\partial t} = \frac{\partial \bv^k(\bx_i)}{\partial t} - \frac{\partial \bvtilde^k(\bx_i)}{\partial t}.
\end{align*}
Applying \sautoref{Lemma}{lem:proj error} gives
\begin{align*}
|\sigma^k_i(v)| &\leq c^\gamma h^\gamma  \Choldernorm{\dparder{\bv^k}{t}}{} \leq c^\gamma h^\gamma \sup_t \Choldernorm{\dparder{\bv(t)}{t}}{}.
\end{align*}
Taking the $L^2$ norm and using the estimates given above yields the inequality
\begin{align*}
\Ltwonorm{\sigma^k(v)}{} &\leq h^{\gamma} c^{\gamma} \sqrt{ \abs{D} } \sup_{t} \Choldernorm{\dparder{\bv(t)}{t}}{}.
\end{align*}


We now estimate $\abs{\sigma^k_i(u)}$. Since $\mathcal{L}^\epsilon = \mathcal{L}^\epsilon_T + \mathcal{L}^\epsilon_D$, we have from \autoref{eq:def consistency error}
\begin{align}\label{eq:estimate sigma u 1}
|\sigma^k_i(u)| &\leq |\mathcal{L}^\epsilon_T(\butilde^k)(\bx_i) - \mathcal{L}^\epsilon_T(\bu^k)(\bx_i)| + |\mathcal{L}^\epsilon_D(\butilde^k)(\bx_i) - \mathcal{L}^\epsilon_D(\bu^k)(\bx_i)| \notag \\
&= I_1 + I_2
\end{align}
To expedite the calculations we employ the following notation for $\bxi \in H_1(\bzero)$,
\begin{align}\label{eq:notations consist}
s_\bxi &:= \epsilon |\bxi|, \: \be_\bxi := \frac{\bxi}{|\bxi|}, \notag \\
\omega_\bxi(\bx) &:= \omega(\bx + \epsilon \bxi) \omega(\bx), \notag \\
\bubar_\bxi(\bx) &:= \bu(\bx+ \epsilon \bxi) - \bu(\bx) .
\end{align}
 We also write hydrostatic strain (see \autoref{sphericalstrain}) as follows
\begin{align}\label{eq:hydro strain}
\theta(\bx;\bu) = \frac{1}{\omega_d} \int_{H_1(\bzero)}\omega(\bx + \epsilon\bxi) J(|\bxi|)\bubar_\bxi(\bx)\cdot \be_\bxi d\bxi.
\end{align}

In our calculations we will also encounter various moments of influence function $J$ therefore we define following term 
\begin{align}\label{eq:def bar J}
\bar{J}_\alpha := \frac{1}{\omega_d} \int_{H_1(\bzero)} J(|\bxi|) |\bxi|^{-\alpha} d\bxi, \qquad \text{ for } \alpha \in \bbR. 
\end{align}
Recall that $J(|\bxi|) = 0$ for $\bxi \notin H_1(\bzero)$ and $0\leq J(|\bxi|) \leq M$ for $\bxi\in H_1(\bzero)$.

Applying the notation $\mathcal{L}^\epsilon_T$ becomes
\begin{align}\label{eq:per bondbased simple force expression}
\mathcal{L}^\epsilon_T(\bu)(\bx) &= \frac{2}{\epsilon \omega_d} \int_{H_1(\bzero)} \omega_\bxi(\bx) \frac{J(|\bxi|)}{\sqrt{s_\bxi}} f'(\bubar_\bxi(\bx)\cdot \be_\bxi/\sqrt{s_\bxi}) \be_{\bxi} d\bxi.
\end{align}
On choosing $\bu = \bu^k$ and $\bu= \butilde^k$ in $\mathcal{L}^\epsilon_T$ given by \autoref{eq:per bondbased simple force expression} we get
\begin{align}\label{eq:consist est 0}
I_1 &\leq \frac{2}{\epsilon \omega_d} \int_{H_1(\bzero)} \omega_\bxi(\bx_i) \frac{J(|\bxi|)}{\sqrt{s_\bxi}} |f'(\bar{\butilde}^k_\bxi(\bx_i) \cdot \be_\bxi/\sqrt{s_\bxi}) - f'(\bar{\bu}^k_\bxi(\bx_i) \cdot \be_\bxi/\sqrt{s_\bxi})| d\bxi \notag \\
&\leq \frac{2C^f_2}{\epsilon \omega_d} \int_{H_1(\bzero)} \frac{J(|\bxi|)}{s_\bxi} |\bar{\butilde}^k_\bxi(\bx_i)- \bar{\bu}^k_\bxi(\bx_i)| d\bxi,
\end{align}
where we have applied \autoref{eq:def Cfi} and used the fact that $|f'(r_1) - f'(r_2)| \leq C^f_2 |r_1 - r_2|$ and $0\leq \omega(\bx) \leq 1$. We use \sautoref{Lemma}{lem:proj error} to estimate $|\bar{\butilde}^k_\bxi(\bx_i)- \bar{\bu}^k_\bxi(\bx_i)|$ as follows
\begin{align}
|\bar{\butilde}^k_\bxi(\bx_i)- \bar{\bu}^k_\bxi(\bx_i)| &\leq |\butilde^k(\bx_i + \epsilon \bxi)- \bu^k(\bx_i+ \epsilon\bxi)| + |\butilde^k(\bx_i )- \bu^k(\bx_i)| \notag \\
&\leq 2 c^\gamma ||\bu(t^k)||_{\Cholder{}} h^\gamma \leq 2 c^\gamma \sup_t ||\bu(t)||_{\Cholder{}} h^\gamma.
\end{align}
From this we get
\begin{align}
I_1 &\leq \left[\frac{4C^f_2 c^\gamma \bar{J}_1}{\epsilon^2} \sup_t ||\bu(t)||_{\Cholder{}} \right] h^\gamma,
\end{align}
where $\bar{J}_\alpha$ for $\alpha \in \bbR$ is defined in \autoref{eq:def bar J}. Clearly, 
\begin{align}\label{eq:consist est 1}
\sum_{i,\bx_i \in D} h^d I_1^2  &\leq \left[\frac{4C^f_2 c^\gamma \bar{J}_1 \sqrt{|D|}}{\epsilon^2} \sup_t ||\bu(t)||_{\Cholder{}} \right]^2 h^{2\gamma} .
\end{align}

We now estimate $I_2$ in \autoref{eq:estimate sigma u 1}. We will consider $g$ of convex-concave type satisfying $C^g_i < \infty$ for $i=0,1,2,3$ where $C^g_0 = \sup |g(r)|$ and $C^g_i = \sup |g^{(i)}(r)|$ for $i=1,2,3$. It is noted that the upper bound for the choice of quadratic $g$ is also found using the steps presented here. We can write $\mathcal{L}^\epsilon_D(\bu)(\bx)$ (see \autoref{nonlocforcedevia}) as follows
\begin{align}\label{eq:per statebased simple force expression}
\mathcal{L}^\epsilon_D(\bu)(\bx) &= \frac{1}{\epsilon^2 \omega_d} \int_{H_1(\bzero)} \omega_\bxi(\bx) J(|\bxi|) [g'(\theta(\bx+\epsilon\bxi;\bu)) + g'(\theta(\bx;\bu))] \be_\bxi d\bxi.
\end{align}
Using this expression we have the upper bound 
\begin{align}\label{eq:consist est 2}
I_2 &= \bigg\vert \frac{1}{\epsilon^2 \omega_d} \int_{H_1(\bzero)} \omega_\bxi(\bx_i) J(|\bxi|) [g'(\theta(\bx_i+\epsilon\bxi;\butilde^k)) + g'(\theta(\bx_i;\butilde^k)) \notag \\
&\quad \quad - g'(\theta(\bx_i+\epsilon\bxi;\bu^k)) + g'(\theta(\bx_i;\bu^k))] \be_\bxi d\bxi \bigg\vert \notag \\
&\leq \frac{1}{\epsilon^2 \omega_d} \int_{H_1(\bzero)} J(|\bxi|) (|g'(\theta(\bx_i+\epsilon\bxi;\butilde^k)) - g'(\theta(\bx_i+\epsilon \bxi;\bu^k))| \notag \\
&\quad \quad + |g'(\theta(\bx_i;\butilde^k)) - g'(\theta(\bx_i;\bu^k))|) d\bxi \notag \\
&\leq \frac{C^g_2}{\epsilon^2 \omega_d} \int_{H_1(\bzero)} J(|\bxi|) (|\theta(\bx_i+\epsilon\bxi;\butilde^k) - \theta(\bx_i+\epsilon \bxi;\bu^k)| \notag \\
&\quad \quad + |\theta(\bx_i;\butilde^k) - \theta(\bx_i;\bu^k)|) d\bxi .
\end{align}
We proceed further as follows using expression of $\theta$ in \autoref{eq:hydro strain}
\begin{align}
& |\theta(\bx_i+\epsilon\bxi;\butilde^k) - \theta(\bx_i+\epsilon \bxi;\bu^k)| \notag \\
&\leq \bigg\vert \frac{1}{\omega_d} \int_{H_1(\bzero)} \omega(\bx_i + \epsilon \bxi + \epsilon \boldsymbol{\eta}) J(|\boldsymbol{\eta}|) (\butilde^k(\bx_i + \epsilon \bxi + \epsilon \boldsymbol{\eta}) \notag \\
&\quad\quad - \bu^k(\bx_i + \epsilon \bxi + \epsilon \boldsymbol{\eta}) - \butilde^k(\bx_i + \epsilon \bxi) + \bu^k(\bx_i + \epsilon \bxi))\cdot \be_{\boldsymbol{\eta}} d\boldsymbol{\eta} \bigg\vert \notag \\
&\leq \frac{1}{\omega_d} \int_{H_1(\bzero)} J(|\boldsymbol{\eta}|) (|\butilde^k(\bx_i + \epsilon \bxi + \epsilon \boldsymbol{\eta}) - \bu^k(\bx_i + \epsilon \bxi + \epsilon \boldsymbol{\eta})| \notag \\
&\quad \quad + |\butilde^k(\bx_i + \epsilon \bxi) - \bu^k(\bx_i + \epsilon \bxi)|)d\boldsymbol{\eta} \notag \\
&\leq 2 c^\gamma h^\gamma \sup_t ||\bu(t)||_{\Cholder{}} \bar{J}_0
\end{align}
where we used \sautoref{Lemma}{lem:proj error} in last step. We combine above estimate in \autoref{eq:consist est 2} to get
\begin{align}
I_2 &\leq \left[ \frac{4C^g_2c^\gamma \bar{J}_0^2 }{\epsilon^2} \sup_t ||\bu(t)||_{\Cholder{}} \right] h^\gamma
\end{align}
and
\begin{align}\label{eq:consist est 3}
\sum_{i,\bx_i \in D} h^d I_2^2 &\leq \left[ \frac{4C^g_2c^\gamma \bar{J}_0^2 \sqrt{|D|} }{\epsilon^2} \sup_t ||\bu(t)||_{\Cholder{}} \right]^2 h^{2\gamma} .
\end{align}
Applying \autoref{eq:consist est 1}, \autoref{eq:consist est 3} and \autoref{eq:estimate sigma u 1} gives
\begin{align*}
||\sigma^k(u)||_{L^2} &\leq \sqrt{\sum_{i,\bx_i \in D} h^d I_1^2} + \sqrt{\sum_{i,\bx_i \in D} h^d I_2^2} \notag \\
&\leq \left[ \frac{4(C^g_2 \bar{J}^2_0 + C^f_2 \bar{J}_1)c^\gamma \sqrt{|D|} }{\epsilon^2} \sup_t ||\bu(t)||_{\Cholder{}} \right] h^\gamma.
\end{align*}
Here we define the constant
\begin{equation}\label{eq:def L3}
L_3 = \left\{\begin{array}{l l}4 (C^f_1 \bar{J}_1 + C^g_2 \bar{J}^2_0),& \hbox{ for $g$ convex-concave}\\4 (C^f_1 \bar{J}_1 + g''(0) \bar{J}^2_0),& \hbox {for $g$ quadratic} \end{array} \right.
\end{equation}
this is also the Lipschitz constant related to Lipschitz continuity of peridynamic force in $L^2$, see  \sautoref{Proposition}{prop:lipschitz L2}. Thus, we have shown for $g$ convex-concave that
\begin{align}
||\sigma^k(u)||_{L^2} &\leq \left[ \frac{L_3 c^\gamma \sqrt{|D|} }{\epsilon^2} \sup_t ||\bu(t)||_{\Cholder{}} \right] h^\gamma.
\end{align}
The same arguments show that an identical inequality holds for quadratic $g$ using the other definition of $L_3$ and this completes the estimation of the consistency errors.

\subsubsection{Stability}\label{sss:stability}
In this subsection we establish estimates that ensure stability of the evolution and apply the consistency estimates of the previous subsection to establish \autoref{thm:convergence}.
Let $e^k$ be the total error at the $k^{\text{th}}$ time step. It is defined as
\begin{align*}
e^k &:=  \Ltwonorm{\be^k(u)}{} + \Ltwonorm{\be^k(v)}{}.
\end{align*}
To simplify the calculations, we collect all the consistency errors  and write them as 
\begin{align*}
\tau &:= \sup_t \left(\Ltwonorm{\tau^k(u)}{}  + \Ltwonorm{\tau^k(v)}{}  + \Ltwonorm{\sigma^k(u)}{} + \Ltwonorm{\sigma^k(v)}{}\right),
\end{align*}
and from our consistency analysis, we know that to leading order in $\Delta t$  and $h$ that
\begin{align}\label{eq:estimate tau}
\tau &\leq C_t \Delta t + \dfrac{C_s}{\epsilon^2} h^\gamma
\end{align}
where,
\begin{align}
C_t &:= \frac{1}{2} \sup_{t} \Ltwonorm{\dsecder{\bu(t)}{t}}{} + \frac{1}{2} \sup_{t} \Ltwonorm{\dfrac{\partial^3 \bu(t)}{\partial t^3}}{}, \label{eq:const Ct} \\
C_s &:= c^\gamma \sqrt{|D|} \left[ \epsilon^2 \sup_{t} \Choldernorm{\dfrac{\partial^2 \bu(t)}{\partial t^2}}{} + L_3 \sup_t  \Choldernorm{\bu(t)}{} \right]. \label{eq:const Cs}
\end{align}

We take the $\Ltwo{}$ norm of \autoref{eq:error eqn in u} and \autoref{eq:error eqn in v} and add them. Using the definition of $\tau$ we get
\begin{align}\label{eq:error k ineq 1}
e^{k+1} &\leq e^k + \Delta t \Ltwonorm{\be^{k+1}(v)}{\perdrd} + \Delta t \tau \notag \\
&\quad + \Delta t \sqrt{ \sum_{i, \bx_i \in D} h^d \abs{ \mathcal{L}^{\epsilon}(\buhat^k)(\bx_i) - \mathcal{L}^{\epsilon}(\butilde^k)(\bx_i)}^2 }.
\end{align}

It now remains to estimate the last term in the above equation. We illustrate the calculations for convex-concave $g$ noting the identical steps apply to quadratic $g$ as well. Let 
\begin{align}\label{eq:consist est 3.1}
H &:= \sqrt{ \sum_{i, \bx_i \in D} h^d \abs{ \mathcal{L}^{\epsilon}(\buhat^k)(\bx_i) - \mathcal{L}^{\epsilon}(\butilde^k)(\bx_i)}^2} \notag \\
&\leq \sqrt{ \sum_{i, \bx_i \in D} h^d \abs{ \mathcal{L}^{\epsilon}_T(\buhat^k)(\bx_i) - \mathcal{L}^{\epsilon}_T(\butilde^k)(\bx_i)}^2 } \notag \\
&\quad+ \sqrt{ \sum_{i, \bx_i \in D} h^d \abs{ \mathcal{L}^{\epsilon}_D(\buhat^k)(\bx_i) - \mathcal{L}^{\epsilon}_D(\butilde^k)(\bx_i)}^2 } \notag \\
&=: H_1 + H_2.
\end{align}
Choosing $\bu = \buhat^k$ and $\bu= \butilde^k$ with $\mathcal{L}^\epsilon_T$ given by \autoref{eq:per bondbased simple force expression} we get
\begin{align}\label{eq:consist est 4}
H_1^2 &\leq \sum_{i, \bx_i \in D} h^d \bigg\vert \frac{2C^f_2}{\epsilon^2 \omega_d} \int_{H_1(\bzero)} \frac{J(|\bxi|)}{|\bxi|} |\bar{\buhat}^k_\bxi(\bx_i)- \bar{\butilde}^k_\bxi(\bx_i)| d\bxi \bigg\vert^2 ,
\end{align}
where $\bar{\buhat}^k_\bxi(\bx) = \buhat^k(\bx + \epsilon \bxi) - \buhat^k(\bx )$. 

We will make use of the following inequality in the sequel.  Let $p(\bxi)$ be a scalar valued function of $\bxi$ and $\alpha \in \bbR$ then
\begin{align}\label{eq:ineq symm square}
&\bigg\vert \frac{C}{\omega_d} \int_{H_1(\bzero)} \frac{J(|\bxi|)}{|\bxi|^\alpha} p(\bxi) d\bxi \bigg\vert^2 \notag \\
&\leq \left(\frac{C}{\omega_d} \right)^2 \int_{H_1(\bzero)} \int_{H_1(\bzero)} \frac{J(|\bxi|)}{|\bxi|^\alpha} \frac{J(|\boldsymbol{\eta}|)}{|\boldsymbol{\eta}|^\alpha} p(\bxi) p(\boldsymbol{\eta}) d\bxi d\boldsymbol{\eta}  \notag \\
&\leq \left(\frac{C}{\omega_d} \right)^2 \int_{H_1(\bzero)} \int_{H_1(\bzero)} \frac{J(|\bxi|)}{|\bxi|^\alpha} \frac{J(|\boldsymbol{\eta}|)}{|\boldsymbol{\eta}|^\alpha} \frac{p(\bxi)^2 + p(\boldsymbol{\eta})^2}{2} d\bxi d\boldsymbol{\eta} \notag \\
&=  C^2 \frac{\bar{J}_\alpha}{\omega_d } \int_{H_1(\bzero)} \frac{J(|\bxi|)}{|\bxi|^\alpha} p(\bxi)^2 d\bxi.
\end{align}

On applying \autoref{eq:ineq symm square} in \autoref{eq:consist est 4} with $C = \frac{2C^f_2}{\epsilon^2}$, $\alpha = 1$, $p(|\bxi|) = |\bar{\buhat}^k_\bxi(\bx_i)- \bar{\butilde}^k_\bxi(\bx_i)|$ we get
\begin{align}
H_1^2 &\leq \sum_{i, \bx_i \in D} h^d \left(\frac{2C^f_2}{\epsilon^2}\right)^2 \frac{\bar{J}_1}{\omega_d} \int_{H_1(\bzero)} \frac{J(|\bxi|)}{|\bxi|} |\bar{\buhat}^k_\bxi(\bx_i)- \bar{\butilde}^k_\bxi(\bx_i)|^2 d\bxi \notag \\
&\leq \left(\frac{2C^f_2}{\epsilon^2}\right)^2 \frac{\bar{J}_1}{\omega_d} \int_{H_1(\bzero)} \frac{J(|\bxi|)}{|\bxi|} \notag \\
&\quad \left[ \sum_{i, \bx_i \in D} h^d 2(|\buhat^k(\bx_i + \epsilon \bxi) - \butilde^k(\bx_i + \epsilon \bxi)|^2 + |\buhat^k(\bx_i) - \butilde^k(\bx_i)|^2) \right] d\bxi \notag \\
&\leq \left(\frac{2C^f_2}{\epsilon^2}\right)^2 \frac{\bar{J}_1}{\omega_d} \int_{H_1(\bzero)} \frac{J(|\bxi|)}{|\bxi|} \notag \\
&\quad \left[ \sum_{i, \bx_i \in D} h^d 2(|\be^k(u)(\bx_i + \epsilon \bxi)|^2 + |\be^k(u)(\bx_i)|^2) \right] d\bxi,
\end{align}
where we substituted definition of $\bar{\buhat}^k_\bxi$ and $\bar{\butilde}^k_\bxi$ and used inequality $(a+b)^2 \leq 2 a^2 + 2 b^2$ in third step, and identified terms as $\be^k(u)$ in last step. Since $\be^k(u)(\bx) = \sum_{i,\bx_i \in D} \be^k_i(u) \chi_{U_i}(\bx)$, we have
\begin{align}\label{eq:consist est 5}
H_1^2 &\leq \left(\frac{2C^f_2}{\epsilon^2}\right)^2 \frac{\bar{J}_1}{\omega_d} \int_{H_1(\bzero)} \frac{J(|\bxi|)}{|\bxi|}  4 ||\be^k(u)||_{L^2}^2 d\bxi, \notag
\end{align}
so
\begin{equation}\label{eq:consist est for H1}
H_1 \leq \frac{4C^f_2 \bar{J}_1}{\epsilon^2} ||\be^k(u)||_{L^2} .
\end{equation}

We now estimate $H_2$. Note that for $I_2 = |\mathcal{L}^\epsilon_D(\butilde^k)(\bx_i) - \mathcal{L}^\epsilon_D(\bu^k)(\bx_i)|$, we have the inequality given by \autoref{eq:consist est 2}. We now use \autoref{eq:consist est 2} but with $\butilde^k$ replaced by $\buhat^k$ and $\bu^k$ replaced by $\butilde^k$ together with the identity $\theta(\bx;\buhat^k) - \theta(\bx; \butilde^k) = \theta(\bx;\buhat^k - \butilde^k) = \theta(\bx; \be^k(u))$, to see that
\begin{align}
H_2^2 &\leq \sum_{i, \bx_i \in D} h^d \bigg( \frac{C^g_2}{\epsilon^2 \omega_d} \int_{H_1(\bzero)} J(|\bxi|) (|\theta(\bx_i+\epsilon\bxi;\be^k(u))|+ |\theta(\bx_i;\be^k(u))|) d\bxi \bigg)^2.
\end{align}
We use inequality \autoref{eq:ineq symm square} with $C = C^g_2/\epsilon^2$, $\alpha = 0$, and $p(\bxi ) = |\theta(\bx_i+\epsilon\bxi;\be^k(u))|+ |\theta(\bx_i;\be^k(u))|$ to get
\begin{align}\label{eq:consist est 6}
H_2^2 &\leq \sum_{i, \bx_i \in D} h^d  \left( \frac{C^g_2}{\epsilon^2}\right)^2 \frac{\bar{J}_0}{\omega_d } \int_{H_1(\bzero)} J(|\bxi|) ( |\theta(\bx_i+\epsilon\bxi;\be^k(u))|+ |\theta(\bx_i;\be^k(u))| )^2 d\bxi \notag \\
&\leq \left( \frac{C^g_2}{\epsilon^2}\right)^2 \frac{\bar{J}_0}{\omega_d } \int_{H_1(\bzero)} J(|\bxi|) \notag \\
& \quad \left[ \sum_{i, \bx_i \in D} h^d 2 (|\theta(\bx_i+\epsilon\bxi;\be^k(u))|^2+ |\theta(\bx_i;\be^k(u))|^2 ) \right] d\bxi,
\end{align}
where we used inequality $(a+b)^2 \leq 2a^2 + 2b^2$ in the second step. We now proceed to estimate the first sum in the last line of \autoref{eq:consist est 6},
\begin{align}
& \sum_{i, \bx_i \in D} h^d |\theta(\bx_i+\epsilon\bxi;\be^k(u))|^2 \notag \\
&\leq \sum_{i, \bx_i \in D} h^d \bigg\vert \frac{1}{\omega_d }\int_{H_1(\bzero)} J(|\boldsymbol{\eta}|) \omega(\bx_i + \epsilon \bxi + \epsilon\boldsymbol{\eta}) \notag \\
& \qquad ( \be^k(u)(\bx_i + \epsilon \bxi + \epsilon\boldsymbol{\eta}) - \be^k(u) (\bx_i + \epsilon \bxi)) \cdot \be_{\boldsymbol{\eta}} d\boldsymbol{\eta} \bigg\vert^2 \notag \\
&\leq \sum_{i, \bx_i \in D} h^d \bigg( \frac{1}{\omega_d }\int_{H_1(\bzero)} J(|\boldsymbol{\eta}|) (|\be^k(u)(\bx_i + \epsilon \bxi + \epsilon\boldsymbol{\eta})| + |\be^k(u)(\bx_i + \epsilon \bxi)|) d\boldsymbol{\eta} \bigg)^2,
\end{align}
where we used expression of $\theta$ from \autoref{eq:hydro strain} in first step, and used $0\leq \omega(\bx) \leq 1$ in the second step. The second summation on the last line of \autoref{eq:consist est 6} is  also bounded above the same way.  We apply inequality \autoref{eq:ineq symm square} with $C =1$, $\alpha = 0$, and $p(\boldsymbol{\eta}) = |\be^k(u)(\bx_i + \epsilon \bxi + \epsilon\boldsymbol{\eta})| + |\be^k(u)(\bx_i + \epsilon \bxi)|$ to get
\begin{align}\label{eq:consist est 7}
& \sum_{i, \bx_i \in D} h^d |\theta(\bx_i+\epsilon\bxi;\be^k(u))|^2 \notag \\
&\leq \sum_{i, \bx_i \in D} h^d \frac{\bar{J}_0}{\omega_d} \int_{H_1(\bzero)} J(|\boldsymbol{\eta}|) (|\be^k(u)(\bx_i + \epsilon \bxi + \epsilon\boldsymbol{\eta})| + |\be^k(u)(\bx_i + \epsilon \bxi)|)^2 d \boldsymbol{\eta} \notag \\
&\leq \frac{\bar{J}_0}{\omega_d} \int_{H_1(\bzero)} J(|\boldsymbol{\eta}|) 2 \sum_{i, \bx_i \in D} h^d(|\be^k(u)(\bx_i + \epsilon \bxi + \epsilon\boldsymbol{\eta})|^2 + |\be^k(u)(\bx_i + \epsilon \bxi)|^2) d \boldsymbol{\eta} \notag \\ 
&\leq \frac{\bar{J}_0}{\omega_d} \int_{H_1(\bzero)} J(|\boldsymbol{\eta}|) 4 ||\be^k(u)||^2_{L^2} d \boldsymbol{\eta} \notag \\
& = 4 \bar{J}^2_0 ||\be^k(u)||^2_{L^2},
\end{align}
where as before we have used the Cauchy inequality. We next apply the estimate \autoref{eq:consist est 7} to \autoref{eq:consist est 6} to see that
\begin{align}\label{eq:consist est 8}
&H_2^2 \leq  16 \bar{J}^2_0 ||\be^k(u)||^2_{L^2} \left( \frac{C^g_2}{\epsilon^2}\right)^2 \frac{\bar{J}_0}{\omega_d } \int_{H_1(\bzero)} J(|\bxi|) d\bxi, \notag
\end{align}
so
\begin{equation}\label{eq:consist est H2}
H_2 \leq  \frac{4 C^g_2 \bar{J}^2_0}{\epsilon^2} ||\be^k(u)||_{L^2} .
\end{equation}

Finally, we apply the inequalities  given by \autoref{eq:consist est for H1} and \autoref{eq:consist est H2} to \autoref{eq:consist est 3.1} and obtain
\begin{align}\label{eq:consist est 9}
H &=  \sqrt{\sum_{i, \bx_i \in D} h^d \abs{ \mathcal{L}^{\epsilon}(\buhat^k)(\bx_i) - \mathcal{L}^{\epsilon}(\butilde^k)(\bx_i)}^2 } \notag \\
&\leq H_1 + H_2 \notag \\
&\leq  \frac{4(C^f_2 \bar{J}_1 + C^g_2 \bar{J}^2_0)}{\epsilon^2 } ||\be^k(u)||_{L^2}  \notag \\
&\leq \left[ \frac{L_3}{\epsilon^2 } ||\be^k(u)||_{L^2} \right]^2,
\end{align}
where $L_3=4(C^f_2 \bar{J}_1 + C^g_2 \bar{J}^2_0)$ for convex-concave $g$. For the case of quadratic $g$ we have the same inequality but with $L_3=4(C^f_2 \bar{J}_1 + g''(0) \bar{J}^2_0)$.

Applying the inequality given by \autoref{eq:consist est 9} to \autoref{eq:error k ineq 1} gives
\begin{align*}
e^{k+1} &\leq e^k + \Delta t \Ltwonorm{\be^{k+1}(v)}{\perdrd} + \Delta t \tau + \Delta t \dfrac{L_3}{\epsilon^2}\Ltwonorm{\be^k(u)}{\perdrd}
\end{align*}
We now add $\Delta t ||e^{k+1}(u)||_{L^2(D;\mathbb{R}^d)} + \Delta t \dfrac{L_3}{\epsilon^2} ||\be^k(v)||_{L^2(D;\mathbb{R}^d)} $ to the right side of the equation above to get
\begin{align*}
&e^{k+1} \leq ( 1 + \Delta t \dfrac{L_3}{\epsilon^2})  e^k + \Delta t e^{k+1} + \Delta t \tau \\
\Rightarrow & e^{k+1} \leq \dfrac{( 1 + \Delta t {L_3}/{\epsilon^2})}{1 - \Delta t} e^{k} + \dfrac{\Delta t }{1 - \Delta t} \tau.
\end{align*}
We now recursively substitute $e^{j}$ as follows
\begin{align}
e^{k+1} &\leq \dfrac{( 1 + \Delta t {L_3}/\epsilon^2)}{1 - \Delta t} e^k + \dfrac{\Delta t }{1 - \Delta t} \tau  \notag \\
&\leq \left(\dfrac{( 1 + \Delta t {L_3}/\epsilon^2)}{1 - \Delta t} \right)^2 e^{k-1} + \dfrac{\Delta t }{1 - \Delta t} \tau  \left(1 + \dfrac{( 1 + \Delta t {L_3}/\epsilon^2)}{1 - \Delta t}\right) \notag \\
&\leq ...\notag \\
&\leq \left(\dfrac{( 1 + \Delta t {L_3}/\epsilon^2)}{1 - \Delta t} \right)^{k+1} e^0 + \dfrac{\Delta t }{1 - \Delta t} \tau  \sum_{j=0}^k \left(\dfrac{( 1 + \Delta t {L_3}/\epsilon^2)}{1 - \Delta t} \right)^{k-j}. \label{eq:ek estimnate}
\end{align}

Since $1/(1-\Delta t)= 1 + \Delta t + \Delta t^2 + O(\Delta t^3)$, we have
\begin{align*}
\dfrac{( 1 + \Delta t {L_3}/\epsilon^2)}{1 - \Delta t} &\leq 1 + (1 + {L_3}/\epsilon^2) \Delta t + (1 + {L_3}/\epsilon^2) \Delta t^2 + O({L_3}/\epsilon^2) O(\Delta t^3).
\end{align*}
Now, for any $k \leq T /\Delta t$ and using the identity $(1+ a)^k \leq \exp [ka]$ for $a\leq 0$, we have
\begin{align*}
&\left( \dfrac{ 1 + \Delta t {L_3}/\epsilon^2}{1 - \Delta t} \right)^k \\
&\leq \exp \left[k (1 + {L_3}/\epsilon^2) \Delta t + k(1 + {L_3}/\epsilon^2) \Delta t^2 + k O(L_3/\epsilon^2) O(\Delta t^3) \right] \\
&\leq \exp \left[T (1 + {L_3}/\epsilon^2) + T (1 + {L_3}/\epsilon^2) \Delta t + O(T{L_3}/\epsilon^2) O(\Delta t^2) \right].
\end{align*}
We write above equation in more compact form as follows
\begin{align*}
&\left( \dfrac{ 1 + \Delta t {L_3}/\epsilon^2}{1 - \Delta t} \right)^k \\
&\leq \exp \left[T (1 + {L_3}/\epsilon^2) (1 + \Delta t + O(\Delta t^2)) \right]. 
\end{align*}

We use above estimate in \autoref{eq:ek estimnate} and get following inequality for $e^k$
\begin{align*}
e^{k+1} &\leq \exp \left[T (1 + {L_3}/\epsilon^2) (1 + \Delta t + O(\Delta t^2)) \right] \left( e^0 + (k+1) \tau \Delta t/(1- \Delta t) \right) \notag \\
&\leq \exp \left[T (1 + {L_3}/\epsilon^2) (1 + \Delta t + O(\Delta t^2)) \right] \left( e^0 + T\tau (1 + \Delta t + O(\Delta t^2) \right).
\end{align*}
where we used the fact that $1/(1-\Delta t) = 1+ \Delta t + O(\Delta t^2)$.

Assuming the error in initial data is zero, i.e. $e^0= 0$, and noting the estimate of $\tau$ in \autoref{eq:estimate tau}, we have
\begin{align*}
&\sup_k e^k \leq \exp \left[T (1 + L_3/\epsilon^2) \right] T \tau 
\end{align*}
and we conclude to leading order that
\begin{align}\label{eq: fund est}
\sup_k e^k \leq \exp \left[T (1 + L_3/\epsilon^2) \right] T \left[ C_t \Delta t + (C_s/\epsilon^2) h^\gamma \right],
\end{align}
Here the constants $C_t$ and $C_s$ are given by \autoref{eq:const Ct} and \autoref{eq:const Cs}. 
This shows the stability of the numerical scheme. We note that constant $L_3 = 4 (C^f_1 \bar{J}_1 + C^g_2 \bar{J}_0^2)$, where $C^f_2 = \sup |f''(r)|, C^g_2 = \sup |g''(r)|$, corresponds to the case when $g$ is convex-concave type. For quadratic $g$ the constant is given by $L_3 = 4 (C^f_1 \bar{J}_1 + g''(0) \bar{J}_0^2)$. 

\section{Numerical results}\label{s:numerical}
In this section, we present numerical simulations that support the theoretical upper bound on the convergence rate and to illustrate the displacement field and fracture set under different loading conditions. 

We specify the density $\rho = 1200 \,kg/m^3$, bulk modulus $K = 25 \, GPa$, and critical energy release rate $G_c = 500 \,Jm^{-2}$. The pairwise interaction and the hydrostatic interaction are characterized by potentials $f(r) = c (1-\exp [-\beta r^2])$ and $g(r) = \bar{C} r^2/2$ respectively. The influence function is $J(r) = 1-r$. Equations 94, 95, and 97 of  \cite{CMPer-Lipton4} relate parameters $c, \beta, \bar{C}$ to the Lam\`e parameters $\lambda, \mu$ and the critical energy release rate $G_c$. In \autoref{tab:mat props}, we list the value of constants for Poisson's ratio $0.22$ and $0.245$. These ratios are computed using the relation established in \cite{CMPer-Lipton4}. The critical bond strain between material point $y$ and $x$ is $S_c = \bar{r}/\sqrt{|y-x|}$ where $\bar{r} = 1/\sqrt{2\beta}$. 

We consider the central difference time discretization described by \autoref{eq:finite diff eqn u} and \autoref{eq:finite diff eqn v} on a uniform square mesh of mesh size $h$.
We can write the peridynamic force $\mathcal{L}^\epsilon(\hat{\bu}^k)(\bx_i)$ as follows
\begin{align}
\mathcal{L}^\epsilon(\hat{\bu}^k)(\bx_i) = \int_{H_\epsilon(\bx_i)} (w_1(\by, \bx_i) + w_2(\by, \bx_i))  d\by,
\end{align}
where $w_1$ and $w_2$ can be determined from expression of $\mathcal{L}^\epsilon$ in \autoref{eq:total peri force}. In the simulation we approximate $\mathcal{L}^\epsilon(\hat{\bu}^k)(\bx_i)$ as below
\begin{align}
\mathcal{L}^\epsilon(\hat{\bu}^k)(\bx_i) \approx \sum_{\bx_j \in D_h\cap H_\epsilon(\bx_i)} (w_1(\bx_j, \bx_i) +w_w(\bx_j, \bx_i)) V_j \bar{V}_{ij},
\end{align}
where $V_j = h^2$ for uniform mesh in $2$-d and $\bar{V}_{ij}$ is the volume correction.  

The numerical results are presented in the following section.


\begin{table}
\centering

\addtolength{\tabcolsep}{+2pt}

\renewcommand{\arraystretch}{1.3}

\begin{tabular}{|l||r|r|}  
\hline
Parameters $\setminus$ Poisson's ratio & $\nu = 0.22$ & $\nu = 0.245$\\
\hline\hline
$c$ & $4712.4$ & $4712.4$\\  
\hline
$\bar{C}$ & $-1.0623\times 10^{12}$ & $-1.7349\times 10^{11}$\\
\hline
$\beta$ & $1.7533\times 10^{8}$ & $1.5647\times 10^{8}$\\
\hline
$\bar{r}$ & $5.3402\times 10^{-5}$ & $5.6529\times 10^{-5}$\\
\hline
\end{tabular}

\caption{Peridynamic material parameters assuming bulk modulus $K = 25 \, GPa$ and critical energy release rate $G_c = 500\, J/m^{-2}$. Density is $\rho = 1200\, kg/m^3$.}
\label{tab:mat props}
\end{table}

\begin{figure}
\centering
\includegraphics[scale=0.15]{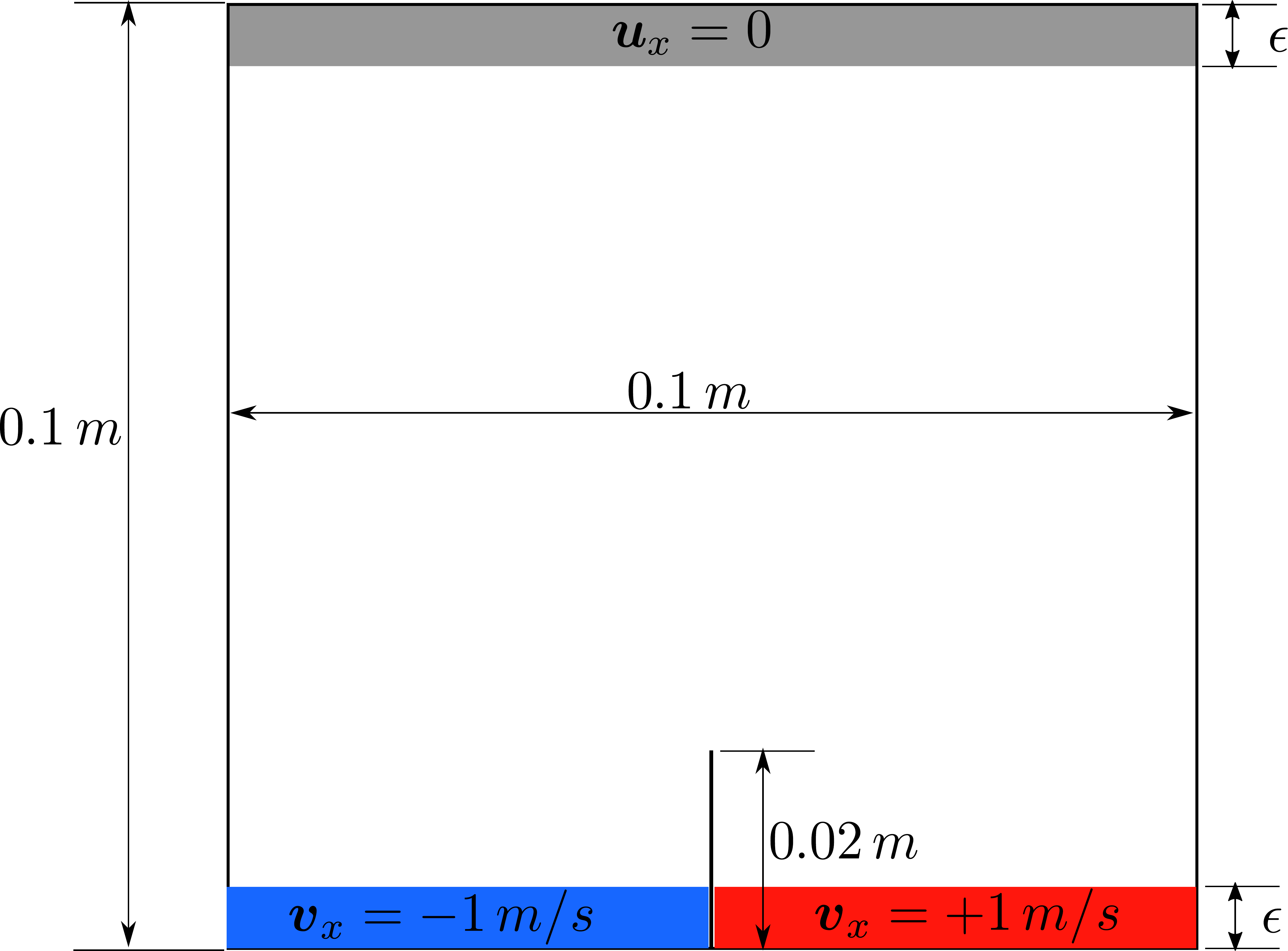}
\caption{Material domain $D = [0,0.1\, m]^2$ with crack of length $0.02\,m$. The x-component of displacement is fixed along a collar of thickness equal to the horizon on top. On the bottom the velocity $\bv_x = \pm 1\, m/s$ along x-direction is specified on either side of the crack to make the crack propagate upwards.}\label{fig:setup crack prop}
\end{figure}

\begin{figure}
    \centering
    \begin{subfigure}{.47\linewidth}
        \includegraphics[scale=0.22]{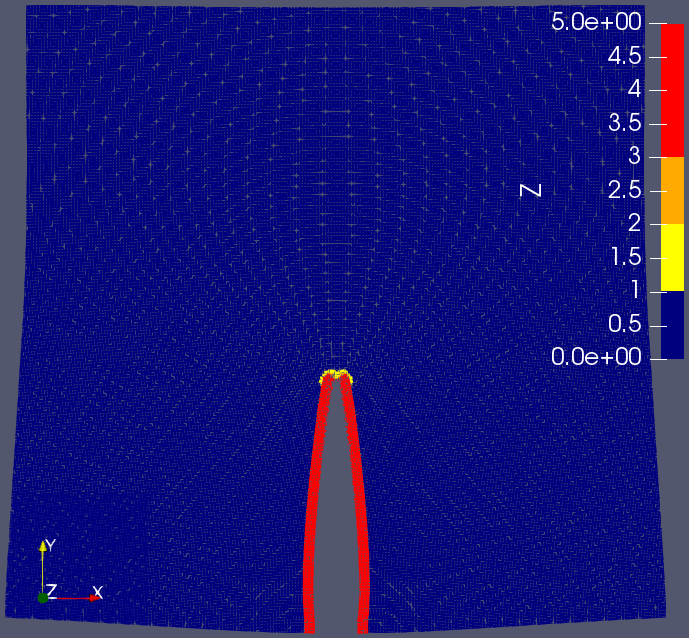}
        \caption{}
    \end{subfigure}
    \begin{subfigure}{.47\linewidth}
	    \includegraphics[scale=0.14]{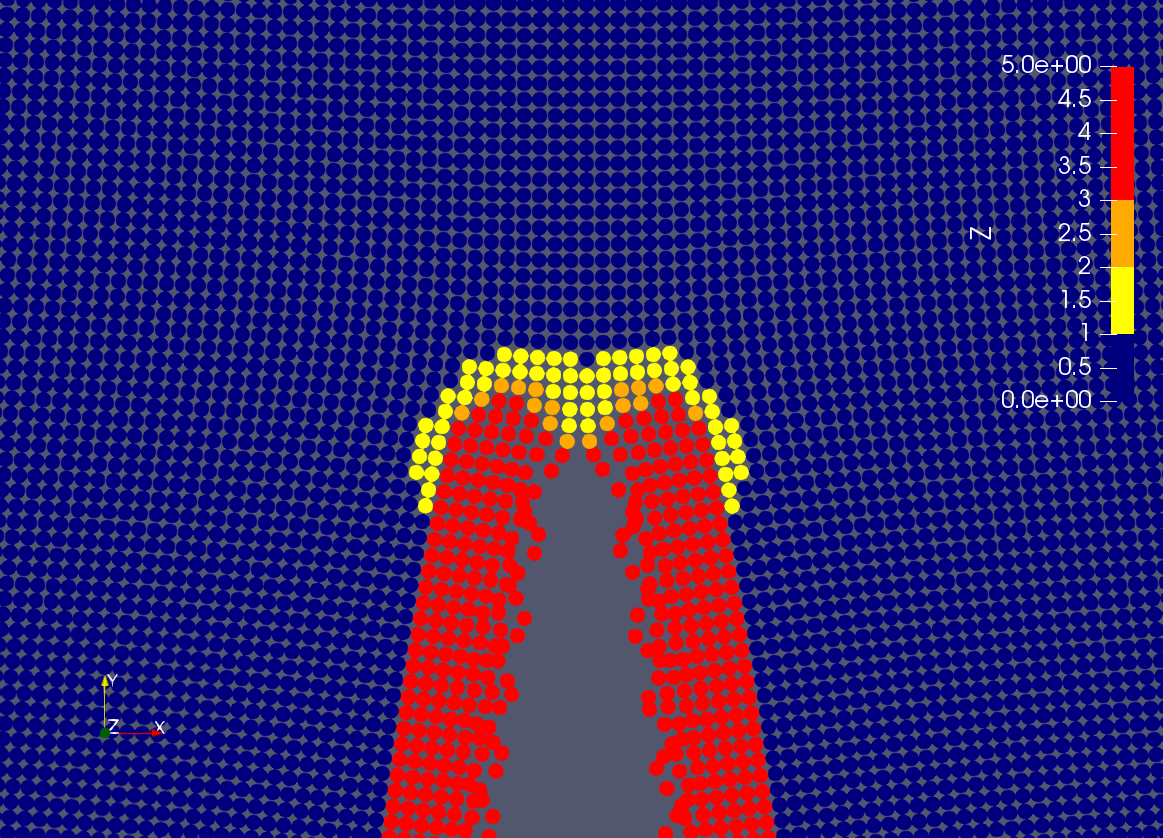}
        \caption{}
    \end{subfigure}
    
    \caption{(a) Color plot of damage function $Z$ on deformed material domain at time $t=34\,\mu s$. Dark blue represents undamaged material $Z<1$, $Z\approx 1$ is yellow at crack tip, red is softening material. The plot is for a horizon $\epsilon = 2\, mm$ and $h = \epsilon/8$. Here, the displacements are scaled by $100$ and damage function is cut off at $5$ to highlight the crack zone. The maximum displacement is $4.4 \, mm$ and the maximum value of $Z(x)$ is $82$ at $t=34\,\mu s$. (b) View near the crack tip.}
    \label{fig:damage plot}
\end{figure}

\subsection{Crack propagation: Fracture energy and numerical convergence study}\label{ss:crack prop}
The problem is intentionally similar to the problem given in the simulation presented in\cite{CMPer-Lipton2}. We consider a 2-d domain $D=[0,0.1\, m]^2$ (with unit thickness in third direction) with vertical crack of length $0.02\, m$. Boundary conditions are described in \autoref{fig:setup crack prop}. 
The simulation time is $T = 34\, \mu s$ and the time step is $\Delta t = 0.004\,\mu s$. 

We run simulations for four different horizons $\epsilon = 8\, mm, 4\, mm, 2\, mm, 1\, mm$. For each horizon, we obtain the results for mesh sizes $h = \epsilon/2, \epsilon/4, \epsilon/8$. We take uniform square mesh of size $h$. 
Material properties correspond to the Poisson's ration $\nu = 0.245$, see \autoref{tab:mat props}.


For the coarsest horizon $\epsilon = 8\, mm$, number of mesh nodes are (approximately) $0.9\times 10^3, 3.5\times 10^3 , 13.7\times 10^3$ for $h=4, 2, 1\, mm$ respectively. The memory consumed are $10$ MB, $16$ MB, $95$ MB respectively. For $\epsilon = 1\, mm$, number of nodes are $42\times 10^3, 167\times 10^3, 668\times 10^3$ for $h=0.5, 0.25, 0.125\, mm$ respectively. The memory consumed are $44$ MB, $370$ MB, $4400$ MB respectively. All computations were performed on a single workstation in parallel using $20$ threads.
%

\subsubsection{Fracture energy of crack zone}

The extent of damage at material point $\bx$ is given by the function $Z(\bx)$
\begin{align}\label{eq:damage}
Z(\bx) &= \max_{\by \in H_\epsilon(\bx) \cap D} \frac{S(\by,\bx;\bu)}{S_c}.
\end{align}
We define the crack zone as set of material points which have $Z > 1$. We compute the peridynamic energy of crack zone and compare it with the Griffith's fracture energy. For a crack of length $l$, the Griffith's fracture energy (G.E.) will be $G.E. = G_c \times l$. The peridynamic fracture energy (P.E.) is given by 
\begin{align*}
P.E. = \int_{\substack{\bx \in D,\\ Z(\bx) \geq 1}} \left[ \frac{1}{\epsilon^d \omega_d} \int_{H_\epsilon(\bx)} |\by - \bx| \mathcal{W}^\epsilon(S(\by,\bx,\bu))\,d\by \right] d\bx,
\end{align*}
where $\mathcal{W}^\epsilon(S(\by,\bx,\bu))$ is the bond-based potential, see \autoref{tensilepot}. For the choice of $f(r) $ and $g(r)$, only bond-based potential $f$ contributes to the fracture energy, therefore $P.E.$ is computed only from bond-based interaction. 

\begin{figure}[h]
    \centering
    \begin{subfigure}{.47\linewidth}
        \includegraphics[scale=0.35]{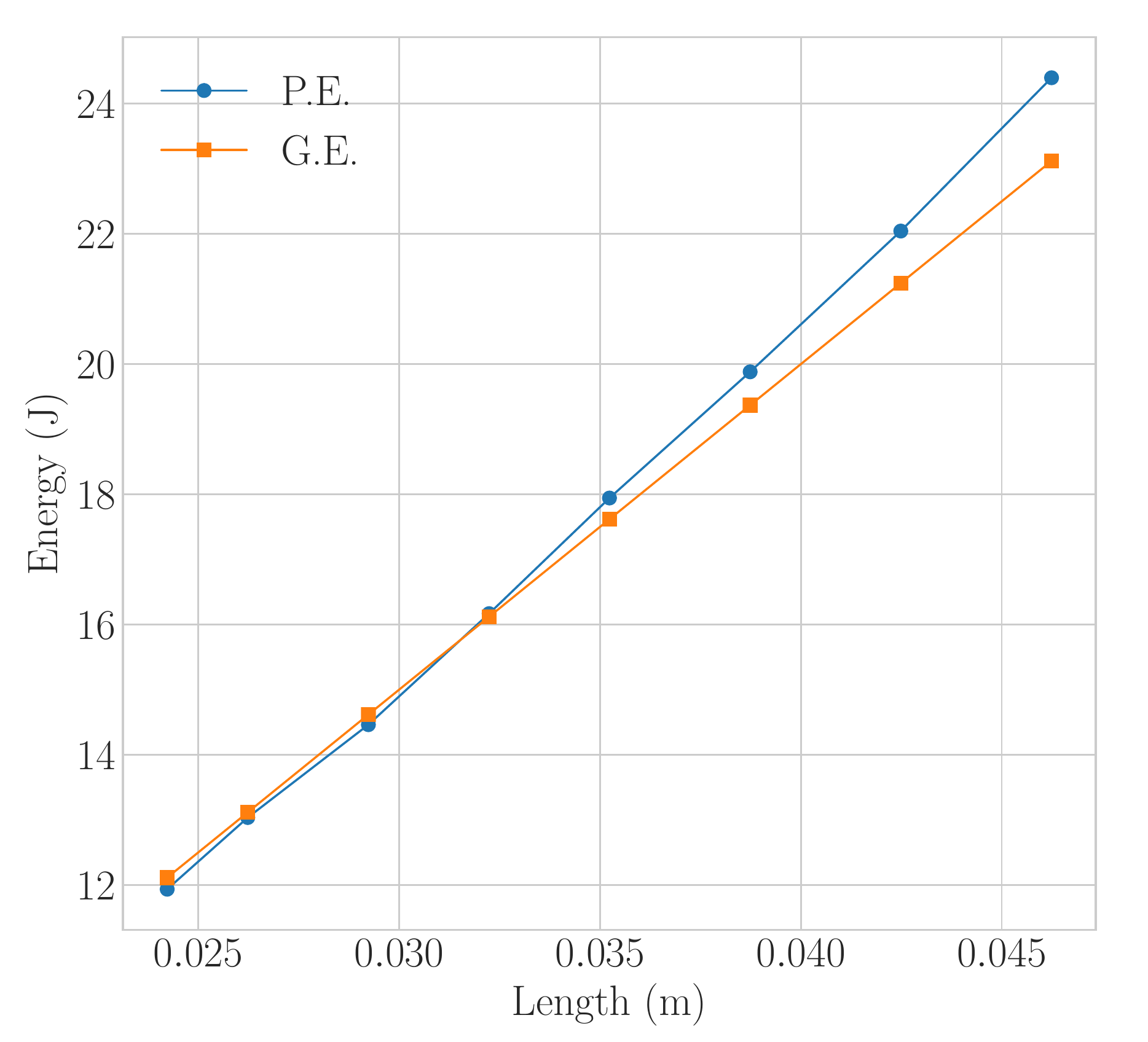}
        \caption{Horizon $\epsilon = 2\, mm$}
    \end{subfigure}
    \begin{subfigure}{.47\linewidth}
	    \includegraphics[scale=0.35]{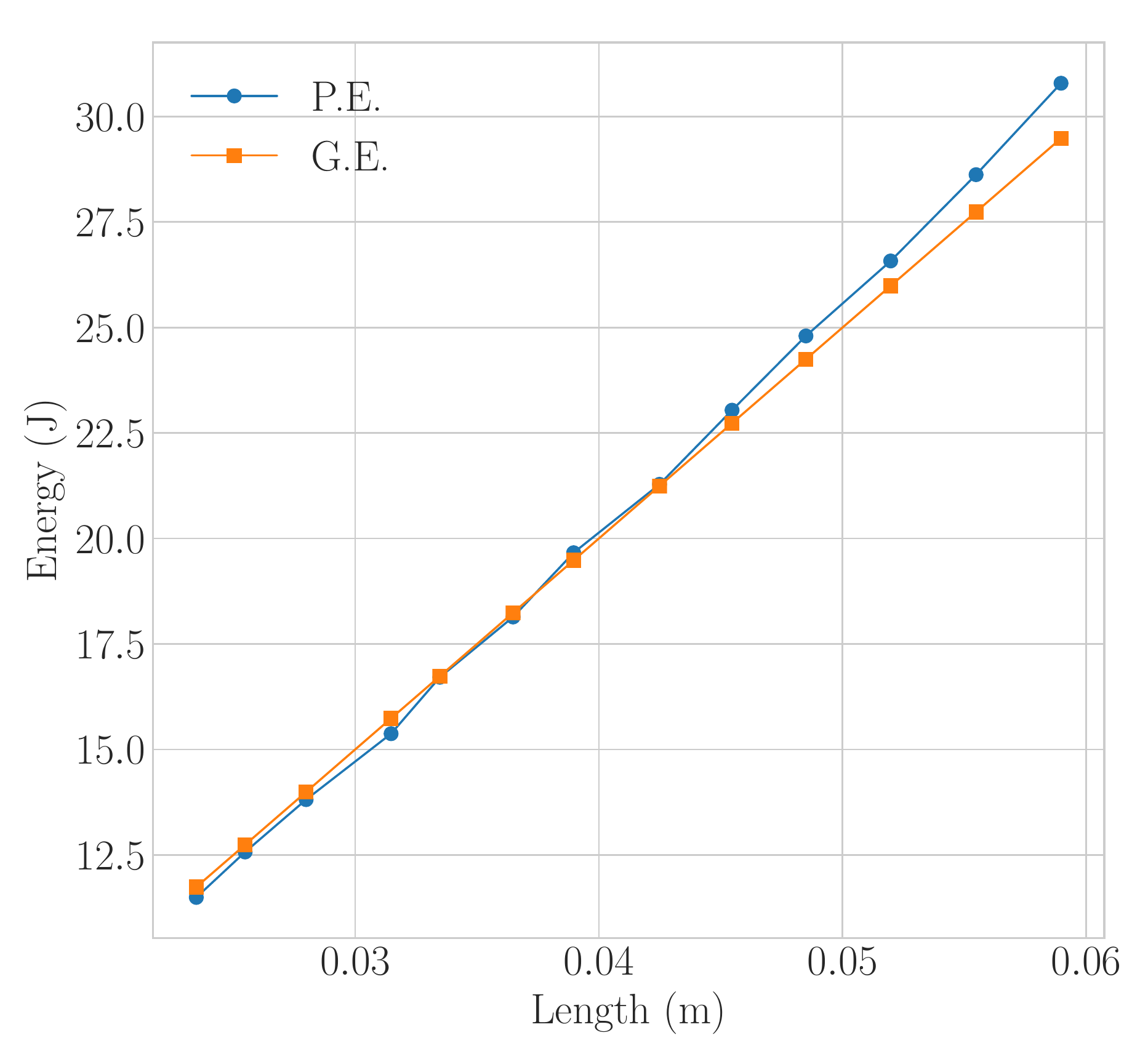}
        \caption{Horizon $\epsilon = 4\, mm$}
    \end{subfigure}
    
   \caption{Crack length vs peridynamic fracture energy (P.E.) and Griffith's fracture energy (G.E.). G.E. is simply $G_c \times l$ where $G_c = 500$.}
    \label{fig:crack zone energy}
\end{figure}

\begin{figure}
\centering
\includegraphics[scale=0.3]{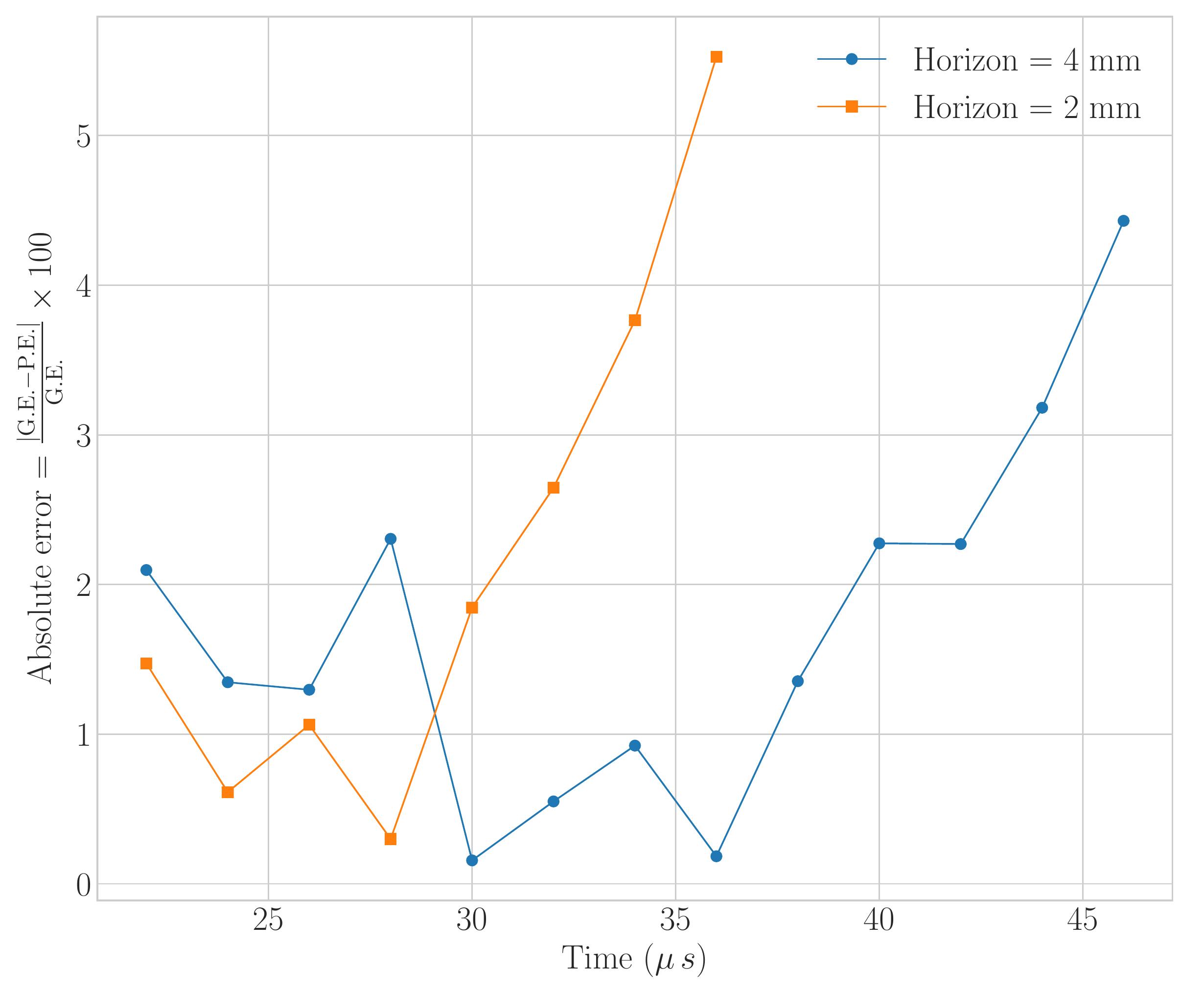}
\caption{Percentage error between peridynamic fracture energy and Griffith's fracture energy.}\label{fig:error crack zone energy}
\end{figure}

\autoref{fig:damage plot} shows the plot of $Z$ at time $t=34\,\mu s$ for horizon $\epsilon = 2\, mm$. The figure on the right shows the $Z$ field near a crack tip. 
In \autoref{fig:crack zone energy} we plot the peridynamic and Griffith's fracture energy as a function of crack length. We see better agreement between the two energies up to larger length of crack for coarse horizon. In \autoref{fig:error crack zone energy} we plot the error in fracture energy at different times. 

\begin{figure}[h]
\centering
\includegraphics[scale=0.35]{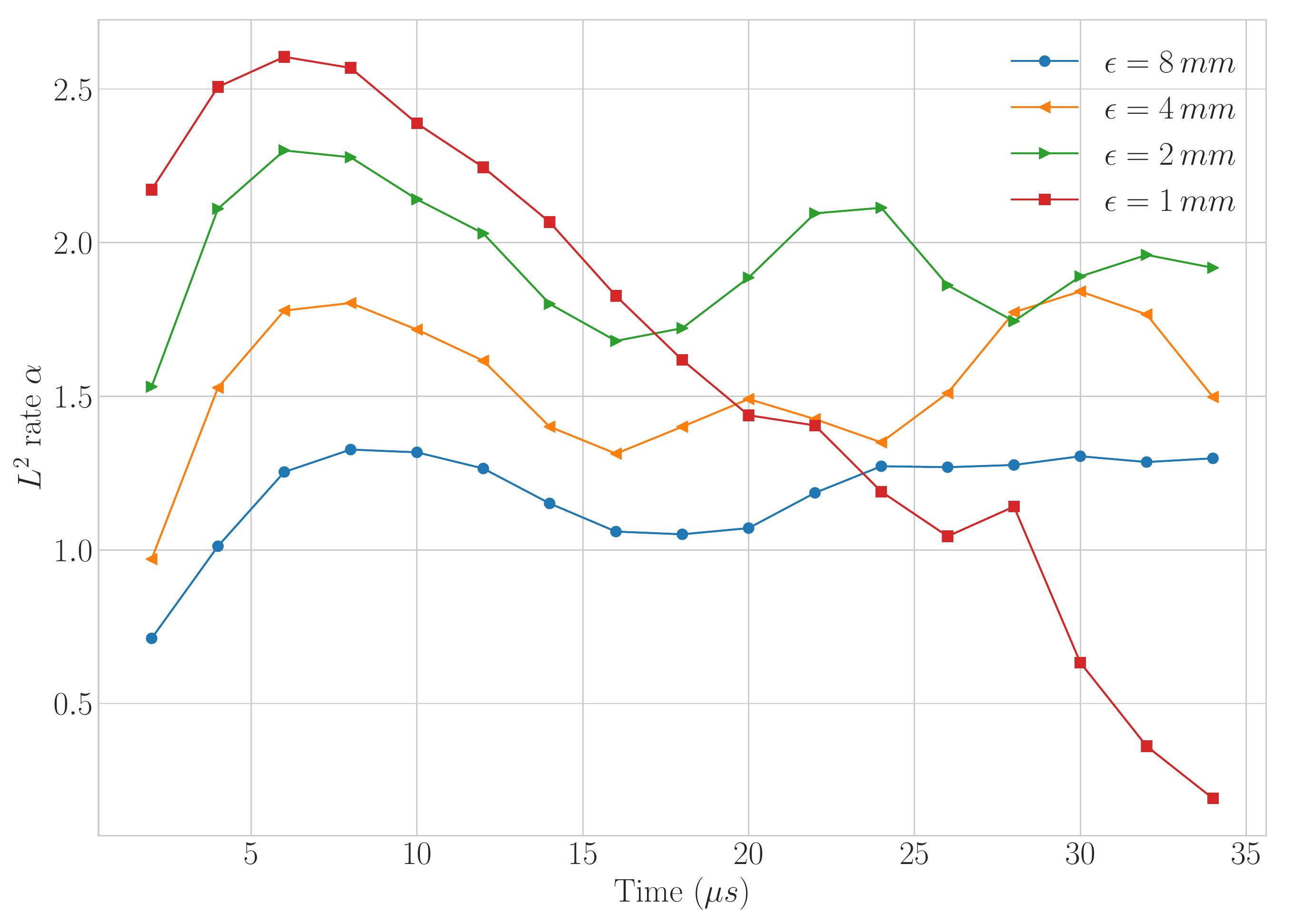}
\caption{Convergence rate with respect to mesh size for different fixed size of horizons.}\label{fig:rate}
\end{figure}

\subsubsection{Convergence rate}
Consider a fixed horizon $\epsilon$ and three different mesh sizes $h = \epsilon/2, \epsilon/4, \epsilon/8$. We compute the convergence rate as follows. Let $\bu_1,\bu_2,\bu_3$ be approximate solutions corresponding to meshes of size  $h_1,h_2,h_3$, and let $\bu$ be the exact solution. We write the error as $||\bu_h - \bu|| =C h^\alpha$ for some constant $C$ and $\alpha>0$, and fix the ratio of mesh size  $h_1/h_2 = h_2/h_3 = r$, to get
\begin{align*}
\log( ||\bu_1 - \bu_2||) &= C + \alpha \log h_2, \\
\log( ||\bu_2 - \bu_3||) &= C + \alpha \log h_3.
\end{align*}
Recall that the norm $||\cdot||$ is $L^2$ norm. From above two equations, it is easy to see that the rate of convergence $\alpha$ is 
\begin{align}\label{eq:rate formula}
\dfrac{\log( ||\bu_1 - \bu_2||) - \log( ||\bu_2 - \bu_3||)}{\log(r)}.
\end{align}

%

The convergence result for four different horizons is shown in \autoref{fig:rate}. From this figure we see that for $\epsilon=1\, mm$ the convergence rate is greater than $1$ for simulation times below $25\mu s$.
For all other horizons the rate is larger than $1$ for simulation times below $35\mu s$. Here the discrepancy is due to the error accumulation at each time step and can be reduced some what by taking smaller time steps. 
The simulations show a rate of convergence that agrees with the a-priori estimates given in \autoref{thm:convergence}.

\subsection{Bending test with pre-crack}\label{ss:bending}
\begin{figure}
\centering
\includegraphics[scale=0.1]{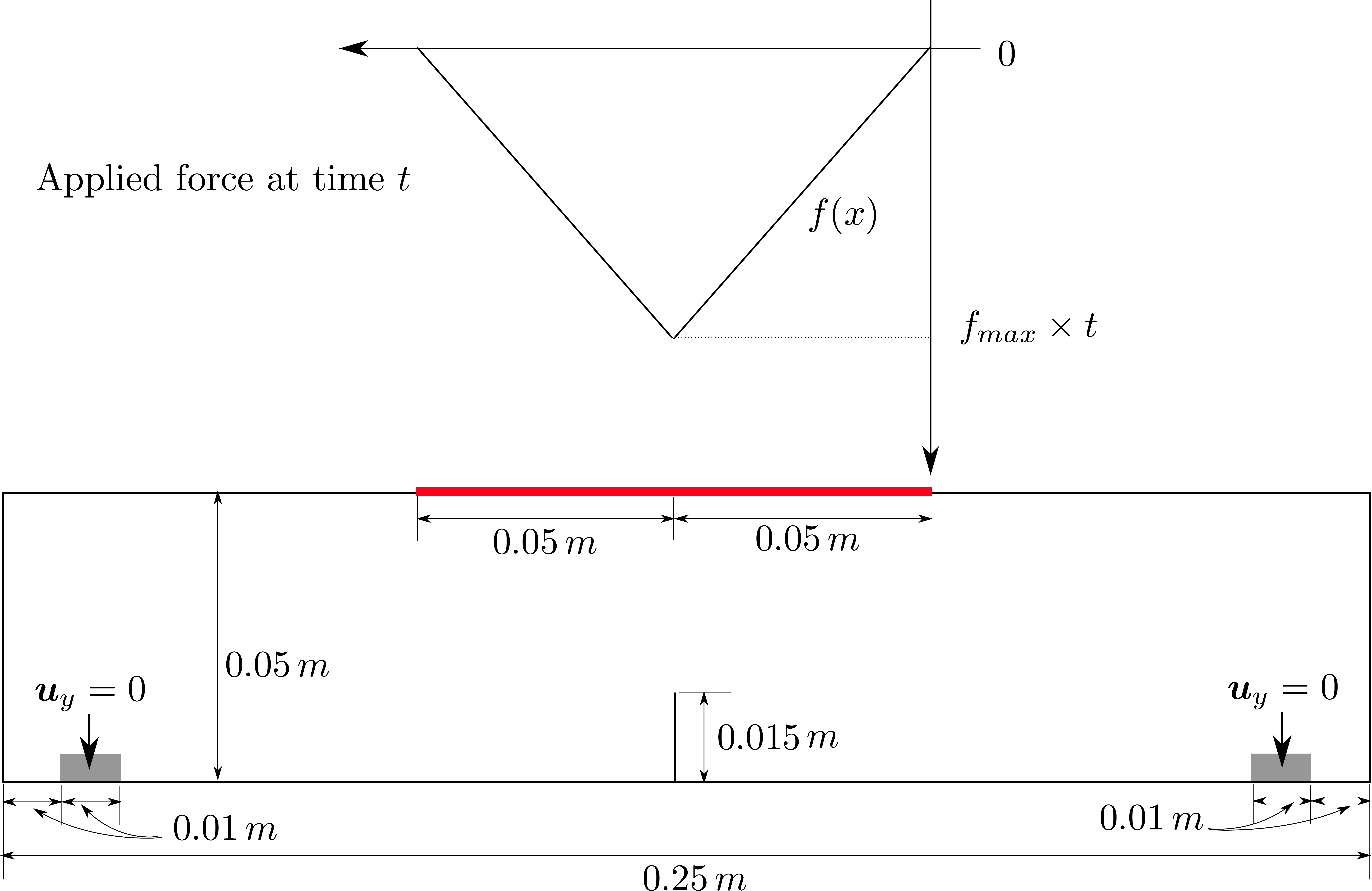}
\caption{Material domain $D = [0,0.25\, m]\times [0,0.05\,m]$ with single verticle crack of length $0.015\,m$ at mid point of bottom edge. We apply linear in time distributed load, along negative y-direction, on part of the top edge. At any time $t$, the load is zero at the end points of loading line (red line) and is $f_{max} \times t$ at the midpoint. We take constant $f_{max} = -1.0\times 10^{13}$. We fix a vertical displacement on two support regions shown in the figure.}\label{fig:setup bending}
\end{figure}

We consider a 2-d material domain (with unit thickness in third direction) $D = [0,0.25\, m]\times [0,0.05\,m]$ with single and double vertical cracks. We fix horizon to $\epsilon = 0.010\, m$ and mesh size $h = \epsilon/4 \, mm$. The boundary conditions are described in \autoref{fig:setup bending} for single crack. For the double crack problem, the two vertical cracks are symmetrically located at distance $0.02\, m$ along x-axis from the mid point $x=0.125 \, m, y=0$. With time step $\Delta t = 0.0014\,\mu s$ we run simulations upto time $T = 350 \, \mu s$.  Material properties correspond to the Poisson's ration $\nu = 0.22$, see \autoref{tab:mat props}.

\begin{figure}
    \centering
    \begin{subfigure}{.47\linewidth}
        \scalebox{1.0}[1.25]{\includegraphics[scale=0.2]{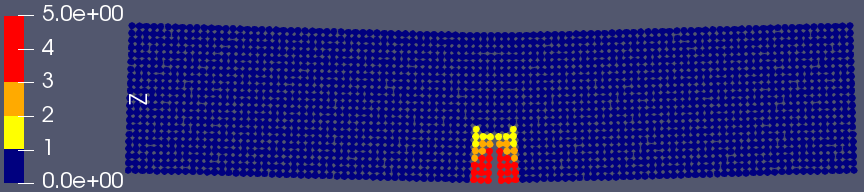}}
        \caption{$t = 130\, \mu s$}
    \end{subfigure}
    \begin{subfigure}{.47\linewidth}
        \scalebox{1.0}[1.25]{\includegraphics[scale=0.2]{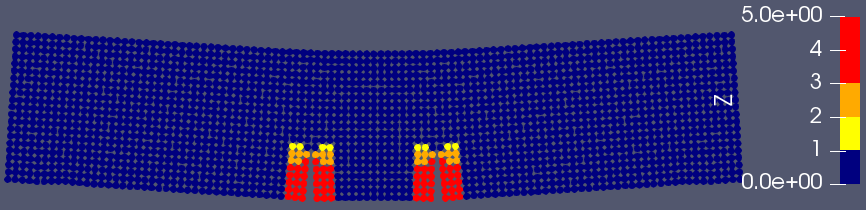}}
        \caption{$t = 180\, \mu s$}
    \end{subfigure}
    
    \begin{subfigure}{.47\linewidth}
        \scalebox{1.0}[1.25]{\includegraphics[scale=0.2]{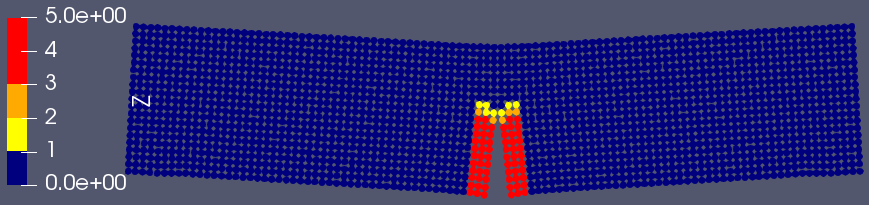}}
        \caption{$t = 180\, \mu s$}
    \end{subfigure}
    \begin{subfigure}{.47\linewidth}
	    \scalebox{1.0}[1.25]{\includegraphics[scale=0.2]{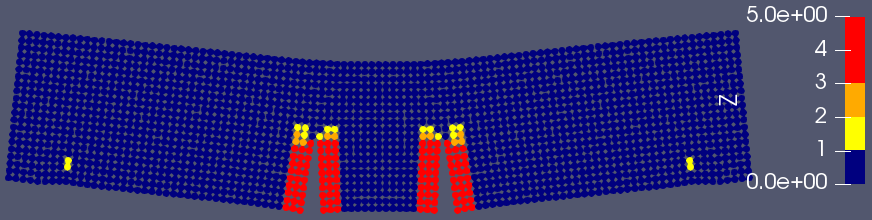}}
        \caption{$t = 220\, \mu s$}
    \end{subfigure}
    
    \begin{subfigure}{.47\linewidth}
        \scalebox{1.0}[1.25]{\includegraphics[scale=0.2]{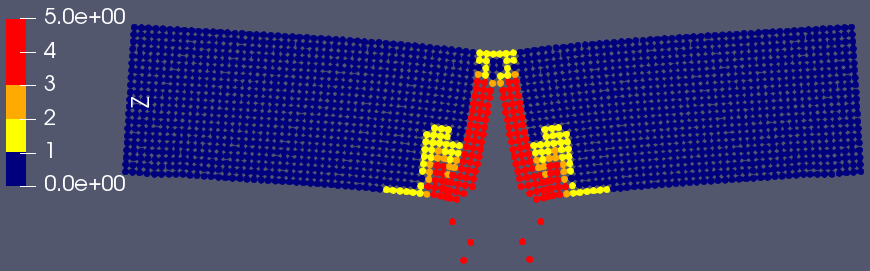}}
        \caption{$t = 190\, \mu s$}
    \end{subfigure}
    \begin{subfigure}{.47\linewidth}
	    \scalebox{1.0}[1.25]{\includegraphics[scale=0.2]{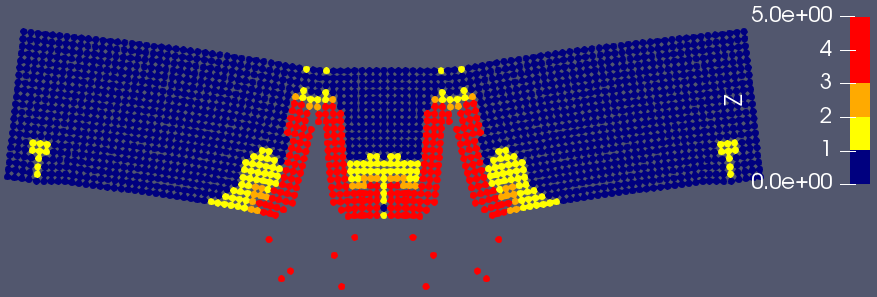}}
        \caption{$t = 240\, \mu s$}
    \end{subfigure}
    
   \caption{Damage profile under bending load. Plots on left are for single crack and plots on right are for double crack.}
    \label{fig:damage plot bending}
\end{figure}

\begin{figure}
    \centering
    \begin{subfigure}{.46\linewidth}
        \includegraphics[scale=0.35]{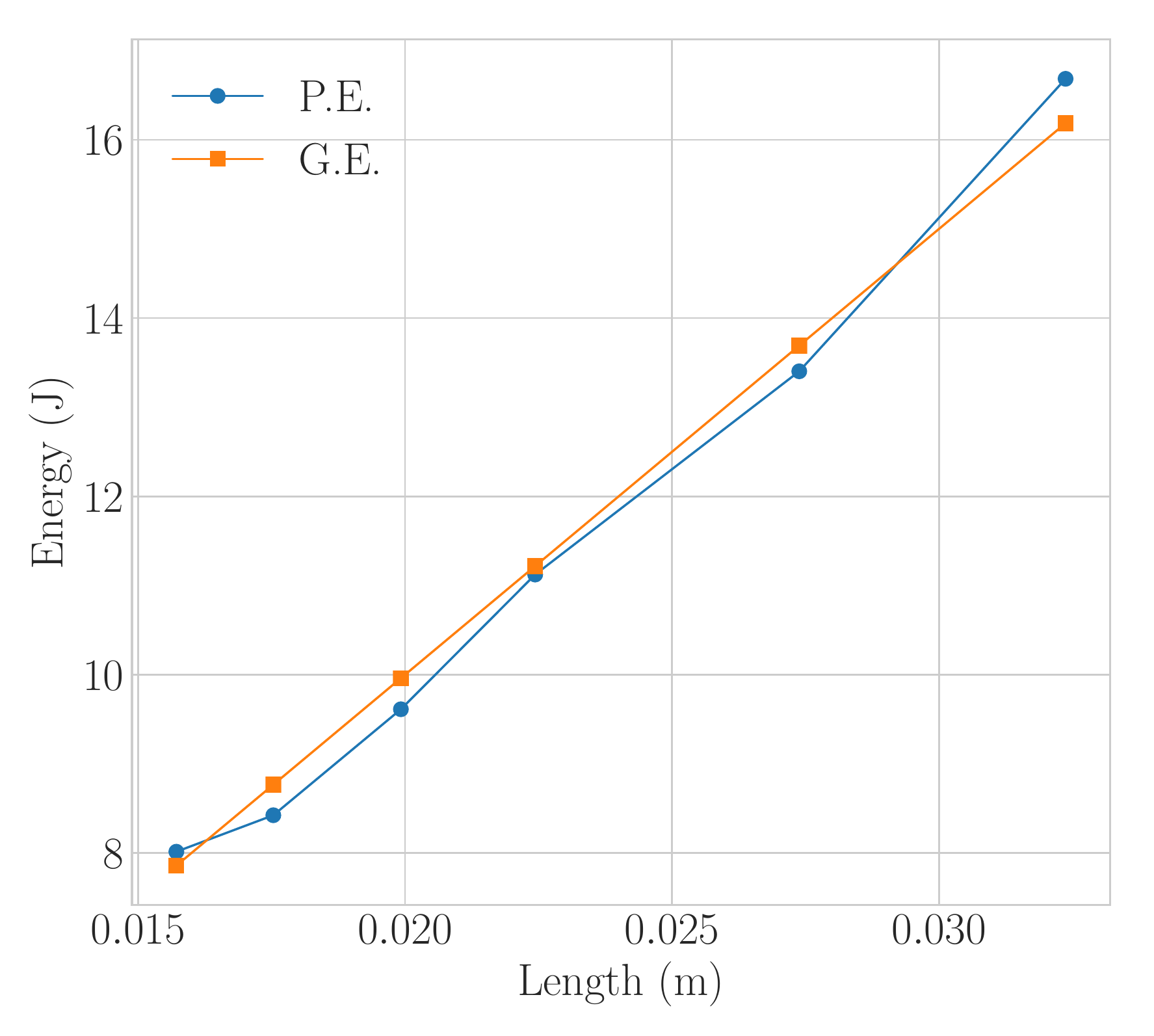}
        \caption{Single crack}
    \end{subfigure}
    \begin{subfigure}{.46\linewidth}
	    \includegraphics[scale=0.35]{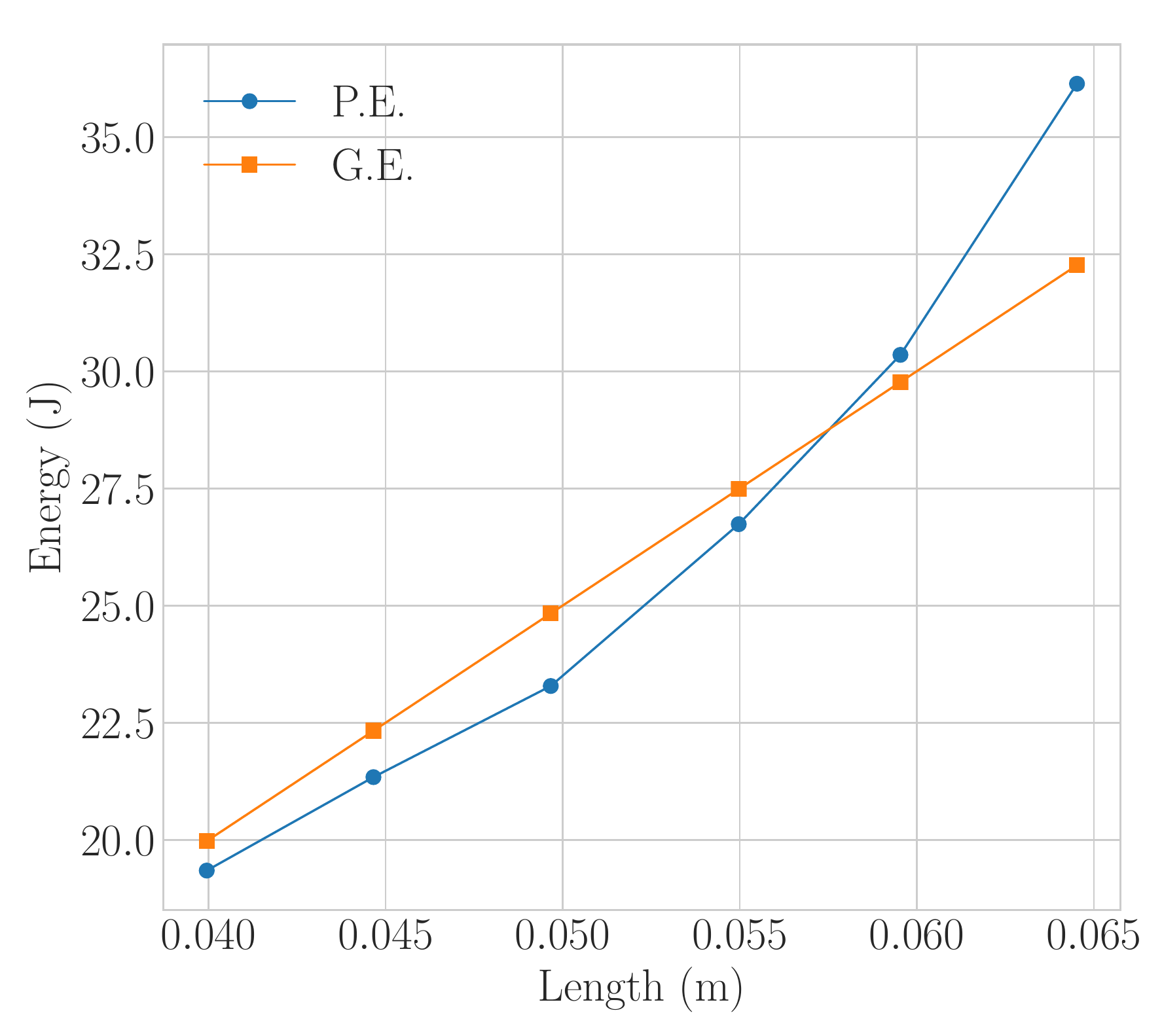}
        \caption{Double crack}
    \end{subfigure}
    
   \caption{Crack length vs peridynamic fracture energy (P.E.) and Griffith's fracture energy (G.E.).}
    \label{fig:crack zone energy bending}
\end{figure}

In \autoref{fig:damage plot bending} damage profile at various times are shown for both single and double crack problem. In \autoref{fig:crack zone energy bending} we plot the fracture energy as a function of total crack length. The error in energy remain below $5\%$ till $185 \, \mu s$ for single crack problem and $232\, \mu s$ for double crack problem. As we can see from \autoref{fig:damage plot bending}, after time $185\, \mu s$ for single crack and $232\, \mu s$ for double crack, the spread of damage around crack is higher and therefore peridynamic fracture energy is higher.

\section{Proof of Lipschitz continuity for the non-local force}\label{s:proofs}
In this section, we prove \sautoref{Proposition}{prop:lipschitz} and \sautoref{Proposition}{prop:lipschitz L2}. 

\subsection{Proof of Proposition 1}\label{ss:proof prop 2}
Recall that $I = [0,T]$ is the time domain, $X=\Cholderz{\perdrd}\times\Cholderz{\perdrd}$, and $F^\epsilon(y,t) = (F^\epsilon_1(y,t), F^\epsilon_2(y,t))$, where $F^\epsilon_1(y,t) = y^2$ and $F^\epsilon_2(y,t) = \mathcal{L}^\epsilon(y^1) + \bb(t)$. Given $t\in I$ and $y=(y^1,y^2), z=(z^1, z^2)\in X$, we have
\begin{align}\label{eq:norm of F diff in X}
&\normX{F^\epsilon(y,t) - F^\epsilon(z,t)}{X} \notag \\
&= \Choldernorm{y^2 - z^2}{\perdrd} + \Choldernorm{\mathcal{L}^{\epsilon}(y^1) - \mathcal{L}^\epsilon(z^1)}{\perdrd}
\end{align}
and
\begin{align}\label{eq:norm of F in X}
\normX{F^\epsilon(y,t)}{X} = \Choldernorm{y^2}{\perdrd} + \Choldernorm{\mathcal{L}^{\epsilon}(y^1)}{\perdrd} + b,
\end{align}
where $b = \sup_t ||\bb(t)||_{\Cholder{}}$. 

Thus, to prove \autoref{eq:lipschitz property of F} and \autoref{eq:bound on F} of  \sautoref{Proposition}{prop:lipschitz} we need to study the terms associated with $\mathcal{L}^\epsilon$ in the equations listed above. The peridynamic force $\mathcal{L}^\epsilon$ is sum of two forces, the tensile force $\mathcal{L}^\epsilon_T$ and the dilatational force $\mathcal{L}^\epsilon_D$. So for $\bu,\bv \in \Cholderz{\perdrd}$ we have
\begin{align}\label{eq:lipschitz norm of per force diff}
&\Choldernorm{\mathcal{L}^{\epsilon}(\bu) - \mathcal{L}^{\epsilon}(\bv)}{\perdrd} \notag \\
&\leq \Choldernorm{\mathcal{L}^{\epsilon}_T(\bu) - \mathcal{L}^{\epsilon}_T(\bv)}{\perdrd} + \Choldernorm{\mathcal{L}^{\epsilon}_D(\bu) - \mathcal{L}^{\epsilon}_D(\bv)}{\perdrd}
\end{align}
and
\begin{align}\label{eq:lipschitz norm of per force}
\Choldernorm{\mathcal{L}^{\epsilon}(\bu)}{\perdrd} \leq \Choldernorm{\mathcal{L}^{\epsilon}_T(\bu)}{\perdrd} + \Choldernorm{\mathcal{L}^{\epsilon}_D(\bu)}{\perdrd}.
\end{align}

We conclude listing estimates that will be used in the sequel.
For $\bu \in \Cholderz{\perdrd}$ and $\omega\in \Cholderz{D;[0,1]}$ one easily deduces the estimates 
\begin{align}
|\bu(\bx + \epsilon \bxi) - \bu(\bx)| &\leq (\epsilon|\bxi|)^\gamma ||\bu||_{\Cholder{}},  \notag \\
|\bu(\bx+ \epsilon \bxi) - \bu(\by + \epsilon\bxi)| &\leq |\bx - \by|^\gamma ||\bu||_{\Cholder{}},\notag \\
|\omega(\bx+ \epsilon \bxi) - \omega(\by + \epsilon\bxi)| &\leq |\bx - \by|^\gamma ||\omega||_{\Cholder{}},
\end{align}
for $\bx, \by \in D$ and $\bxi \in H_1(\bzero)$. Since $\bu$ and $\omega$ are extended by zero outside $D$ these estimates also hold for all points outside $D$.

\subsubsection{Lipschitz continuity in H\"older space}
In this subsection, we provide upper bounds on \autoref{eq:lipschitz norm of per force diff}.

\paragraph{Non-local tensile force}\label{sss:bond-based proof}
For any $\bu,\bv \in \Cholderz{\perdrd}$, we provide upper bounds on
\begin{align}\label{eq:lipshitz norm bond force}
&\Choldernorm{\mathcal{L}^{\epsilon}_T(\bu) - \mathcal{L}^{\epsilon}_T(\bv)}{\perdrd} \notag \\
&= \sup_{\bx\in D} |\mathcal{L}^{\epsilon}_T(\bu)(\bx) - \mathcal{L}^{\epsilon}_T(\bv)(\bx)| \notag \\
&+ \sup_{\bx, \by \in D, \bx \neq \by} \frac{|(\mathcal{L}^{\epsilon}_T(\bu)(\bx) - \mathcal{L}^{\epsilon}_T(\bv)(\bx)) - (\mathcal{L}^{\epsilon}_T(\bu)(\by) - \mathcal{L}^{\epsilon}_T(\bv)(\by))|}{|\bx - \by|^\gamma}.
\end{align}
Applying \autoref{eq:per bondbased simple force expression} and proceeding as in section
\autoref{sss:consistency} we see that
\begin{align}\label{eq:est 0}
& |\mathcal{L}^\epsilon_T(\bu)(\bx) - \mathcal{L}^\epsilon_T(\bv)(\bx) | \notag \\
&= \left\vert \frac{2}{\epsilon \omega_d} \int_{H_1(\bzero)} \omega_\bxi(\bx) \frac{J(|\bxi|)}{\sqrt{s_\bxi}} \left[ f'(\bubar_\bxi(\bx)\cdot \be_\bxi/\sqrt{s_\bxi}) - f'(\bvbar_\bxi(\bx)\cdot \be_\bxi/\sqrt{s_\bxi}) \right] \be_{\bxi} d\bxi \right\vert \notag \\
&\leq \frac{2}{\epsilon \omega_d} \int_{H_1(\bzero)} \omega_\bxi(\bx) \frac{J(|\bxi|)}{\sqrt{s_\bxi}} \left\vert f'(\bubar_\bxi(\bx)\cdot \be_\bxi/\sqrt{s_\bxi}) - f'(\bvbar_\bxi(\bx)\cdot \be_\bxi/\sqrt{s_\bxi}) \right\vert  d\bxi \notag \\
&\leq \frac{2C^f_2}{\epsilon \omega_d} \int_{H_1(\bzero)} \frac{J(|\bxi|)}{s_\bxi} \left\vert \bubar_\bxi(\bx) - \bvbar_\bxi(\bx)\right\vert  d\bxi.
\end{align}
A straightforward calculation gives the estimate
\begin{align*}
&|\bu_\bxi(\bx) - \bv_\bxi(\bx)| = |\bu(\bx + \epsilon \bxi) - \bu(\bx) - (\bv(\bx+ \epsilon\bxi) - \bv(\bx))|  \notag \\
&\leq |\bu(\bx + \epsilon \bxi) - \bv(\bx+ \epsilon\bxi)| + |\bu(\bx) - \bv(\bx)| \leq 2 ||\bu - \bv||_{\Cholder{}}
\end{align*}
and on applying this \autoref{eq:est 0} we get
\begin{align}\label{eq:lipshitz bound bondbased 1}
|\mathcal{L}^\epsilon_T(\bu)(\bx) - \mathcal{L}^\epsilon_T(\bv)(\bx) | &\leq \frac{4 C^f_2 \bar{J}_1}{\epsilon^2} ||\bu - \bv||_{\Cholder{}},
\end{align}
where $\bar{J}_1$ is given by \autoref{eq:def bar J}.
Next we derive a bound on
\begin{align*}
&\frac{|(\mathcal{L}^{\epsilon}_T(\bu)(\bx) - \mathcal{L}^{\epsilon}_T(\bv)(\bx)) - (\mathcal{L}^{\epsilon}_T(\bu)(\by) - \mathcal{L}^{\epsilon}_T(\bv)(\by))|}{|\bx - \by|^\gamma} \notag \\
&= \frac{1}{|\bx - \by|^\gamma} \bigg\vert\frac{2}{\epsilon\omega_d} \int_{H_1(\bzero)} \frac{J(|\bxi|)}{\sqrt{s_\bxi}} \left[\omega_\bxi(\bx) ( f'(\bubar_\bxi(\bx)\cdot \be_\bxi/\sqrt{s_\bxi}) - f'(\bvbar_\bxi(\bx)\cdot \be_\bxi/\sqrt{s_\bxi}) ) \right. \notag \\
&\qquad \left. - \omega_\bxi(\by) ( f'(\bubar_\bxi(\by)\cdot \be_\bxi/\sqrt{s_\bxi}) - f'(\bvbar_\bxi(\by)\cdot \be_\bxi/\sqrt{s_\bxi}) )\right] \be_{\bxi} d\bxi \bigg\vert.
\end{align*}
Let
\begin{align}\label{eq:def H}
H &:= \frac{1}{|\bx - \by|^\gamma} \bigg\vert \omega_\bxi(\bx) ( f'(\bubar_\bxi(\bx)\cdot \be_\bxi/\sqrt{s_\bxi}) - f'(\bvbar_\bxi(\bx)\cdot \be_\bxi/\sqrt{s_\bxi}) ) \notag \\
& \qquad - \omega_\bxi(\by) ( f'(\bubar_\bxi(\by)\cdot \be_\bxi/\sqrt{s_\bxi}) - f'(\bvbar_\bxi(\by)\cdot \be_\bxi/\sqrt{s_\bxi}) )\bigg\vert.
\end{align}
Then
\begin{align}\label{eq:bondbased ineq 1}
&\frac{|(\mathcal{L}^{\epsilon}_T(\bu)(\bx) - \mathcal{L}^{\epsilon}_T(\bv)(\bx)) - (\mathcal{L}^{\epsilon}_T(\bu)(\by) - \mathcal{L}^{\epsilon}_T(\bv)(\by))|}{|\bx - \by|^\gamma} \notag \\
&\leq \frac{2}{\epsilon\omega_d} \int_{H_1(\bzero)} \frac{J(|\bxi|)}{\sqrt{s_\bxi}} H d\bxi.
\end{align}
To analyze $H$ we consider the function $\br: [0,1] \times D \to \bbR^d$ given by
\begin{align}
\br(l,\bx) &:= \bvbar_\bxi(\bx) + l(\bubar_\bxi(\bx) - \bvbar_\bxi(\bx)),
\end{align}
and $\partial \br(l,\bx)/\partial l = \bubar_\bxi(\bx) - \bvbar_\bxi(\bx)$. We write
\begin{align}\label{eq:est 1}
& f'(\bubar_\bxi(\bx)\cdot \be_\bxi/\sqrt{s_\bxi}) - f'(\bvbar_\bxi(\bx)\cdot \be_\bxi/\sqrt{s_\bxi}) \notag \\
&= \int_0^1 \frac{\partial f'(\br(l,\bx)\cdot\be_\bxi/\sqrt{s_\bxi})}{\partial l} dl \notag \\
&=\int_0^1 \frac{\partial f'(\br\cdot\be_\bxi/\sqrt{s_\bxi})}{\partial \br} \bigg\vert_{\br = \br(l,\bx)} \cdot \frac{\partial \br(l,\bx)}{\partial l} dl  \notag \\
&=\int_0^1 f''(\br(l,\bx)\cdot\be_\bxi/\sqrt{s_\bxi}) \frac{\be_\bxi}{\sqrt{s_\bxi}} \cdot (\bubar_\bxi(\bx) - \bvbar_\bxi(\bx)) dl.
\end{align}
and similarly we have
\begin{align}\label{eq:est 2}
&f'(\bubar_\bxi(\by)\cdot \be_\bxi/\sqrt{s_\bxi}) - f'(\bvbar_\bxi(\by)\cdot \be_\bxi/\sqrt{s_\bxi}) \notag \\
&= \int_0^1 f''(\br(l,\by)\cdot\be_\bxi/\sqrt{s_\bxi}) \frac{\be_\bxi}{\sqrt{s_\bxi}} \cdot (\bubar_\bxi(\by) - \bvbar_\bxi(\by)) dl.
\end{align}
Substituting \autoref{eq:est 1} and \autoref{eq:est 2} into \autoref{eq:def H} gives
\begin{align*}
H &= \frac{1}{|\bx - \by|^\gamma} \bigg\vert \int_0^1 \left[ \omega_\bxi(\bx) f''(\br(l,\bx)\cdot \be_\bxi/\sqrt{s_\bxi})(\bubar_\bxi(\bx) - \bvbar_\bxi(\bx)) \right. \notag \\
& \qquad \left. - \omega_\bxi(\by) f''(\br(l,\by)\cdot \be_\bxi/\sqrt{s_\bxi})(\bubar_\bxi(\by) - \bvbar_\bxi(\by)) \right] \cdot \frac{\be_\bxi}{\sqrt{s_\bxi}} d\bxi \bigg\vert. \notag \\
&\leq \frac{1}{|\bx - \by|^\gamma} \frac{1}{\sqrt{s_\bxi}} \int_0^1 \bigg\vert \omega_\bxi(\bx) f''(\br(l,\bx)\cdot \be_\bxi/\sqrt{s_\bxi})(\bubar_\bxi(\bx) - \bvbar_\bxi(\bx)) \notag \\
& \qquad - \omega_\bxi(\by) f''(\br(l,\by)\cdot \be_\bxi/\sqrt{s_\bxi})(\bubar_\bxi(\by) - \bvbar_\bxi(\by)) \bigg\vert  d\bxi.
\end{align*}
We now add and subtract $\omega_\bxi(\bx) f''(\br(l,\bx)\cdot \be_\bxi/\sqrt{s_\bxi})(\bubar_\bxi(\by) - \bvbar_\bxi(\by))$, and note $0\leq \omega_\bxi \leq 1$, to get
\begin{align}\label{eq:est 3}
H &\leq \frac{1}{|\bx - \by|^\gamma} \frac{1}{\sqrt{s_\bxi}} \int_0^1 |f''(\br(l,\bx) \cdot \be_\bxi/\sqrt{s_\bxi})| |\bubar_\bxi(\bx) - \bvbar_\bxi(\bx) - \bubar_\bxi(\by) + \bvbar_\bxi(\by)| dl \notag \\
&+ \frac{1}{|\bx - \by|^\gamma} \frac{1}{\sqrt{s_\bxi}} \int_0^1 |\omega_\bxi(\bx) f''(\br(l,\bx)\cdot \be_\bxi/\sqrt{s_\bxi}) - \omega_\bxi(\by) f''(\br(l,\by)\cdot \be_\bxi/\sqrt{s_\bxi})| \notag \\
&\qquad \qquad|\bubar_\bxi(\by) - \bvbar_\bxi(\by)|dl \notag \\
&=: H_1 + H_2,
\end{align}
where we denoted first and second term on right hand side as $H_1$ and $H_2$. 
Using the estimate
\begin{align*}
\frac{|\bubar_\bxi(\bx) - \bvbar_\bxi(\bx) - \bubar_\bxi(\by) + \bvbar_\bxi(\by)|}{|\bx - \by|^\gamma} &\leq 2 ||\bu - \bv||_{\Cholder{}}.
\end{align*}
and $|f''(r)| \leq C^f_2$ we see that
\begin{align}\label{eq:est 4}
H_1 &\leq \frac{C^f_2}{|\bx - \by|^\gamma \sqrt{s_\bxi}} \int_0^1 |\bubar_\bxi(\bx) - \bvbar_\bxi(\bx) - \bubar_\bxi(\by) + \bvbar_\bxi(\by)| dl \notag \\
&= \frac{C^f_2}{|\bx - \by|^\gamma \sqrt{s_\bxi}} |\bubar_\bxi(\bx) - \bvbar_\bxi(\bx) - \bubar_\bxi(\by) + \bvbar_\bxi(\by)| \notag \\
&\leq \frac{2C^f_2}{\sqrt{s_\bxi}} ||\bu - \bv||_{\Cholder{}}.
\end{align}

To  bound $H_2$, we add and subtract $\omega_\bxi(\bx) f''(\br(l,\by)\cdot \be_\bxi/\sqrt{s_\bxi})$ and further split the terms
\begin{align}\label{eq:est 5}
H_2 &\leq  \int_0^1  \frac{| f''(\br(l,\bx)\cdot \be_\bxi/\sqrt{s_\bxi}) - f''(\br(l,\by)\cdot \be_\bxi/\sqrt{s_\bxi})|}{|\bx - \by|^\gamma \sqrt{s_\bxi}} |\bubar_\bxi(\by) - \bvbar_\bxi(\by)| dl \notag \\
&+  \int_0^1  \frac{|\omega_\bxi(\bx) - \omega_\bxi(\by)|}{|\bx - \by|^\gamma \sqrt{s_\bxi}} | f''(\br(l,\by)\cdot \be_\bxi/\sqrt{s_\bxi})| |\bubar_\bxi(\by) - \bvbar_\bxi(\by)| dl \notag \\
&=: H_3 + H_4,
\end{align}
where we used the fact that $0\leq\omega_\bxi \leq 1$ in first term. 

We consider $H_3$ first. With $|f'''(r)| \leq C^f_3$ and $0\leq l,1-l \leq 1$ for $l\in [0,1]$, we have
\begin{align*}
&\frac{| f''(\br(l,\bx)\cdot \be_\bxi/\sqrt{s_\bxi}) - f''(\br(l,\by)\cdot \be_\bxi/\sqrt{s_\bxi})|}{|\bx - \by|^\gamma} \notag \\
&\leq \frac{C^f_3}{\sqrt{s_\bxi}} \frac{|\br(l,\bx) - \br(l,\by)|}{|\bx - \by|^\gamma} \notag \\
&\leq \frac{C^f_3}{\sqrt{s_\bxi}} \frac{|1-l||\bvbar_\bxi(\bx) - \bvbar_\bxi(\by)| + |l||\bubar_\bxi(\bx) - \bubar_\bxi(\by)|}{|\bx - \by|^\gamma} \notag \\
&\leq \frac{C^f_3}{\sqrt{s_\bxi}} \left( \frac{|\bvbar_\bxi(\bx) - \bvbar_\bxi(\by)|}{|\bx - \by|^\gamma} + \frac{|\bubar_\bxi(\bx) - \bubar_\bxi(\by)|}{|\bx - \by|^\gamma} \right)
\end{align*}
Following estimates
\begin{align*}
 \frac{|\bvbar_\bxi(\bx) - \bvbar_\bxi(\by)|}{|\bx - \by|^\gamma} &\leq 2 ||\bv||_{\Cholder{}},\quad    \frac{|\bubar_\bxi(\bx) - \bubar_\bxi(\by)|}{|\bx - \by|^\gamma} \leq 2 ||\bu||_{\Cholder{}}
\end{align*}
delivers
\begin{align}\label{eq:est 6}
\frac{| f''(\br(l,\bx)\cdot \be_\bxi/\sqrt{s_\bxi}) - f''(\br(l,\by)\cdot \be_\bxi/\sqrt{s_\bxi})|}{|\bx - \by|^\gamma} &\leq \frac{2C^f_3}{\sqrt{s_\bxi}} (||\bu||_{\Cholder{}} + ||\bv||_{\Cholder{}}).
\end{align}
We use the inequality above together with the estimate
\begin{align*}
|\bubar_\bxi(\by) - \bvbar_\bxi(\by)| &\leq 2 s_\bxi^\gamma ||\bu - \bv||_{\Cholder{}}
\end{align*}
to get
\begin{align}\label{eq:est 7}
H_3 &\leq \frac{4C^f_3}{s^{1-\gamma}_\bxi} (||\bu||_{\Cholder{}} + ||\bv||_{\Cholder{}}) ||\bu - \bv||_{\Cholder{}}.
\end{align}

We now consider $H_4$ in \autoref{eq:est 5}. Using $|f''(r)|\leq C^f_2$, $|\bubar_\bxi(\by) - \bvbar_\bxi(\by)|\leq 2 ||\bu - \bv||_{\Cholder{}}$, and the following estimate
\begin{align}
\frac{|\omega_\bxi(\bx) - \omega_\bxi(\by)|}{|\bx - \by|^\gamma} &= \frac{|\omega(\bx + \epsilon \bxi) \omega(\bx) - \omega(\by + \epsilon\bxi) \omega(\by)|}{|\bx - \by|^\gamma} \notag \\
&\leq \frac{|\omega(\bx+\epsilon\bxi)| |\omega(\bx) - \omega(\by)|}{|\bx - \by|^\gamma} + \frac{|\omega(\by)| |\omega(\bx + \epsilon \by) - \omega(\by + \epsilon\bxi)|}{|\bx - \by|^\gamma} \notag \\
&\leq 2 ||\omega||_{\Cholder{}},
\end{align}
we have
\begin{align}\label{eq:est 8}
H_4 &\leq \frac{4 C^f_2 ||\omega||_{\Cholder{}}}{\sqrt{s_\bxi}}||\bu - \bv||_{\Cholder{}}.
\end{align}

Applying the inequalities \autoref{eq:est 7} and \autoref{eq:est 8} to \autoref{eq:est 5} gives
\begin{align}\label{eq:est 9}
H_2 &\leq \left[ \frac{4C^f_3}{s^{1-\gamma}_\bxi} (||\bu||_{\Cholder{}} + ||\bv||_{\Cholder{}}) +  \frac{4 C^f_2 ||\omega||_{\Cholder{}}}{\sqrt{s_\bxi}} \right] ||\bu - \bv||_{\Cholder{}}.
\end{align}
Applying the upper bounds on $H_1$ and $H_2$ shows that
\begin{align}\label{eq:est 10}
H &\leq \left[ \frac{4C^f_3}{s^{1-\gamma}_\bxi} (||\bu||_{\Cholder{}} + ||\bv||_{\Cholder{}}) +  \frac{4 C^f_2 (1+||\omega||_{\Cholder{}})}{\sqrt{s_\bxi}} \right] ||\bu - \bv||_{\Cholder{}}.
\end{align}
We substitute the upper bound on $H$ in \autoref{eq:bondbased ineq 1} to find that
\begin{align}\label{eq:lipshitz bound bondbased 2}
&\frac{|(\mathcal{L}^{\epsilon}_T(\bu)(\bx) - \mathcal{L}^{\epsilon}_T(\bv)(\bx)) - (\mathcal{L}^{\epsilon}_T(\bu)(\by) - \mathcal{L}^{\epsilon}_T(\bv)(\by))|}{|\bx - \by|^\gamma} \notag \\
&\leq \frac{2}{\epsilon \omega_d} \int_{H_1(\bzero)} \frac{J(|\bxi|)}{\sqrt{s_\bxi}} \left[ \frac{4C^f_3}{s^{1-\gamma}_\bxi} (||\bu||_{\Cholder{}} + ||\bv||_{\Cholder{}}) \right. \notag \\
&\quad \left. +  \frac{4 C^f_2 (1+||\omega||_{\Cholder{}})}{\sqrt{s_\bxi}} \right] ||\bu - \bv||_{\Cholder{}} d\bxi \notag \\
&= \left[ \frac{8C^f_3 \bar{J}_{3/2-\gamma}}{\epsilon^{5/2 - \gamma}} (||\bu||_{\Cholder{}} + ||\bv||_{\Cholder{}}) + \frac{8C^f_2(1+||\omega||_{\Cholder{}}) \bar{J}_{1}}{\epsilon^2} \right] ||\bu - \bv||_{\Cholder{}},
\end{align}
where $\bar{J}_\alpha$ is defined in \autoref{eq:def bar J}. Application of \autoref{eq:lipshitz bound bondbased 1} and \autoref{eq:lipshitz bound bondbased 2} deliver
\begin{align}\label{eq:lipshitz bound bondbased 3}
&\Choldernorm{\mathcal{L}^{\epsilon}_T(\bu) - \mathcal{L}^{\epsilon}_T(\bv)}{} \notag \\
&\leq \left[ \frac{8C^f_3 \bar{J}_{3/2-\gamma}}{\epsilon^{5/2 - \gamma}} (||\bu||_{\Cholder{}} + ||\bv||_{\Cholder{}}) + \frac{8C^f_2(2+||\omega||_{\Cholder{}}) \bar{J}_{1}}{\epsilon^2} \right] ||\bu - \bv||_{\Cholder{}},
\end{align}
and we have established the Lipschitz continuity of the non-local force due to tensile strain.

Now we establish the Lipschitz continuity for the non-local dilatational force. 
For any $\bu,\bv \in \Cholderz{\perdrd}$ we write
\begin{align}\label{eq:lipshitz norm state force}
&\Choldernorm{\mathcal{L}^{\epsilon}_D(\bu) - \mathcal{L}^{\epsilon}_D(\bv)}{\perdrd} \notag \\
&= \sup_{\bx\in D} |\mathcal{L}^{\epsilon}_D(\bu)(\bx) - \mathcal{L}^{\epsilon}_D(\bv)(\bx)| \notag \\
&+ \sup_{\bx, \by \in D, \bx \neq \by} \frac{|(\mathcal{L}^{\epsilon}_D(\bu)(\bx) - \mathcal{L}^{\epsilon}_D(\bv)(\bx)) - (\mathcal{L}^{\epsilon}_D(\bu)(\by) - \mathcal{L}^{\epsilon}_D(\bv)(\by))|}{|\bx - \by|^\gamma}.
\end{align}

The potential function $g$ can either be a quadratic function, e.g.,  $g(r) =\beta r^2/2$ or it can be a convex-concave function, see \autoref{ConvexConcaveFunctionGa}. Here we present the derivation of  Lipschitz continuity for the convex-concave type $g$. The proof for the quadratic potential functions $g$ is identical. 

Let $g$ be a bounded convex-concave potential function with bounded derivatives expressed by \autoref{eq:def Cgi}.
As in previous sections we use the notation \autoref{eq:notations consist} and \autoref{eq:def bar J} and begin by estimating $|\theta(\bx;\bu) - \theta(\bx;\bv)|$ where $\theta(\bx;\bu)$ is given by \autoref{eq:hydro strain}. Application of the  inequality $|\bubar_\bxi(\bx) - \bvbar_\bxi(\bx)| \leq 2 ||\bu - \bv||_{\Cholder{}}$, and a straightforward calculation shows that
\begin{align}\label{eq:est 11}
|\theta(\bx;\bu) - \theta(\bx;\bv)| &\leq 2\bar{J}_0 ||\bu - \bv||_{\Cholder{}}.
\end{align}
We now bound $|\theta(\bx;\bu) - \theta(\by; \bu)|$ as follows
\begin{align}
|\theta(\bx;\bu) - \theta(\by; \bu)| &= \bigg\vert\frac{1}{\omega_d}\int_{H_1(\bzero)}J(|\bxi|) \left[\omega(\bx+\epsilon\bxi) \bubar_\bxi(\bx) - \omega(\by + \epsilon\bxi) \bubar_\bxi(\by) \right] \cdot \be_\bxi d\bxi \bigg\vert \notag \\
&\leq \frac{1}{\omega_d}\int_{H_1(\bzero)}J(|\bxi|) \vert \omega(\bx+\epsilon\bxi) \bubar_\bxi(\bx) - \omega(\by + \epsilon\bxi) \bubar_\bxi(\by) \vert d\bxi \notag \\
&\leq \frac{1}{\omega_d}\int_{H_1(\bzero)}J(|\bxi|) |\omega(\bx+\epsilon\bxi)| \vert \bubar_\bxi(\bx) -  \bubar_\bxi(\by) \vert d\bxi \notag \\
& \; + \frac{1}{\omega_d}\int_{H_1(\bzero)}J(|\bxi|) |\omega(\bx + \epsilon \bxi) - \omega(\by+ \epsilon\bxi)| \vert \bubar_\bxi(\by) \vert d\bxi,
\end{align}
where we used $|\be_\bxi|=1$ and Cauchy's inequality in the first equation, added and subtracted $\omega(\bx+\epsilon \bxi)\bubar_\bxi(\by)$ in the second equation and used the triangle inequality. Applying $|\bubar_\bxi(\bx) - \bubar_\bxi(\by)| \leq 2|\bx - \by|^\gamma ||\bu||_{\Cholder{}}$, $|\omega(\bx+\epsilon \bxi) - \omega(\by + \epsilon\bxi)| \leq |\bx - \by|^\gamma ||\omega||_{\Cholder{}}$, and $|\bubar_\bxi(\by)|\leq 2||\bu||_{\Cholder{}}$ gives
\begin{align*}
|\theta(\bx;\bu) - \theta(\by; \bu)| &\leq \frac{1}{\omega_d} \int_{H_1(\bzero)} J(|\bxi|) 2 |\bx - \by|^\gamma ||\bu||_{\Cholder{}} d\bxi \notag \\
&\; +  \frac{1}{\omega_d} \int_{H_1(\bzero)} J(|\bxi|) |\bx - \by|^\gamma ||\omega||_{\Cholder{}} 2 ||\bu||_{\Cholder{}}, 
\end{align*}
i.e.,
\begin{align}\label{eq:est 12}
|\theta(\bx;\bu) - \theta(\by; \bu)|
&\leq 2\bar{J}_0(1+||\omega||_{\Cholder{}}) ||\bu||_{\Cholder{}} |\bx - \by|^\gamma.
\end{align}
We note that estimate \autoref{eq:est 11} and \autoref{eq:est 12} holds for all $\bx,\by \in D$ as well as for $\bx$ and $\by$ in the layer of thickness $2\epsilon$ surrounding $D$.

Using \autoref{eq:per statebased simple force expression} we have
\begin{align}
|\mathcal{L}^\epsilon_D(\bu)(\bx) - \mathcal{L}^\epsilon_D(\bv)(\bx)| &= \bigg\vert \frac{1}{\epsilon^2 \omega_d} \int_{H_1(\bzero)} \omega_\bxi(\bx) J(|\bxi|) [g'(\theta(\bx+\epsilon\bxi;\bu)) + g'(\theta(\bx;\bu))  \notag \\
&\qquad \qquad -g'(\theta(\bx+\epsilon\bxi;\bv)) - g'(\theta(\bx;\bv)) ] \be_\bxi d\bxi \bigg\vert \notag \\
&\leq \frac{1}{\epsilon^2 \omega_d} \int_{H_1(\bzero)} J(|\bxi|) \bigg\vert g'(\theta(\bx+\epsilon\bxi;\bu)) + g'(\theta(\bx;\bu)) \notag \\
&\qquad \qquad -g'(\theta(\bx+\epsilon\bxi;\bv)) - g'(\theta(\bx;\bv)) \bigg\vert d\bxi \notag \\
&\leq \frac{1}{\epsilon^2 \omega_d} \int_{H_1(\bzero)} J(|\bxi|) \left\{ \bigg\vert g'(\theta(\bx+\epsilon\bxi;\bu)) -g'(\theta(\bx+\epsilon\bxi;\bv)) \bigg\vert \right. \notag \\
&\qquad \qquad \left.+ \bigg\vert g'(\theta(\bx;\bu))  - g'(\theta(\bx;\bv)) \bigg\vert \right\} d\bxi.
\end{align}
Since $|g'(r_1) - g'(r_2)|\leq C^g_2 |r_1 - r_2|$, we have
\begin{align*}
|g'(\theta(\bx;\bu))  - g'(\theta(\bx;\bv))| &\leq C^g_2 |\theta(\bx;\bu)  - \theta(\bx;\bv)| \notag \\
&\leq 2 C^g_2 \bar{J}_0 ||\bu - \bv||_{\Cholder{}},
\end{align*}
where we used \autoref{eq:est 11}. Similarly we have
\begin{align*}
|g'(\theta(\bx+\epsilon\bxi;\bu))  - g'(\theta(\bx+\epsilon \bxi;\bv))| &\leq 2 C^g_2 \bar{J}_0 ||\bu - \bv||_{\Cholder{}}.
\end{align*}
and we arrive at the estimate
\begin{align}\label{eq:lipshitz bound statebased 1}
|\mathcal{L}^\epsilon_D(\bu)(\bx) - \mathcal{L}^\epsilon_D(\bv)(\bx)| &\leq \frac{4C^g_2\bar{J}_0^2}{\epsilon^2} ||\bu - \bv||_{\Cholder{}}.
\end{align}

Now we estimate
\begin{align*}
& \frac{|(\mathcal{L}^{\epsilon}_D(\bu)(\bx) - \mathcal{L}^{\epsilon}_D(\bv)(\bx)) - (\mathcal{L}^{\epsilon}_D(\bu)(\by) - \mathcal{L}^{\epsilon}_D(\bv)(\by))|}{|\bx - \by|^\gamma}.
\end{align*}
We write
\begin{align*}
\mathcal{L}^\epsilon_D(\bu)(\bx) - \mathcal{L}^\epsilon_D(\bv)(\bx) &=  \frac{1}{\epsilon^2 \omega_d} \int_{H_1(\bzero)} \omega_\bxi(\bx) J(|\bxi|) [g'(\theta(\bx+\epsilon\bxi;\bu)) + g'(\theta(\bx;\bu))  \notag \\
&\qquad \qquad -g'(\theta(\bx+\epsilon\bxi;\bv)) - g'(\theta(\bx;\bv)) ] \be_\bxi d\bxi
\end{align*}
and
\begin{align*}
\mathcal{L}^\epsilon_D(\bu)(\by) - \mathcal{L}^\epsilon_D(\bv)(\by) &=  \frac{1}{\epsilon^2 \omega_d} \int_{H_1(\bzero)} \omega_\bxi(\by) J(|\bxi|) [g'(\theta(\by+\epsilon\bxi;\bu)) + g'(\theta(\by;\bu))  \notag \\
&\qquad \qquad -g'(\theta(\by+\epsilon\bxi;\bv)) - g'(\theta(\by;\bv)) ] \be_\bxi d\bxi.
\end{align*}
to find
\begin{align}
&|(\mathcal{L}^\epsilon_D(\bu)(\bx) - \mathcal{L}^\epsilon_D(\bv)(\bx)) - (\mathcal{L}^\epsilon_D(\bu)(\by) - \mathcal{L}^\epsilon_D(\bv)(\by))| \notag \\
&= \bigg\vert \frac{1}{\epsilon^2 \omega_d} \int_{H_1(\bzero)} J(|\bxi|)\notag \\
&\; \bigg( \omega_\bxi(\bx)[g'(\theta(\bx+\epsilon\bxi;\bu)) + g'(\theta(\bx;\bu)) -g'(\theta(\bx+\epsilon\bxi;\bv)) - g'(\theta(\bx;\bv)) ]  \notag \\
&\; - \omega_\bxi(\by)[g'(\theta(\by+\epsilon\bxi;\bu)) + g'(\theta(\by;\bu)) -g'(\theta(\by+\epsilon\bxi;\bv)) - g'(\theta(\by;\bv)) ]\bigg) \be_\bxi d\bxi \bigg\vert \notag \\
&\leq \frac{1}{\epsilon^2 \omega_d} \int_{H_1(\bzero)} J(|\bxi|) \notag \\
&\; \bigg\vert \omega_\bxi(\bx)[g'(\theta(\bx+\epsilon\bxi;\bu)) + g'(\theta(\bx;\bu)) -g'(\theta(\bx+\epsilon\bxi;\bv)) - g'(\theta(\bx;\bv)) ]  \notag \\
&\; - \omega_\bxi(\by)[g'(\theta(\by+\epsilon\bxi;\bu)) + g'(\theta(\by;\bu)) -g'(\theta(\by+\epsilon\bxi;\bv)) - g'(\theta(\by;\bv)) ] \bigg\vert d\bxi \notag \\
&= \frac{1}{\epsilon^2 \omega_d} \int_{H_1(\bzero)} J(|\bxi|) \bigg\vert \bigg( \omega_\bxi(\bx)[g'(\theta(\bx+\epsilon\bxi;\bu)) -g'(\theta(\bx+\epsilon\bxi;\bv))] \notag \\
&\qquad \qquad - \omega_\bxi(\by)[g'(\theta(\by+\epsilon\bxi;\bu))  -g'(\theta(\by+\epsilon\bxi;\bv))] \bigg) \notag \\
&\qquad \qquad + \bigg( \omega_\bxi(\bx)[g'(\theta(\bx;\bu)) -g'(\theta(\bx;\bv))] \notag \\
&\qquad \qquad - \omega_\bxi(\by)[g'(\theta(\by;\bu))  -g'(\theta(\by;\bv))] \bigg) \bigg\vert d\bxi,
\end{align}
where we have rearranged the terms in last step. Application of the triangle inequality gives
\begin{align}
&|(\mathcal{L}^\epsilon_D(\bu)(\bx) - \mathcal{L}^\epsilon_D(\bv)(\bx)) - (\mathcal{L}^\epsilon_D(\bu)(\by) - \mathcal{L}^\epsilon_D(\bv)(\by))| \notag \\
&\leq \frac{1}{\epsilon^2 \omega_d} \int_{H_1(\bzero)} J(|\bxi|) \bigg( \bigg\vert \omega_\bxi(\bx)[g'(\theta(\bx+\epsilon\bxi;\bu)) -g'(\theta(\bx+\epsilon\bxi;\bv))] \notag \\
&\qquad \qquad - \omega_\bxi(\by)[g'(\theta(\by+\epsilon\bxi;\bu))  -g'(\theta(\by+\epsilon\bxi;\bv))] \bigg\vert \notag \\
&\qquad \qquad + \bigg\vert \omega_\bxi(\bx)[g'(\theta(\bx;\bu)) -g'(\theta(\bx;\bv))] \notag \\
&\qquad \qquad - \omega_\bxi(\by)[g'(\theta(\by;\bu))  -g'(\theta(\by;\bv))] \bigg\vert \bigg)  d\bxi.
\end{align}
Now write $h_\bxi:\bbR^d \times \bbR^d \to \bbR^+$ given by
\begin{align}\label{eq:def h}
h_\bxi(\bx,\by) &:= \bigg\vert \omega_\bxi(\bx)[g'(\theta(\bx;\bu)) -g'(\theta(\bx;\bv))] - \omega_\bxi(\by)[g'(\theta(\by;\bu))  -g'(\theta(\by;\bv))] \bigg\vert.
\end{align}
and 
\begin{align}\label{eq:est 13}
&|(\mathcal{L}^\epsilon_D(\bu)(\bx) - \mathcal{L}^\epsilon_D(\bv)(\bx)) - (\mathcal{L}^\epsilon_D(\bu)(\by) - \mathcal{L}^\epsilon_D(\bv)(\by))| \notag \\
&\leq \frac{1}{\epsilon^2 \omega_d} \int_{H_1(\bzero)} J(|\bxi|) (h_\bxi(\bx+ \epsilon\bxi, \by+ \epsilon\bxi) + h(\bx,\by) ) d\bxi.
\end{align}
We now estimate $h_\bxi(\bx,\by)$ for any $\bx,\by $ in $D$ and in the layer of thickness $\epsilon$ surrounding $D$. 

Proceeding as before we define $r:[0,1]\times D \to \bbR$ as follows
\begin{align}\label{eq:def r}
r(l,\bx) &:= \theta(\bx;\bv) + l(\theta(\bx;\bu) - \theta(\bx;\bv)),
\end{align}
so $\frac{\partial r(l,\bx)}{\partial l} = \theta(\bx;\bu) - \theta(\bx;\bv)$. We also have
\begin{align}\label{eq:est 14}
g'(\theta(\bx;\bu)) -g'(\theta(\bx;\bv)) &= g'(r(1,\bx)) - g'(r(0,\bx)) \notag \\
&= \int_0^1 \frac{\partial g'(r(l,\bx))}{\partial l} dl \notag \\
&= \int_0^1 g''(r(l,\bx)) (\theta(\bx;\bu) - \theta(\bx;\bv)) dl.
\end{align}
Similarly,
\begin{align}\label{eq:est 15}
g'(\theta(\by;\bu)) -g'(\theta(\by;\bv)) &= \int_0^1 g''(r(l,\by)) (\theta(\by;\bu) - \theta(\by;\bv)) dl.
\end{align}
Substitution of \autoref{eq:est 14} and \autoref{eq:est 15} in $h_\bxi(\bx,\by)$  gives
\begin{align*}
h_\bxi(\bx,\by) &= \bigg\vert \int_0^1 ( \omega_\bxi(\bx)  g''(r(l,\bx)) (\theta(\bx;\bu) - \theta(\bx;\bv)) \notag \\
&\qquad \qquad - \omega_\bxi(\by)  g''(r(l,\by)) (\theta(\by;\bu) - \theta(\by;\bv)) ) dl \bigg\vert \notag \\
&\leq \int_0^1 \bigg\vert \omega_\bxi(\bx)  g''(r(l,\bx)) (\theta(\bx;\bu) - \theta(\bx;\bv)) \notag \\
&\qquad \qquad - \omega_\bxi(\by)  g''(r(l,\by)) (\theta(\by;\bu) - \theta(\by;\bv)) \bigg\vert dl.
\end{align*}
Adding and subtracting $\omega_\bxi(\bx)g''(r(l,\bx)) (\theta(\by;\bu) - \theta(\by;\bv))$ gives
\begin{align}\label{eq:est 16}
h_\bxi(\bx,\by) &\leq \int_0^1 |\omega_\bxi(\bx)|\, |g''(r(l,\bx))|\, |(\theta(\bx;\bu) - \theta(\bx;\bv)) - (\theta(\by;\bu) - \theta(\by;\bv))| dl \notag \\
&+ \int_0^1 |\omega_\bxi(\bx)g''(r(l,\bx)) - \omega_\bxi(\by)g''(r(l,\by))|\, |\theta(\by;\bu) - \theta(\by;\bv)| dl \notag \\
&=: I_1 + I_2,
\end{align}
For $I_1$, we note that $0\leq \omega(\bx) \leq 1$ and $|g''(r)| \leq C^g_2$ and proceed further to find that
\begin{align}
I_1 &\leq C^g_2 |(\theta(\bx;\bu) - \theta(\bx;\bv)) - (\theta(\by;\bu) - \theta(\by;\bv))| \notag \\
&= C^g_2 |\theta(\bx;\bu - \bv) - \theta(\by;\bu-\bv)|.
\end{align}
Using the estimate given in \autoref{eq:est 12} we see that
\begin{align}\label{eq:est 17}
I_1 &\leq 2 \bar{J}_0 C^g_2 (1+||\omega||_{\Cholder{}})||\bu - \bv||_{\Cholder{}} |\bx - \by|^\gamma.
\end{align}
Now we apply the inequality given in \autoref{eq:est 11} to $I_2$ to find that
\begin{align*}
I_2 &\leq 2 \bar{J}_0 ||\bu - \bv||_{\Cholder{}} \int_0^1 |\omega_\bxi(\bx)g''(r(l,\bx)) - \omega_\bxi(\by)g''(r(l,\by))| dl.
\end{align*}
Adding and subtracting $\omega_\bxi(\bx)g''(r(l,\by))$ gives
\begin{align*}
I_2 &\leq 2 \bar{J}_0 ||\bu - \bv||_{\Cholder{}} \int_0^1 |\omega_\bxi(\bx)|\, |g''(r(l,\bx)) - g''(r(l,\by))| dl \notag \\
&\quad + 2 \bar{J}_0 ||\bu - \bv||_{\Cholder{}} \int_0^1 |\omega_\bxi(\bx) - \omega_\bxi(\by)|\, |g''(r(l,\by))| dl \notag \\
&\leq 2 C^g_3 \bar{J}_0 ||\bu - \bv||_{\Cholder{}} \int_0^1 |r(l,\bx) - r(l,\by)| dl \notag \\
&\quad + 2 C^g_2 \bar{J}_0 ||\bu - \bv||_{\Cholder{}} \int_0^1 |\omega_\bxi(\bx) - \omega_\bxi(\by)| dl.
\end{align*}
The quantity  $|r(l,\bx) - r(l,\by)| $ (see \autoref{eq:def r}) can be estimated as follows
\begin{align}
&|r(l,\bx) - r(l,\by)| \notag \\
&= | (1-l) \theta(\bx;\bv) + l\theta(\bx;\bu) - ((1-l) \theta(\by;\bv) + l\theta(\by;\bu))| \notag \\
&\leq |1-l|\, |\theta(\bx;\bv) - \theta(\by;\bv)| + |l|\, |\theta(\bx;\bu) - \theta(\by;\bu)| \notag \\
&\leq |\theta(\bx;\bv) - \theta(\by;\bv)| + |\theta(\bx;\bu) - \theta(\by;\bu)| \notag \\
&\leq 2\bar{J}_0 (1+||\omega||_{\Cholder{}}) ||\bv||_{\Cholder{}} |\bx - \by|^\gamma + 2\bar{J}_0 (1+||\omega||_{\Cholder{}}) ||\bu||_{\Cholder{}} |\bx - \by|^\gamma \notag \\
&= 2\bar{J}_0 (1+||\omega||_{\Cholder{}})(||\bu||_{\Cholder{}} + ||\bv||_{\Cholder{}}) |\bx - \by|^\gamma,
\end{align}
where we used the fact that $l\in [0,1]$ and \autoref{eq:est 12}. Using the inequality above and $|\omega_\bxi(\bx) - \omega_\bxi(\by)| \leq 2|\bx - \by|^\gamma ||\omega||_{\Cholder{}}$ we get
\begin{align}\label{eq:est 18}
I_2 &\leq 2 C^g_3 \bar{J}_0 ||\bu - \bv||_{\Cholder{}} 2\bar{J}_0 (1+||\omega||_{\Cholder{}})(||\bu||_{\Cholder{}} + ||\bv||_{\Cholder{}}) \, |\bx - \by|^\gamma \notag \\
&\quad + 2 C^g_2 \bar{J}_0 ||\bu - \bv||_{\Cholder{}} 2|\bx - \by|^\gamma |\,|\omega||_{\Cholder{}} \notag \\
&\leq 4\bar{J}_0 (1+||\omega||_{\Cholder{}})\, [ C^g_3 \bar{J}_0(||\bu||_{\Cholder{}} + ||\bv||_{\Cholder{}}) + C^g_2] \, ||\bu - \bv||_{\Cholder{}} \,|\bx - \by|^\gamma.
\end{align}
Substituting \autoref{eq:est 17} and \autoref{eq:est 18} into \autoref{eq:est 16} gives
\begin{align}\label{eq:est 19}
&h_\bxi(\bx, \by) \notag \\
&\leq 6\bar{J}_0 (1+||\omega||_{\Cholder{}})\, [ C^g_3 \bar{J}_0(||\bu||_{\Cholder{}} + ||\bv||_{\Cholder{}}) + C^g_2] \, ||\bu - \bv||_{\Cholder{}} \,|\bx - \by|^\gamma.
\end{align}

We now apply \autoref{eq:est 19} to \autoref{eq:est 13} and divide both sides by $|\bx - \by|^\gamma$ to see that
\begin{align}\label{eq:lipshitz bound statebased 2}
&\frac{|(\mathcal{L}^\epsilon_D(\bu)(\bx) - \mathcal{L}^\epsilon_D(\bv)(\bx)) - (\mathcal{L}^\epsilon_D(\bu)(\by) - \mathcal{L}^\epsilon_D(\bv)(\by))|}{|\bx - \by|^\gamma} \notag \\
&\leq \frac{1}{\epsilon^2 \omega_d} \int_{H_1(\bzero)} J(|\bxi|) \notag \\
&\quad 2\times 6\bar{J}_0 (1+||\omega||_{\Cholder{}})\, [ C^g_3 \bar{J}_0(||\bu||_{\Cholder{}} + ||\bv||_{\Cholder{}}) + C^g_2] \, ||\bu - \bv||_{\Cholder{}} \, d\bxi \notag \\
&= \frac{12\bar{J}^2_0 (1+||\omega||_{\Cholder{}})\, [ C^g_3 \bar{J}_0(||\bu||_{\Cholder{}} + ||\bv||_{\Cholder{}}) + C^g_2]}{\epsilon^2} ||\bu - \bv||_{\Cholder{}}.
\end{align}

Collecting results inequalities \autoref{eq:lipshitz bound statebased 1} and \autoref{eq:lipshitz bound statebased 2} deliver the upper bound given by
\begin{align}\label{eq:lipshitz bound statebased 3}
&||\mathcal{L}^\epsilon_D(\bu) - \mathcal{L}^\epsilon_D(\bv)||_{\Cholder{}} \notag \\
&\leq \frac{16\bar{J}^2_0 (1+||\omega||_{\Cholder{}})\, [ C^g_3 \bar{J}_0(||\bu||_{\Cholder{}} + ||\bv||_{\Cholder{}}) + C^g_2]}{\epsilon^2} ||\bu - \bv||_{\Cholder{}}.
\end{align}

\paragraph{Lipschitz continuity for $\mathcal{L}^\epsilon(\bu)$}
Using \autoref{eq:lipshitz bound bondbased 3} and \autoref{eq:lipshitz bound statebased 3} we get
\begin{align}
&||\mathcal{L}^\epsilon(\bu) - \mathcal{L}^\epsilon(\bv)||_{\Cholder{}} \notag \\
&\leq \bigg( \frac{8C^f_3 \bar{J}_{3/2-\gamma}}{\epsilon^{5/2 - \gamma}} (||\bu||_{\Cholder{}} + ||\bv||_{\Cholder{}}) + \frac{8C^f_2(2+||\omega||_{\Cholder{}}) \bar{J}_{1}}{\epsilon^2} \notag \\
&\: + \frac{16\bar{J}^2_0 (1+||\omega||_{\Cholder{}})\, [ C^g_3 \bar{J}_0(||\bu||_{\Cholder{}} + ||\bv||_{\Cholder{}}) + C^g_2]}{\epsilon^2} \bigg)||\bu - \bv||_{\Cholder{}}.
\end{align}
Let $\alpha(\gamma)$ defined as follows: $\alpha(\gamma) = 0$ if $\gamma \geq 1/2$ and $\alpha(\gamma) = 1/2 - \gamma $ if $\gamma \leq 1/2$. It is easy to verify that, for all $\gamma \in (0,1]$ and $0< \epsilon \leq 1$
\begin{align}
\max \bigg\{ \frac{1}{\epsilon^2}, \frac{1}{\epsilon^{5/2 - \gamma}} \bigg\} \leq \frac{1}{\epsilon^{2+\alpha(\gamma)}}.
\end{align}
Using $\alpha(\gamma)$ and renaming the constants we have
\begin{align}
&||\mathcal{L}^\epsilon(\bu) - \mathcal{L}^\epsilon(\bv)||_{\Cholder{}} \notag \\
&\leq \frac{L_1 (1+||\omega||_{\Cholder{}}) (1+ ||\bu||_{\Cholder{}} + ||\bv||_{\Cholder{}})}{\epsilon^{2+\alpha(\gamma)}} ||\bu - \bv||_{\Cholder{}}.
\end{align}
To complete the proof of \autoref{eq:lipschitz property of F}, we substitute the inequality above into \autoref{eq:norm of F in X} to obtain
\begin{align}
&\normX{F^\epsilon(y,t) - F^\epsilon(z,t)}{X} \notag \\
&\leq ||y^2 - z^2||_{\Cholder{}} + \frac{L_1 (1+||\omega||_{\Cholder{}}) (1+ ||y^1||_{\Cholder{}} + ||z^1||_{\Cholder{}})}{\epsilon^{2+\alpha(\gamma)}} ||y^1 - z^1||_{\Cholder{}} \notag \\
&\leq \frac{L_1 (1+||\omega||_{\Cholder{}}) (1+ ||y||_X + ||z||_X)}{\epsilon^{2+\alpha(\gamma)}} ||y - z||_X,
\end{align}
and \autoref{eq:lipschitz property of F} is proved.

\subsubsection{Bound on the non-local force in the H\"older norm}
In this subsection, we bound  $||\mathcal{L^\epsilon}(\bu)||_{\Cholder{}}$ from above. It follows from \autoref{eq:per bondbased simple force expression} and a straightforward calculation similar to the previous sections that
\begin{align}
|\mathcal{L}^\epsilon_T(\bu)(\bx)| &\leq \frac{2C^f_1 \bar{J}_{1/2}}{\epsilon^{3/2}}, \notag \\
\frac{|\mathcal{L}^\epsilon_T(\bu)(\bx) - \mathcal{L}^\epsilon_T(\bu)(\by)|}{|\bx-\by|^\gamma} &\leq \frac{4C^f_2 \bar{J}_1 ||\bu||_{\Cholder{}} + 4 C^f_1 \bar{J}_{1/2} ||\omega||_{\Cholder{}}}{\epsilon^2}.
\end{align}

Next we consider the non-local dilatational force $\mathcal{L}^\epsilon_D$. We show how to calculate the bounds for the case of a convex-concave potential function $g$. When $g$ is quadratic we can still proceed along identical lines. We use the formula for  $\mathcal{L}^\epsilon_D(\bu)(\bx)$ given by \autoref{eq:per statebased simple force expression} and perform a straightforward calculation to obtain the upper bound given by
\begin{align}
|\mathcal{L}^\epsilon_D(\bu)(\bx)| &\leq \frac{2C^g_1 \bar{J}_{0}}{\epsilon^{2}}.
\end{align}
We have the estimate
\begin{align}
&|\mathcal{L}^\epsilon_D(\bu)(\bx) - \mathcal{L}^\epsilon_D(\bu)(\by)|  \notag \\
&\leq \frac{1}{\epsilon^2\omega_d} \int_{H_1(\bzero)} J(|\bxi|) \bigg\vert \omega_\bxi(\bx) (g'(\theta(\bx+\epsilon \bxi;\bu)) + g'(\theta(\bx;\bu)))) \notag \\
&\qquad \qquad - \omega_\bxi(\by) (g'(\theta(\by+\epsilon \bxi;\bu)) + g'(\theta(\by;\bu)))) \bigg\vert d\bxi \notag \\
&\leq \frac{1}{\epsilon^2\omega_d} \int_{H_1(\bzero)} J(|\bxi|) \bigg\vert \omega_\bxi(\bx) g'(\theta(\bx+\epsilon \bxi;\bu)) - \omega_\bxi(\by) g'(\theta(\by+\epsilon \bxi;\bu))\bigg\vert d\bxi \notag \\
&+ \frac{1}{\epsilon^2\omega_d} \int_{H_1(\bzero)} J(|\bxi|) \bigg\vert \omega_\bxi(\bx) g'(\theta(\bx;\bu)) - \omega_\bxi(\by) g'(\theta(\by;\bu))\bigg\vert d\bxi.
\end{align}
Using $|\omega_\bxi(\bx) - \omega_\bxi(\by)| \leq 2|\bx - \by|^\gamma ||\omega||_{\Cholder{}}$, $|g'(r_1) - g'(r_2)| \leq C^g_2 |r_1 - r_2|$, $|g'(r)| \leq C^g_1$, and the estimate on $|\theta(\bx;\bu) - \theta(\by;\bu)|$ given by \autoref{eq:est 12}, we obtain
\begin{align}
&|\mathcal{L}^\epsilon_D(\bu)(\bx) - \mathcal{L}^\epsilon_D(\bu)(\by)|  \notag \\
&\leq\frac{ [2 \bar{J}_0 C^g_2 (1+||\omega||_{\Cholder{}})||\bu||_{\Cholder{}} + 2C^g_1||\omega||_{\Cholder{}} ||\bu||_{\Cholder{}} ]}{\epsilon^2}\, |\bx - \by|^\gamma.
\end{align}
Last we combine results and rename the constants to get
\begin{align}
\label{rhs est}
||\mathcal{L}^\epsilon(\bu)||_{\Cholder{}} &\leq \frac{L_2(1+||\omega||_{\Cholder{}}) (1+||\bu||_{\Cholder{}})}{\epsilon^2}.
\end{align}
This completes the proof of \autoref{eq:bound on F}.

\subsection{Proof of Proposition 2}\label{ss:proof prop 4}
Given $\bu,\bv \in L^2_0(D;\bbR^d)$ we find upper bounds on the Lipschitz continuity of the nonlocal force with respect to the $L^2$ norm. Motivated by the inequality
\begin{align}
||\mathcal{L}^\epsilon(\bu) - \mathcal{L}^\epsilon(\bv)||_{L^2} &\leq ||\mathcal{L}^\epsilon_T(\bu) - \mathcal{L}^\epsilon_T(\bv)||_{L^2} + ||\mathcal{L}^\epsilon_D(\bu) - \mathcal{L}^\epsilon_D(\bu)||_{L^2}.,
\end{align}
we bound the Lipschitz continuity of the nonlocal forces due to tensile strain and dilatational strain separately. 
We study $\mathcal{L}^\epsilon_T$ first. It is evident from \autoref{eq:per bondbased simple force expression} and using the estimate $|f'(r_1) - f'(r_2)|\leq C^f_2 |r_1 - r_2|$, and arguments similar to previous sections  that we have
\begin{align}\label{eq:est 19.1}
&|\mathcal{L}^\epsilon_T(\bu)(\bx) - \mathcal{L}^\epsilon_T(\bv)(\bx)| \notag \\
&\leq \frac{2}{\epsilon\omega_d} \int_{H_1(\bzero)} \frac{J(|\bxi|)}{\sqrt{s_\bxi}} |f'(\bubar_\bxi(\bx).\be_\bxi/\sqrt{s_\bxi}) - f'(\bvbar_\bxi(\bx).\be_\bxi/\sqrt{s_\bxi})| d\bxi \notag \\
&\leq \frac{2C^f_2}{\epsilon^2 \omega_d} \int_{H_1(\bzero)} \frac{J(|\bxi|)}{|\bxi|} |\bubar_\bxi(\bx) - \bvbar_\bxi(\bx)| d\bxi,
\end{align}
where we also substituted $s_\bxi = \epsilon |\bxi|$. 

We apply \autoref{eq:ineq symm square} to  \autoref{eq:est 19.1} with  $C=\frac{2C^f_2}{\epsilon^2}$, $\alpha = 1$, and $p(\bxi) = |\bubar_\bxi(\bx) - \bvbar_\bxi(\bx)|$ to get
\begin{align}
& ||\mathcal{L}^\epsilon_T(\bu) - \mathcal{L}^\epsilon_T(\bv)||_{L^2}^2 \notag \\
&\leq \int_D |\mathcal{L}^\epsilon_T(\bu)(\bx) - \mathcal{L}^\epsilon_T(\bv)(\bx)|^2 d\bx \notag \\
&\leq  \int_D \left(\frac{2C^f_2}{\epsilon^2} \right)^2 \frac{\bar{J}_1}{\omega_d} \int_{H_1(\bzero)} \frac{J(|\bxi|)}{|\bxi|} |\bubar_\bxi(\bx) - \bvbar_\bxi(\bx)|^2 d\bxi d\bx \notag \\
&= \left(\frac{2C^f_2}{\epsilon^2} \right)^2 \frac{\bar{J}_1}{\omega_d} \int_{H_1(\bzero)} \frac{J(|\bxi|)}{|\bxi|} \left[  \int_D |\bubar_\bxi(\bx) - \bvbar_\bxi(\bx)|^2 d\bx \right] d\bxi,
\end{align}
where we interchanged integration in last step. Using
\begin{align}
\int_D |\bubar_\bxi(\bx) - \bvbar_\bxi(\bx)|^2 d\bx &\leq 2 ||\bu - \bv||^2_{L^2}
\end{align}
we conclude that
\begin{align}\label{eq:est 20}
||\mathcal{L}^\epsilon_T(\bu) - \mathcal{L}^\epsilon_T(\bv)||_{L^2} \leq \frac{4 C^f_2 \bar{J}_1}{\epsilon^2} ||\bu - \bv||_{L^2}.
\end{align}

In estimating $||\mathcal{L}^\epsilon_D(\bu) -\mathcal{L}^\epsilon_D(\bv)||_{L^2}$ we will consider convex-concave $g$ noting that the case of quadratic $g$ is dealt in a similar fashion. From \autoref{eq:per statebased simple force expression} and using estimate $|g'(r_1) - g'(r_2)|\leq C^g_2 |r_1 - r_2|$, and proceeding as before we have
\begin{align}\label{lipd inequal}
&|\mathcal{L}^\epsilon_D(\bu)(\bx) - \mathcal{L}^\epsilon_D(\bv)(\bx)| \notag \\
&\leq \frac{1}{\epsilon^2 \omega_d} \int_{H_1(\bzero)} J(|\bxi|) [|g'(\theta(\bx+\epsilon\bxi;\bu)) - g'(\theta(\bx+\epsilon\bxi;\bv))| \notag \\
&\qquad \quad + |g'(\theta(\bx;\bu)) - g'(\theta(\bx;\bv))|] d\bxi \notag \\
&\leq \frac{C^g_2}{\epsilon^2 \omega_d} \int_{H_1(\bzero)} J(|\bxi|) [|\theta(\bx+\epsilon\bxi;\bu) - \theta(\bx+\epsilon\bxi;\bv)| + |\theta(\bx;\bu) - \theta(\bx;\bv)|] d\bxi \notag \\
&= \frac{C^g_2}{\epsilon^2 \omega_d} \int_{H_1(\bzero)} J(|\bxi|) [|\theta(\bx+\epsilon\bxi;\bu-\bv) | + |\theta(\bx;\bu-\bv)|] d\bxi.
\end{align}
Squaring \autoref{lipd inequal} and applying inequality \autoref{eq:ineq symm square} with $C=\frac{C^g_2}{\epsilon^2}$, $\alpha = 0$, and $p(\bxi) = |\theta(\bx+\epsilon\bxi;\bu-\bv) | + |\theta(\bx;\bu-\bv)|$ gives
\begin{align}\label{l2bound}
&||\mathcal{L}^\epsilon_D(\bu) - \mathcal{L}^\epsilon_D(\bv)||^2_{L^2} \notag \\
&\leq \int_D \left(\frac{C^g_2}{\epsilon^2} \right)^2  \frac{\bar{J}_0}{\omega_d}\int_{H_1(\bzero)}  J(|\bxi|) (|\theta(\bx+\epsilon\bxi;\bu-\bv) | + |\theta(\bx;\bu-\bv)|)^2 d\bxi d\bx \notag \\
&\leq \left(\frac{C^g_2}{\epsilon^2} \right)^2  \frac{\bar{J}_0}{\omega_d}\int_{H_1(\bzero)}  J(|\bxi|) \left[ \int_D 2(|\theta(\bx+\epsilon\bxi;\bu-\bv) |^2 + |\theta(\bx;\bu-\bv)|^2) d\bx \right] d\bxi,
\end{align}
where we used Cauchy's inequality and exchanged integration in the last step. It is easy to verify that
\begin{align*}
\int_D |\theta(\bx+\epsilon \bxi; \bu)|^2 d\bx &\leq 2 \bar{J}_0^2 ||\bu||_{L^2}^2 
\end{align*}
holds for all $\bxi \in H_1(\bzero)$. Combining this estimate and \autoref{l2bound} we see that
\begin{align}\label{eq:est 21}
||\mathcal{L}^\epsilon_D(\bu) - \mathcal{L}^\epsilon_D(\bv)||_{L^2} &\leq \frac{4 C^g_2 \bar{J}_0^2}{\epsilon^2 }||\bu - \bv||_{L^2}.
\end{align}
Estimates \autoref{eq:est 20} and \autoref{eq:est 21} together delivers (after renaming the constants)
\begin{align}
||\mathcal{L}^\epsilon(\bu) - \mathcal{L}^\epsilon(\bv)||_{L^2} &\leq \frac{L_3}{\epsilon^2} ||\bu - \bv||_{L^2},
\end{align}
where $L_3$ is given by \autoref{eq:def L3}.
This completes the proof of \sautoref{Proposition}{prop:lipschitz L2}.

\section{Energy stability of the semi-discrete scheme}\label{ss:stab proof}
In this section, we establish \autoref{thm:stab semi} for convex-concave potential functions $g$ as well as for quadratic potential functions. 
We recall the semi-discrete problem introduced in \autoref{semidiscrete}.
We first introduce the semi-discrete boundary condition by setting $\hat{\bu}_i(t) = \bzero$ for all $t$ and for all $\bx_i \notin D$.
Let $\{\hat{\bu}_i(t)\}_{i,\bx_i\in D}$ denote the semi-discrete approximate solution which satisfies the following evolution, for all $t\in [0,T]$ and $i$ such that $\bx_i\in D$,
\begin{align}\label{eq:stab 1}
\ddot{\hat{\bu}}_i(t) = \mathcal{L}^\epsilon(\hat{\bu}(t))(\bx_i) + \bb(\bx_i,t),
\end{align}
where $\hat{\bu}(t)$ is the piecewise constant extension of $\{\buhat(t)\}_{i,\bx_i \in D}$, given by
\begin{align*}
\hat{\bu}(t,\bx) = \sum_{i, \bx_i \in D} \hat{\bu}_i(t) \chi_{U_i}(\bx).
\end{align*}
Let $\hat{\mathcal{L}^\epsilon}(\buhat(t))(\bx)$ be defined as
\begin{align*}
\hat{\mathcal{L}^\epsilon}(\buhat(t))(\bx) = \sum_{i, \bx_i \in D} \mathcal{L}^\epsilon(\buhat(t))(\bx_i) \chi_{U_i}(\bx)
\end{align*}
and define $\hat{\bb}(t)$ similarly. From \autoref{eq:stab 1} noting the definition of piecewise constant extension
\begin{align}\label{eq:stab 2}
\ddot{\hat{\bu}}(\bx,t) &= \hat{\mathcal{L}^\epsilon}(\hat{\bu}(t))(\bx) + \hat{\bb}(\bx,t) \notag \\
&= \mathcal{L}^\epsilon(\hat{\bu}(t))(\bx) + \hat{\bb}(\bx,t) + \sigma(\bx,t),
\end{align}
where the error term $\sigma(\bx,t)$ is given by
\begin{align}
\sigma(\bx,t) &:= \hat{\mathcal{L}^\epsilon}(\hat{\bu}(t))(\bx) - \mathcal{L}^\epsilon(\hat{\bu}(t))(\bx).
\end{align}
We split $\sigma$ into two parts
\begin{align}
\sigma(\bx,t) &=  \left[ \hat{\mathcal{L}^\epsilon_T}(\hat{\bu}(t))(\bx) - \mathcal{L}^\epsilon_T(\hat{\bu}(t))(\bx) \right] + \left[\hat{\mathcal{L}^\epsilon_D}(\hat{\bu}(t))(\bx) - \mathcal{L}^\epsilon_D(\hat{\bu}(t))(\bx) \right] \notag \\
&=: \sigma_T(\bx,t) + \sigma_D(\bx,t).
\end{align}
Multiplying both sides of \autoref{eq:stab 2} by $\dot{\buhat}(t)$ and integrating over $D$ gives
\begin{align}\label{eq:stab 2.1}
(\ddot{\hat{\bu}}(t), \dot{\buhat}(t)) &= (\mathcal{L}^\epsilon(\buhat(t)), \dot{\buhat}(t)) + (\hat{\bb}(t), \dot{\buhat}(t)) + (\sigma(t), \dot{\buhat}(t)),
\end{align}
where $(\cdot, \cdot)$ denotes the $L^2$-inner product. 

\subsection{Estimating $\sigma$}
We proceed by estimating $L^2$-norm of $\sigma(t)$. It follows easily from \autoref{eq:per bondbased simple force expression} that
\begin{align}
|\sigma_T(\bx,t)| &\leq \frac{4C^f_1 \bar{J}_{1/2}}{\epsilon^{3/2}} \quad \Rightarrow ||\sigma_T(t)||_{L^2} \leq \frac{4 C^f_1 \bar{J}_{1/2}\sqrt{|D|}}{\epsilon^{3/2}}.
\end{align}
We now deal with two cases of $g$ separately. 

\textbf{1. Convex-concave type $g$: }In this case, we can easily show from \autoref{eq:per statebased simple force expression} that
\begin{align}
|\sigma_D(\bx,t)| &\leq \frac{4C^g_1 \bar{J}_0}{\epsilon^2} \quad \Rightarrow ||\sigma_D(t)||_{L^2} \leq \frac{4 C^g_1 \bar{J}_0\sqrt{|D|}}{\epsilon^2}.
\end{align}

\textbf{2. Quadratic type $g$: }In this case we have $g'(r) = g''(0) r$. Let $\bx \in U_i$, i.e. in the unit cell of the $i^{th}$ mesh node. To simplify the calculations let $\bu = \buhat(t)$ (and later we will use the fact that $\buhat$ is piecewise constant function). From \autoref{eq:per statebased simple force expression}, we have
\begin{align}
|\sigma_D(\bx,t)| &= |\mathcal{L}^\epsilon_D(\bu)(\bx_i) - \mathcal{L}^\epsilon_D(\bu)(\bx)| \notag \\
&= \bigg\vert \frac{g''(0)}{\epsilon^2 \omega_d} \int_{H_1(\bzero)} J(|\bxi|) \bigg[ \omega_{\bxi}(\bx_i) (\theta(\bx_i + \epsilon \bxi; \bu) + \theta(\bx_i;\bu)) \notag \\
&\qquad -  \omega_{\bxi}(\bx) (\theta(\bx + \epsilon \bxi; \bu) + \theta(\bx;\bu)) \bigg]\be_\bxi d\bxi \bigg\vert.
\end{align}

Now consider the function $\ba(\bx,\bxi)$ defined as 
\begin{align}\label{eq:def a}
\ba(\bx, \bxi) = \theta(\bx+\epsilon\bxi; \bu).
\end{align}
We then have
\begin{align}
&|\sigma_D(\bx,t)| \notag \\
&= \bigg\vert \frac{g''(0)}{\epsilon^2 \omega_d} \int_{H_1(\bzero)} J(|\bxi|) \bigg[ \omega_{\bxi}(\bx_i) (\ba(\bx_i,\bxi) + \ba(\bx_i, \bzero)) \notag \\
&\qquad - \omega_{\bxi}(\bx) (\ba(\bx,\bxi) + \ba(\bx, \bzero)) \bigg]\be_\bxi d\bxi \bigg\vert \notag \\
&\leq \frac{g''(0)}{\epsilon^2 \omega_d} \int_{H_1(\bzero)} J(|\bxi|) (|\ba(\bx_i,\bxi)| + |\ba(\bx_i,\bzero)| + |\ba(\bx,\bxi)| + |\ba(\bx,\bzero)|) d\bxi.
\end{align}
Let 
\begin{align}\label{eq:def b}
b_\bxi := |\ba(\bx_i,\bxi)| + |\ba(\bx_i,\bzero)| + |\ba(\bx,\bxi)| + |\ba(\bx,\bzero)|
\end{align}
then using the inequality \autoref{eq:ineq symm square} with $C=\frac{g''(0)}{\epsilon^2}$, $\alpha = 0$, and $p(\bxi) = b_\bxi$, we get
\begin{align}\label{eq:stab 3}
|\sigma_D(\bx,t)|^2 &\leq\left(\frac{g''(0)}{\epsilon^2}  \right)^2 \frac{\bar{J}_0}{\omega_d} \int_{H_1(\bzero)} J(|\bxi|) b_\bxi^2 d\bxi.
\end{align}
Thus on an interchange of integration we have
\begin{align}\label{eq:stab 4}
||\sigma_D(t)||^2_{L^2} &= \int_D |\sigma_D(\bx,t)|^2 d\bx \notag \\
&\leq \sum_{i,\bx_i \in D} \int_{U_i} |\sigma_D(\bx,t)|^2 d\bx \notag \\
&\leq \left(\frac{g''(0)}{\epsilon^2}  \right)^2 \frac{\bar{J}_0}{\omega_d} \int_{H_1(\bzero)} J(|\bxi|) \left[\sum_{i,\bx_i \in D} \int_{U_i} b^2_\bxi d\bx \right] d\bxi.
\end{align}
We denote the term inside square bracket as $I$ and estimate it next. Recalling the definition of $b_\bxi$ in \autoref{eq:def b} and using the identity $(\sum_{n=1}^4 c_n)^2 \leq 4 \sum_{n=1}^4 c^2_n$ we have
\begin{align}\label{eq:stab 5}
I &\leq 4 \sum_{i,\bx_i \in D} \int_{U_i} (|\ba(\bx_i,\bxi)|^2 + |\ba(\bx_i,\bzero)|^2 + |\ba(\bx,\bxi)|^2 + |\ba(\bx,\bzero)|^2) d\bx. 
\end{align}
For $\bx$ either in $D$ or in layer of thickness $\epsilon$ surrounding $D$ take $\bxi \in H_1(\bzero)$ and from the definition of $a(\bx,\bxi)$ we have
\begin{align}
|\ba(\bx,\bxi)|^2 &= |\theta(\bx+\epsilon \bxi, \bu)|^2 \notag \\
&= \bigg\vert \frac{1}{\omega_d} \int_{H_1(\bzero)} \omega(\bx + \epsilon \bxi + \epsilon \boldsymbol{\eta}) J(|\boldsymbol{\eta}|) \bubar_{\boldsymbol{\eta}}(\bx + \epsilon\bxi) \cdot \be_{\boldsymbol{\eta}} d\boldsymbol{\eta} \bigg\vert^2 \notag \\
&\leq \bigg\vert \frac{1}{\omega_d} \int_{H_1(\bzero)}  J(|\boldsymbol{\eta}|) (|\bu(\bx + \epsilon\bxi + \epsilon \boldsymbol{\eta})| +  |\bu(\bx + \epsilon \bxi)|) d\boldsymbol{\eta} \bigg\vert^2,
\end{align}
where we used the fact that $0\leq \omega(x) \leq 1$ and definition of $\bubar_{\boldsymbol{\eta}}(\bx+\epsilon\bxi)$. We now apply inequality \autoref{eq:ineq symm square} with $C = 1$, $\alpha = 0$ and $p(\boldsymbol{\eta}) = |\bu(\bx + \epsilon\bxi + \epsilon \boldsymbol{\eta})| +  |\bu(\bx + \epsilon \bxi)|$ to obtain
\begin{align}
|\ba(\bx,\bxi)|^2 &\leq \frac{\bar{J}_0}{\omega_d} \int_{H_1(\bzero)}  J(|\boldsymbol{\eta}|) (|\bu(\bx + \epsilon\bxi + \epsilon \boldsymbol{\eta})| +  |\bu(\bx + \epsilon \bxi)|)^2 d\boldsymbol{\eta} \notag \\
&\leq \frac{2\bar{J}_0}{\omega_d} \int_{H_1(\bzero)}  J(|\boldsymbol{\eta}|) (|\bu(\bx + \epsilon\bxi + \epsilon \boldsymbol{\eta})|^2 +  |\bu(\bx + \epsilon \bxi)|^2) d\boldsymbol{\eta},
\end{align}
where we have also used the inequality $(a+b)^2 \leq 2 a^2 + 2b^2$. This inequality holds for all $\bx$ and $\bxi$ which includes $\bx = \bx_i$ and $\bxi = \bzero$. 

With estimate on $|\ba(\bx,\bxi)|^2$ and the fact that $\bu$ is a piecewise constant function defined over unit cells $U_i$, we immediately have
\begin{align}
\sum_{i,\bx_i \in D} \int_{U_i} |\ba(\bx,\bxi)|^2 d\bx 
&\leq 4 \bar{J}_0^2 ||\bu||_{L^2}^2 = 4 \bar{J}_0^2 ||\buhat(t)||_{L^2}^2
\end{align}
where we substituted $\buhat(t)$ for $\bu$. 
Combining above estimate with \autoref{eq:stab 5} we get
\begin{align}
I &\leq 64 \bar{J}^2_0||\buhat(t)||_{L^2}^2.
\end{align}
Finally, we use e the bound on $I$ and substitute it into \autoref{eq:stab 4} to show
\begin{align}
&||\sigma_D(t)||^2_{L^2} \leq \left(\frac{g''(0)}{\epsilon^2}  \right)^2 \frac{\bar{J}_0}{\omega_d} \int_{H_1(\bzero)} J(|\bxi|) 64 \bar{J}^2_0||\buhat(t)||_{L^2}^2 d\bxi \notag \\
\Rightarrow &||\sigma_D(t)||_{L^2} \leq \frac{8 g''(0) \bar{J}^2_0}{\epsilon^2} ||\buhat(t)||_{L^2}.
\end{align}

On renaming the constants the bound on $\sigma(t)$ can be summarized as
\begin{align}\label{eq:bound on sigma}
&||\sigma(t)||_{L^2} \notag \\
&\leq \begin{cases}
\frac{4 C^f_1 \bar{J}_{1/2}\sqrt{|D|}}{\epsilon^{3/2}} + \frac{4 C^g_1 \bar{J}_0\sqrt{|D|}}{\epsilon^2} \leq \frac{C}{\epsilon^2} \quad \text{for convex-concave }g,\\
\frac{4 C^f_1 \bar{J}_{1/2}\sqrt{|D|}}{\epsilon^{3/2}} + \frac{8 g''(0) \bar{J}^2_0}{\epsilon^2} ||\buhat(t)||_{L^2} \leq \frac{C_1 + C_2 ||\buhat(t)||_{L^2}}{\epsilon^2} \quad \text{for quadratic }g \,.
\end{cases}
\end{align}

\subsection{Energy inequality}
\label{energy inequality}
From \autoref{eq:stab 2.1} and noting the identity 
\begin{align}
\frac{d}{dt}\mathcal{E}^\epsilon(\buhat)(t) = (\ddot{\buhat}(t), \dot{\buhat}(t)) - (\mathcal{L}^\epsilon(\buhat(t)), \dot{\buhat}(t))
\end{align}
we have
\begin{align}\label{eq:stab 6}
\frac{d\mathcal{E}^\epsilon(\buhat)(t)}{dt} &=  (\hat{\bb}(t), \dot{\buhat}(t)) + (\sigma(t), \dot{\buhat}(t)) \notag \\
&\leq (||\hat{\bb}(t)||_{L^2} + ||\sigma(t)||_{L^2}) \,||\dot{\buhat}(t)||_{L^2}.
\end{align}

When $g$ is convex-concave we can apply identical steps as in the proof of Theorem 5 of \cite{CMPer-JhaLipton} together with the estimate \autoref{eq:bound on sigma} to obtain
\begin{align}
\sqrt{\mathcal{E}^\epsilon(\buhat)(t)} &\leq \sqrt{\mathcal{E}^\epsilon(\buhat)(0)} + \frac{t C}{\epsilon^2} + \int_0^t ||\hat{\bb}(s)||_{L^2} ds
\end{align}
for all $t\in [0,T]$. This completes the proof of energy stability for convex-concave potential functions $g$. 

We now address the case of quadratic potential functions $g$.  We introduce the energy $\bar{\mathcal{E}}^\epsilon(\buhat)(t)$ given by
\begin{align*}
\bar{\mathcal{E}}^\epsilon(\bu)(t) := \mathcal{E}^\epsilon(\bu)(t) + \frac{1}{2}||\bu(t)||^2_{L^2}. 
\end{align*}
Differentiation shows that
\begin{align*}
\frac{d\mathcal{E}^\epsilon(\buhat)(t)}{dt} &= \frac{d\bar{\mathcal{E}}^\epsilon(\buhat)(t)}{dt}  - (\buhat(t), \dot{\buhat}(t)).
\end{align*}
Thus from \autoref{eq:stab 6} we get
\begin{align}\label{eq:stab 7}
\frac{d\bar{\mathcal{E}}^\epsilon(\buhat)(t)}{dt} &\leq (||\hat{\bb}(t)||_{L^2} + ||\sigma(t)||_{L^2}) \,||\dot{\buhat}(t)||_{L^2} +  (\buhat(t), \dot{\buhat}(t)) \notag \\
&\leq (||\hat{\bb}(t)||_{L^2} + C_1/\epsilon^2) \,||\dot{\buhat}(t)||_{L^2} + (C_2/\epsilon^2 + 1)||\buhat(t)||_{L^2} \, \,||\dot{\buhat}(t)||_{L^2}.
\end{align}
From the definition of energy $\bar{\mathcal{E}}^\epsilon$ we have
\begin{align}
||\buhat(t)||_{L^2} &\leq \sqrt{2 \bar{\mathcal{E}}^\epsilon(\buhat)(t)} \hbox{            and},\notag \\
||\dot{\buhat}(t)||_{L^2} &\leq \sqrt{2 \bar{\mathcal{E}}^\epsilon(\buhat)(t)}.
\end{align}
Using the above inequalities in \autoref{eq:stab 7} along with Cauchy's inequality gives
\begin{align}
\frac{d\bar{\mathcal{E}}^\epsilon(\buhat)(t)}{dt} &\leq ||\hat{\bb}(t)||^2_{L^2} + \frac{C_1^2}{\epsilon^4} + 3(\frac{C_2}{\epsilon^2} + 1) \bar{\mathcal{E}}^\epsilon(\buhat)(t).
\end{align}
Using the integrating factor $\exp[-3(C_2/\epsilon^2 + 1)t]$ we recover the inequality
\begin{align}
\bar{\mathcal{E}}^\epsilon(\buhat)(t) &\leq \exp[3(C_2/\epsilon^2 + 1)t] \bigg( \bar{\mathcal{E}}^\epsilon(\buhat)(0) \notag \\
&\quad \quad + \int_0^t (\frac{C_1^2}{\epsilon^4} + ||\hat{\bb}(s)||^2_{L^2}) \exp[-3(C_2/\epsilon^2 + 1)s]  ds \bigg).
\end{align}
This completes the proof of \autoref{thm:stab semi}. 

\section{Conclusions}
\label{s:conclusions}
In this article, we present an a-priori convergence analysis for a class of nonlinear nonlocal  state based peridynamic models. We have shown that the convergence rate applies, even when the fields do not have well-defined spatial derivatives. The results are valid for two different classes of state-based peridynamic models depending on the potential functions associated with the dilatational  energy. For both models the potential function characterizing  the energy due to tensile strain is of convex-concave type while the potential function for the dilatational strain can be either convex-concave or quadratic.  
The convergence rate of the discrete approximation to the true solution in the mean square norm is given by $C(\Delta t+h^\gamma/\epsilon^2)$. Here the constant depends on the H\"older and $L^2$ norm of the true solution and its time derivatives.
The Lipschitz property of the nonlocal, nonlinear force together with boundedness of the nonlocal kernel plays an important role. It ensures that the error in the nonlocal force remains bounded when replacing the exact solution with its approximation. This, in turn, implies that even in the presence of mechanical instabilities the global approximation error remains controlled by the local truncation error in space and time. This is supported by numerical results with crack propagation. The analysis shows that the method is stable and one can control the error by choosing the time step and spatial discretization sufficiently small. 


\newcommand{\noopsort}[1]{}


\begin{thebibliography}{33}
\expandafter\ifx\csname natexlab\endcsname\relax\def\natexlab#1{#1}\fi
\expandafter\ifx\csname url\endcsname\relax
  \def\url#1{\texttt{#1}}\fi
\expandafter\ifx\csname urlprefix\endcsname\relax\def\urlprefix{URL }\fi

\bibitem[{Agwai et~al.(2011)Agwai, Guven, and Madenci}]{CMPer-Agwai}
Agwai, A., Guven, I., Madenci, E., 2011. Predicting crack propagation with
  peridynamics: a comparative study. International journal of fracture 171~(1),
  65--78.

\bibitem[{Aksoylu and Mengesha(2010)}]{AksMen}
Aksoylu, B., Mengesha, T., 2010. Results on nonlocal boundary value problems.
  Numerical functional analysis and optimization 31~(12), 1301--1317.

\bibitem[{Aksoylu and Parks(2011)}]{AksoyluParks}
Aksoylu, B., Parks, M.~L., 2011. Variational theory and domain decomposition
  for nonlocal problems. Applied Mathematics and Computation 217~(14),
  6498--6515.

\bibitem[{Aksoylu and Unlu(2014)}]{AksoyluUnlu}
Aksoylu, B., Unlu, Z., 2014. Conditioning analysis of nonlocal integral
  operators in fractional sobolev spaces. SIAM Journal on Numerical Analysis
  52, 653--677.

\bibitem[{Bobaru and Hu(2012)}]{BobaruHu}
Bobaru, F., Hu, W., 2012. The meaning, selection, and use of the peridynamic
  horizon and its relation to crack branching in brittle materials.
  International journal of fracture 176~(2), 215--222.

\bibitem[{Bobaru et~al.(2009)Bobaru, Yang, Alves, Silling, Askari, and
  Xu}]{CMPer-Bobaru}
Bobaru, F., Yang, M., Alves, L.~F., Silling, S.~A., Askari, E., Xu, J., 2009.
  Convergence, adaptive refinement, and scaling in 1d peridynamics.
  International Journal for Numerical Methods in Engineering 77~(6), 852--877.

\bibitem[{Chen and Gunzburger(2011)}]{CMPer-Chen}
Chen, X., Gunzburger, M., 2011. Continuous and discontinuous finite element
  methods for a peridynamics model of mechanics. Computer Methods in Applied
  Mechanics and Engineering 200~(9), 1237--1250.

\bibitem[{Driver(2003)}]{MA-Driver}
Driver, B.~K., June 2003. Analysis tools with applications. Lecture Notes.
\newline\urlprefix\url{http://math.ucsd.edu/~driver/240-01-02/Lecture_Notes/anal.pdf}

\bibitem[{Du et~al.(2012)Du, Gunzburger, Lehoucq, and Zhou}]{CMPer-Du1}
Du, Q., Gunzburger, M., Lehoucq, R.~B., Zhou, K., 2012. Analysis and
  approximation of nonlocal diffusion problems with volume constraints. SIAM
  review 54~(4), 667--696.

\bibitem[{D’Elia and Gunzburger(2013)}]{DeEllaGunzberger}
D’Elia, M., Gunzburger, M., 2013. The fractional laplacian operator on
  bounded domains as a special case of the nonlocal diffusion operator.
  Computers \& Mathematics with Applications 66~(7), 1245--1260.

\bibitem[{Emmrich et~al.(2013)Emmrich, Lehoucq, and Puhst}]{CMPer-Emmrich}
Emmrich, E., Lehoucq, R.~B., Puhst, D., 2013. Peridynamics: a nonlocal
  continuum theory. In: Meshfree Methods for Partial Differential Equations VI.
  Springer, pp. 45--65.

\bibitem[{Florin et~al.(2016)Florin, Foster, Geubelle, Geubelle, and
  Silling}]{Handbook}
Florin, B., Foster, J.~T., Geubelle, P.~H., Geubelle, P.~H., Silling, S.~A.,
  2016. Handbook of peridynamic modeling.

\bibitem[{Foster et~al.(2011)Foster, Silling, and Chen}]{CMPer-Silling7}
Foster, J.~T., Silling, S.~A., Chen, W., 2011. An energy based failure
  criterion for use with peridynamic states. International Journal for
  Multiscale Computational Engineering 9~(6).

\bibitem[{Ghajari et~al.(2014)Ghajari, Iannucci, and Curtis}]{CMPer-Ghajari}
Ghajari, M., Iannucci, L., Curtis, P., 2014. A peridynamic material model for
  the analysis of dynamic crack propagation in orthotropic media. Computer
  Methods in Applied Mechanics and Engineering 276, 431--452.

\bibitem[{Ha and Bobaru(2010)}]{HaBobaru}
Ha, Y.~D., Bobaru, F., 2010. Studies of dynamic crack propagation and crack
  branching with peridynamics. International Journal of Fracture 162~(1-2),
  229--244.

\bibitem[{Jha and Lipton(2017)}]{CMPer-JhaLipton3}
Jha, P.~K., Lipton, R., October 2017. Finite element approximation of nonlinear
  nonlocal models. arXiv preprint arXiv:1710.07661.

\bibitem[{Jha and Lipton(2018)}]{CMPer-JhaLipton}
Jha, P.~K., Lipton, R., 2018. Numerical analysis of nonlocal fracture models in
  h\"older space. SIAM Journal on Numerical Analysis 56~(2), 906--941.
\newline\urlprefix\url{https://doi.org/10.1137/17M1112236}

\bibitem[{Lindsay et~al.(2016)Lindsay, Parks, and Prakash}]{LindParks}
Lindsay, P., Parks, M., Prakash, A., 2016. Enabling fast, stable and accurate
  peridynamic computations using multi-time-step integration. Computer Methods
  in Applied Mechanics and Engineering 306, 382--405.

\bibitem[{Lipton(2014)}]{CMPer-Lipton3}
Lipton, R., 2014. Dynamic brittle fracture as a small horizon limit of
  peridynamics. Journal of Elasticity 117~(1), 21--50.

\bibitem[{Lipton(2016)}]{CMPer-Lipton}
Lipton, R., 2016. Cohesive dynamics and brittle fracture. Journal of Elasticity
  124~(2), 143--191.

\bibitem[{Lipton et~al.(2018)Lipton, Said, and Jha}]{CMPer-Lipton4}
Lipton, R., Said, E., Jha, P.~K., 2018. Dynamic brittle fracture from nonlocal
  double-well potentials: A state-based model. Handbook of Nonlocal Continuum
  Mechanics for Materials and Structures, 1--27.
\newline\urlprefix\url{https://doi.org/10.1007/978-3-319-22977-5_33-1}

\bibitem[{Lipton et~al.(2016)Lipton, Silling, and Lehoucq}]{CMPer-Lipton2}
Lipton, R., Silling, S., Lehoucq, R., 2016. Complex fracture nucleation and
  evolution with nonlocal elastodynamics. arXiv preprint arXiv:1602.00247.

\bibitem[{Mengesha and Du(2013)}]{CMPer-Du3}
Mengesha, T., Du, Q., 2013. Analysis of a scalar peridynamic model with a sign
  changing kernel. Discrete Contin. Dynam. Systems B 18, 1415--1437.

\bibitem[{Mengesha and Du(2015)}]{CMPer-Mengesha2}
Mengesha, T., Du, Q., 2015. On the variational limit of a class of nonlocal
  functionals related to peridynamics. Nonlinearity 28~(11), 3999.

\bibitem[{Nochetto et~al.(2015)Nochetto, Ot{\'a}rola, and Salgado}]{Nochetto1}
Nochetto, R.~H., Ot{\'a}rola, E., Salgado, A.~J., 2015. A pde approach to
  fractional diffusion in general domains: a priori error analysis. Foundations
  of Computational Mathematics 15~(3), 733--791.

\bibitem[{Silling et~al.(2010)Silling, Weckner, Askari, and
  Bobaru}]{CMPer-Silling5}
Silling, S., Weckner, O., Askari, E., Bobaru, F., 2010. Crack nucleation in a
  peridynamic solid. International Journal of Fracture 162~(1-2), 219--227.

\bibitem[{Silling(2000)}]{CMPer-Silling}
Silling, S.~A., 2000. Reformulation of elasticity theory for discontinuities
  and long-range forces. Journal of the Mechanics and Physics of Solids 48~(1),
  175--209.

\bibitem[{Silling and Askari(2005)}]{CMPer-Silling8}
Silling, S.~A., Askari, E., 2005. A meshfree method based on the peridynamic
  model of solid mechanics. Computers \& structures 83~(17), 1526--1535.

\bibitem[{Silling and Bobaru(2005)}]{SillBob}
Silling, S.~A., Bobaru, F., 2005. Peridynamic modeling of membranes and fibers.
  International Journal of Non-Linear Mechanics 40~(2), 395--409.

\bibitem[{Silling et~al.(2007)Silling, Epton, Weckner, Xu, and Askari}]{States}
Silling, S.~A., Epton, M., Weckner, O., Xu, J., Askari, E., 2007. Peridynamic
  states and constitutive modeling. Journal of Elasticity 88~(2), 151--184.

\bibitem[{Silling and Lehoucq(2008)}]{CMPer-Silling4}
Silling, S.~A., Lehoucq, R.~B., 2008. Convergence of peridynamics to classical
  elasticity theory. Journal of Elasticity 93~(1), 13--37.

\bibitem[{Tian and Du(2013)}]{CMPer-Du2}
Tian, X., Du, Q., 2013. Analysis and comparison of different approximations to
  nonlocal diffusion and linear peridynamic equations. SIAM Journal on
  Numerical Analysis 51~(6), 3458--3482.

\bibitem[{Weckner and Emmrich(2005)}]{CMPer-Weckner}
Weckner, O., Emmrich, E., 2005. Numerical simulation of the dynamics of a
  nonlocal, inhomogeneous, infinite bar. J. Comput. Appl. Mech 6~(2), 311--319.

\end{thebibliography}
\end{document}